\documentclass{amsart}

\usepackage{amssymb,amsmath,amsthm,amsfonts,framed}
\usepackage{framed,comment,enumerate}
\usepackage{xcolor, framed}
\usepackage{tikz}
\usepackage{color, framed}
\usepackage{graphicx,caption,subcaption}

\usepackage[
  hmarginratio={1:1},     % equal left and right margins
  vmarginratio={1:1},     % equal top and bottom margins
  textwidth=15cm,        % new text width  
  textheight=20cm,  
  heightrounded,          % always useful
]{geometry}

\usepackage[colorlinks]{hyperref}
\hypersetup{linkcolor=blue,citecolor=blue,filecolor=black,urlcolor=blue}

\usepackage{verbatim}

\usepackage[sort,nocompress]{cite}

\def\R {\mathbb{R}}
\def\C {\mathcal{C}}
\def\D {\mathcal{D}}
\def\N {\mathbb{N}}

\def\eps{\varepsilon}
\def\dist{{\rm dist}}

\def\d{\delta}
\def\x{\xi}
\def\dn{\delta_{i,n}}

\newcommand{\wc}{\rightharpoonup}

\newcommand{\pa}{\partial}
\newcommand{\bs}[1]{\boldsymbol{#1}}
\newcommand{\mf}[1]{\mathbf{#1}}

\DeclareMathOperator{\supp}{supp}

\newtheorem{proposition}{Proposition}[section]
\newtheorem{theorem}[proposition]{Theorem}
\newtheorem*{theorem*}{Theorem}
\newtheorem{corollary}[proposition]{Corollary}
\newtheorem{lemma}[proposition]{Lemma}

\theoremstyle{definition}
\newtheorem{remark}[proposition]{Remark}
\numberwithin{equation}{section}

\title{On Coron's problem for weakly coupled elliptic systems}

\author{Angela Pistoia and Nicola Soave}

\address{
\hbox{\parbox{5.7in}{\medskip\noindent
Angela Pistoia\\ 
Dipartimento di Metodi e Modelli Matematici, Universit\`a di Roma ``La Sapienza",\\
Via Antonio Scarpa 16, 00161 Roma (Italy),\\[2pt]
{\em{E-mail address: }}{\tt angela.pistoia@uniroma1.it} \\ [5pt]
Nicola Soave\\
Mathematisches Institut, Justus-Liebig-Universit\"at Giessen, \\
Arndtstrasse 2, 35392 Giessen (Germany),\\[2pt]
{\em{E-mail address: }}{\tt nicola.soave@gmail.com, Nicola.Soave@math.uni-giessen.de}}}}


\thanks{{\em Acknowlegments.} Part of this work was carried out while N. Soave was visiting the University of Rome, ``La Sapienza", which he wishes to thank for the hospitality. N. Soave is partially supported through the project ERC Advanced Grant  2013 n. 339958 ``Complex Patterns for Strongly Interacting Dynamical Systems - COMPAT''}

\begin{document}
\maketitle

\begin{abstract}
We consider the following critical weakly coupled elliptic system 
\[
\begin{cases}
-\Delta u_i = \mu_i |u_i|^{2^*-2}u_i + \sum_{j \neq i} \beta_{ij} |u_j|^{\frac{2^*}{2}} |u_i|^{\frac{2^*-4}{2}} u_i & \text{in $\Omega_\eps$} \\
%u_i >0 & \text{in $\Omega_\eps$} \\
u_i = 0 & \text{on $\pa \Omega_\eps$},
\end{cases} \qquad i =1,\dots,m,
\]
in a domain $\Omega_\eps \subset \R^N$, $N=3,4$, with small shrinking holes as the parameter $\eps \to 0$. We prove the existence of positive solutions of two different types: either each density concentrates around a different hole, or we have groups of components such that all the components within a single group concentrate around the same point, and different groups concentrate around different points.
\end{abstract}

\section{Introduction}

The system of nonlinear elliptic equations
\begin{equation}\label{BE sys}
\begin{cases}
-\Delta u + \lambda_i u_i = \mu_i |u_i|^{p-1}u_i + \sum_{j \neq i} \beta_{ij} |u_j|^{\frac{p+1}{2}} |u_i|^{\frac{p-3}{2}} u_i   & \text{in $\Omega$}\\
u_i \in H_0^1(\Omega),
\end{cases} \qquad i =1,\dots,m,
\end{equation}
where $\Omega$ is either a bounded domain or the whole space $\R^N$, has been studied intensively in the last decade, in light of its relevance in different physical context: \eqref{BE sys} appears when looking for solitary wave solutions $\Phi_i(t,x) = e^{ \iota \lambda_i t} u_i(x)$ for the coupled Gross-Pitaevskii equation
\[
- \iota \pa_t \Phi_i = \Delta \Phi_i + \mu_i |\Phi_i|^{p-1} \Phi_i+ \sum_{j \neq i}\beta_{ij}
 |\Phi_i|^{\frac{p-3}2} |\Phi_j|^{\frac{p+1}2} \Phi_i,
 \]
which is of interest in nonlinear optics and in quantum mechanics, see e.g. \cite{AkAn, Timm}. In the models, $|u_i|$ represents the amplitude of the $i$-th density, and the real parameters $\mu_i$ and $\beta_{ij}$ represent the intra-spaces and inter-species scattering length, describing respectively the interaction between particles of the same component or of different components. In particular, the positive sign of $\mu_i$ (and of $\beta_{ij}$) stays for attractive interaction, while the negative sign stays for repulsive interaction.

From the mathematical point of view, \eqref{BE sys} is one of the simplest, yet highly non-trivial, examples of \emph{weakly coupled system}, i.e. is a system admitting non-trivial solutions $(u_1,\dots,u_m) \not \equiv (0,\dots,0)$ with some trivial components $u_i \equiv 0$. This feature stimulated a lot of research about existence of \emph{fully nontrivial solutions}, i.e. solutions with $u_i \not \equiv 0$ for every $i$. Nowadays, many results in this direction are available, mainly concerning the cubic problem $p=3$ in dimension $N \le 3$, i.e. in a \emph{Sobolev subcritical regime}. A complete review of the result in this framework would be beyond the aim of the present paper, and we refer the interested read to the quite exhaustive  introductions in \cite{So, SoTa} and to the references therein.  In this paper we focus instead on the much less understood \emph{Sobolev critical regime} $p=2^*-1$ in dimension $N=3$ or $N=4$, where $2^*= 2N/(N-2)$ is the critical exponent for the Sobolev embedding $H^1(\R^N) \hookrightarrow L^r(\R^N)$; moreover, from now on we limit ourselves to the \emph{focusing setting} $\mu_i >0$ for every $i$.

The study of the critical system \eqref{BE sys} started with the Chen and Zou's paper \cite{ChenZou1}, where the authors focused on \eqref{BE sys} with $2$ components in bounded domains of $\R^4$ (thus, with $p=3$), and proved existence of least energy positive solutions under suitable assumptions on the parameters $\lambda_i$, $\mu_i>0$, $\beta_{ij}$. In \cite{ChenZou2}, the authors extended their results in higher dimension $N \ge 5$. In both papers it is assumed that $-\nu_1(\Omega)< \lambda_1,\lambda_2<0$ (here $\nu_1(\Omega)$ denotes the first eigenvalue of $(-\Delta)$ with homogeneous Dirichlet boundary conditions in $\Omega$), and this plays a crucial role: indeed, as remarked by Chen and Zou, system \eqref{BE sys} with $\Omega$ bounded, $\mu_i>0$ and $p = 2^*-1$ can be considered a critically coupled version of the \emph{Brezis-Nirenberg problem}
\begin{equation}\label{bre-nir}
\begin{cases}
-\Delta u +  \lambda u=  |u|^{2^*-2} u  & \text{in $\Omega$} \\
u = 0 & \text{on $\pa \Omega$},
\end{cases}
\end{equation}
and it is well known that in any dimension $N \ge 4$ \eqref{bre-nir} admits a positive solution (for an arbitrary bounded domain) if and only if $-\nu_1(\Omega)<\lambda<0$ (see \cite{BreNir} for this result, and the survey \cite{Pis} for a more extended discussion). 

The relation between \eqref{BE sys} and the Brezis-Nirenberg problem has been recently exploited also in \cite{ChenLin1, PisTav}. In \cite{ChenLin1}, Chen and Lin described the blow-up behaviur of least energy positive solutions as $\lambda_i \to 0$, in case of $2$ components system with $\beta_{12}>0$. In \cite{PisTav}, the first author and Tavares constructed, under appropriate assumptions on the domain $\Omega \subset \R^4$ and on the parameters $\beta_{ij} \in \R$, solutions to \eqref{BE sys} with all the components $u_i$ concentrating around different points $a_i \in \Omega$ as $\lambda_i \to 0^-$. 

We finally refer to \cite{GuLiWe}, where the authors proved existence of infinitely many non-radial solutions for \eqref{BE sys} in $\R^3$, with $\lambda_i=0$ and $\beta_{ij}<0$. 
%and to \cite{ChenLi, ChenLin2, GuoLiu}, where the authors studied

In this paper we address \eqref{BE sys} when $\Omega$ is a bounded domain of $\R^3$ or $\R^4$ and $\lambda_i = 0$ for every $i$. If $\Omega$ is star-shaped, the Pohozaev identity for gradient-type systems implies that the problem has no nonnegative solutions (meaning that $u_i \ge 0$ for every $i$)  but the trivial one $(u_1,\dots,u_m) = (0,\dots,0)$, but if $\Omega$ has some hole there is hope to find fully nontrivial positive solutions, in the spirit of the celebrated \emph{Coron and Bahri-Coron results}, which we briefly review: let us consider the critical problem
\begin{equation}\label{critical eq}
\begin{cases}
-\Delta u = |u|^{2^*-2}u & \text{in $\Omega$} \\
u >0 & \text{in $\Omega$} \\
u = 0 & \text{on $\pa \Omega$},
\end{cases} \quad \text{with $\Omega$ bounded domain}.
\end{equation}
If $\Omega$ is star-shaped, non-trivial solutions do not exist, but the situation drastically changes removing this geometric assumption: indeed, as observed by Kazdan and Warner in \cite{KaWa}, \eqref{critical eq} in an annulus admits a positive solution. The result was then improved by Coron, who showed in \cite{Co} that \eqref{critical eq} has a positive solution as long as $\Omega$ has a small hole. A further improvement was achieved by Bahri and Coron, who proved in \cite{BaCo} that a positive solution does exist provided that $\Omega$ has non-trivial topology. For multiplicity results and sign changing solutions, we refer the interested reader to \cite{ClMuPi, ClWe, GeMuPi, LiYaYa, MuPi, Rey}.

The purpose of this paper is to discuss the extension to \eqref{BE sys} of the Coron result. While in \cite{Co} a variational argument is considered, we adopt here a perturbation approach based on the Lyapunov-Schmidt finite dimensional reduction, which has been already fruitfully used to deal with scalar Coron's type problem in \cite{ClMuPi, GeMuPi, MuPiPisa, MuPi}. In order to state our main result, we introduce some notation and recall some basic results.

We consider from now on the following general assumptions, which will not be always recalled. Let $\Omega \subset \R^N$ with $N=3,4$ be a bounded and sufficiently regular domain, and let $m \in \N$. Let $a_1,\dots,a_m \in \Omega$ be $m$ (not necessarily different) points in $\Omega$, $r_1,\dots,r_m>0$, $\mu_1,\dots,\mu_m>0$ and $\beta_{ij}=\beta_{ji}$ for every $i,j =1,\dots,m$, $i \neq j$. We consider the following Coron-type problem:
\begin{equation}\label{system}
\begin{cases}
-\Delta u_i = \mu_i |u_i|^{p-1}u_i + \sum_{j \neq i} \beta_{ij} |u_j|^{\frac{p+1}{2}} |u_i|^{\frac{p-3}{2}} u_i & \text{in $\Omega_\eps$} \\
%u_i >0 & \text{in $\Omega_\eps$} \\
u_i = 0 & \text{on $\pa \Omega_\eps$},
\end{cases} \qquad i =1,\dots,m,
\end{equation}
where
\[
\Omega_\eps:= \Omega \setminus \bigcup_{i=1}^m B_{r_i\eps}(a_i),
\]
and $p:= 2^*-1= (N+2)/(N-2)$. 

For $\delta>0$ and $\xi \in \R^N$, we denote by $U_{\delta,\xi}$ the standard bubble
\begin{equation}\label{def bubble}
U_{\d,\x}(x):= \alpha_N \left( \frac{\delta}{\delta^2 + |x-\xi|^2} \right)^{\frac{N-2}{2}},
\end{equation}
where $\alpha_N>0$ is a suitable constant depending on the space dimension. It is well known (see \cite[Corollary 8.2]{CaGiSp}) that the family $\{U_{\d,\x}: \ \d >0,\ \x \in \R^N\}$ contains all the solutions to the critical problem
\begin{equation}\label{critical single}
\begin{cases}
-\Delta U = U^p & \text{in $\R^N$} \\
U >0 & \text{in $\R^N$} \\
U \in \mathcal{D}^{1,2}(\R^N),
\end{cases} 
\end{equation}
where $\mathcal{D}^{1,2}(\R^N)$ is the completion of $\mathcal{C}^\infty_c(\R^N)$ with respect to the norm $\|u\|_{\mathcal{D}^{1,2}(\R^N)}:= \|\nabla u\|_{L^2(\R^N)}$. We consider the projection $P_\eps U_{\d,\xi}$ of $U_{\d,\xi}$ into $H_0^1(\Omega_\eps)$, i.e. the only solution to 
\begin{equation}\label{def projection PU}
\begin{cases}
-\Delta (P_\eps U_{\d,\xi}) = -\Delta U_{\d,\x} = U_{\d,\x}^p & \text{in $\Omega_\eps$} \\
P_\eps U_{\d,\x}  = 0 & \text{on $\pa \Omega_\eps$}.
\end{cases} 
\end{equation}

%The first of our main result allows interaction of competitive type ($\beta_{ij}<0$) or of weakly cooperative type ($0 \le \beta_{ij}$ small) between the different components.

The first of our main results describe the situation where the components $u_i$ are all concentrating around different points.
 
\begin{theorem}\label{thm: main 1}
Let $N=4$, and let us suppose that $a_i \neq a_j$ for $i \neq j$. Then there exists $\bar \eps, \bar \beta>0$ such that, if $\eps \in (0,\bar \eps)$ and $-\infty<\beta_{ij}< \bar \beta$ for every $i \neq j$, then problem \eqref{system} has a fully nontrivial solution $(u_{1,\eps},\dots,u_{m,\eps})$, where each $u_{i,\eps}$ is positive and is concentrating around $a_i$ as $\eps \to 0$. 

To be precise, we have that
\[
u_{i,\eps} = \mu_i^{-\frac1{p-1}} P_\eps U_{\d_i,\x_i} + \phi_{i,\eps}>0 \qquad \text{in $\Omega_\eps$},
\]
where $\d_i = d_i \sqrt{\eps}$, $\xi_i = a_i + \d_i \tau_i$ for suitably chosen $d_i>0$ and $\tau_i \in \R^N$, and there exists $C>0$ (independent of $\eps$) such that
\[
\|\phi_{i,\eps}\|_{H_0^1(\Omega_\eps)} \le C \eps^{\frac{N-2}2} \qquad i=1,\dots,m.
\]

If $N=3$, the same conclusion holds without any restriction on $\beta_{ij}$ (which can be also positive and large).
\end{theorem}

Notice that in dimension $N=4$ we allow interactions of competitive type ($\beta_{ij}<0$) or of weakly cooperative type ($0 \le \beta_{ij}$ small) between the different components, while in dimension $N=3$ we have no restriction, and for the reason of this difference we refer to the forthcoming Remark \ref{rem: N=4,3}

In Theorem \ref{thm: main 1} we proved the existence of solutions to \eqref{system} with all the components concentrating around different points. On the other hand, it is natural to wonder if it is possible to find solutions with several groups $G_1,\dots,G_q$ of components such that:
\begin{itemize}
\item each component within a given group $G_h$ concentrate around a point $a_h$;
\item the different groups concentrate around different points, i.e. $a_h \neq a_k$ if $h \neq k$.
\end{itemize}
The following theorem gives a positive answer to this question. 

In what follow we focus on system \eqref{system} with $3$ components in dimension $N=4$.

\begin{theorem}\label{thm: main 2}
Let $N=4$, $m=3$, and let us suppose that $a_1=a_2 \neq a_3$. Let us suppose that either $-\sqrt{\mu_1 \mu_2} < \beta_{12} < \min\{\mu_1,\mu_2\}$, or $\beta_{12}> \max\{\mu_1,\mu_2\}$. Then there exists $\bar \eps, \bar \beta>0$ such that, if $\eps \in (0,\bar \eps)$ and $-\infty<\beta_{13},\beta_{23}< \bar \beta$, then problem \eqref{system} has a solution $(u_{1,\eps},u_{2,\eps},u_{3,\eps})$ where $u_{1,\eps}$ and $u_{2,\eps}$ are concentrating around $a_1$, while $u_{3,\eps}$ is concentrating around $a_3$ as $\eps \to 0$. 

To be precise, let $c_1,c_2,c_3>0$ be defined by
\[
c_1^2 = \frac{\beta_{12}-\mu_2}{\beta_{12} -\mu_1 \mu_2}, \quad c_2^2 = \frac{\beta_{12}-\mu_1}{\beta_{12} -\mu_1 \mu_2}, \quad c_3 = \mu_3^{-\frac1{2}}.
\] 
We have that
\[
u_{i,\eps} = c_i P_\eps U_{\d_1,\x_1} + \phi_{i,\eps}, \quad i=1,2, \quad u_{3,\eps} = c_3  P_\eps U_{\d_3,\x_3} + \phi_{3,\eps}, 
\]
where $\d_i = d_i \sqrt{\eps}$, $\xi_i = a_i + \d_i \tau_i$ for suitably chosen $d_i>0$ and $\tau_i \in \R^N$, and there exists $C>0$ (independent of $\eps$) such that
\[
\|\phi_{i,\eps}\|_{H_0^1(\Omega_\eps)} \le C \eps^{\frac{N-2}2} \qquad i=1,\dots,m.
\]

Moreover, there exists $\tilde \beta \le \bar \beta$ such that if $-\infty<\beta_{ij}< \tilde \beta$ for every $i \neq j$, then $u_{i,\eps}>0$ in $\Omega_\eps$ for every $i$.
\end{theorem}

Notice that the interaction between the components $u_1$ and $u_3$ and between $u_2$ and $u_3$ is competitive or weakly cooperative, while the one between $u_1$ and $u_2$ is weakly competitive or cooperative (possibly with a large coupling parameter $\beta_{12}>0$).

Theorem \ref{thm: main 2} is the first result in the literature dealing with concentration of groups of components. As it will be clear from the proof, our method works in a much more general context with respect to the one considered in Theorem \ref{thm: main 2}. One could both extend the result in $\R^3$, and (what is more important) consider more groups of components. This is the content of the forthcoming Theorem \ref{thm: main 2 general}, for which we need some further notation. Let us consider system \eqref{system}, and let $\beta_{ii}:=\mu_i$.

For an arbitrary $1 < q < m$, we say that a vector $\mf{l}=(l_0,\dots,l_q) \in \N^{q+1}$ is \emph{a $q$-decomposition of $m$} if
\[
0=l_0<l_1<\dots<l_{q-1}<l_q=m;
\]
given a $q$-decomposition $\mf{l}$ of $m$, we set, for $h=1,\dots,q$,
\begin{equation}\label{def Ih}
\begin{split}
& I_h:= \{i \in  \{1,\dots,m\}:  l_{h-1} < i \le l_h \}, \\
& \mathcal{K} := \left\{(i,j) \in I_h \times I_k \text{ with $h \neq k$} \right\}. 
\end{split}
\end{equation}
This way, we have partitioned the set $\{1,\ldots, m\}$ into $q$ groups $I_1,\ldots, I_q$, and have consequently splitted the components into $q$ groups: $\{u_i:\ i\in I_h\}$. 

Let us consider the $q$ non-linear systems
\begin{equation}\label{sist1intro}
-\Delta U_i=\sum\limits_{j \in I_h} \beta_{ij} |U_i|^\frac{p-3}2 |U_j|^\frac{p+1}2 U_i \qquad \hbox{in}\ \mathbb R^N,  \ i \in I_h, \ h=1,\dots,q.
\end{equation}
Any such system has a solution $U_i=c_i U,$ with $c_i>0$ for $i \in I_h$, and where $U=U_{1,0}$ is a standard bubble in $\R^N$ centered in $0$ and with $\d=1$, if and only if there exists a vector of positive numbers $(c_1,\dots,c_k)$ such that
\begin{equation}\label{sist2intro}
 \sum\limits_{j=1,\dots,k}\beta_{ij} c_i^\frac{p-1}2 c_j^\frac{p+1}2=c_i, \quad  i=1,\dots,k.
\end{equation}
Supposing that such a vector does exist, we linearize the $q$ systems \eqref{sist1intro} in $\mathcal{D}^{1,2}(\R^N)$ around the solutions $(c_{l_{h-1}+1}U,\dots,c_{l_h} U)$, obtaining
\begin{equation}\label{sist3intro}
-\Delta v_i=\left[\left(p\beta_{ii}c_i^{p-1}+\frac{p-1}2\sum\limits_{j \in I_h \atop j\not=i}\beta_{ij} c_i^\frac{p-3}2 c_j^\frac{p+1}2 \right) v_i+\frac{p+1}2\sum\limits_{j \in I_h \atop j\not=i}\beta_{ij}c_i^{\frac{p-1}2} c_j^\frac{p-1}2 v_j\right]U^{p-1}\ \quad \text{in }\mathbb R^N, 
\end{equation}
for all $i \in I_h$, for all $h =1,\dots,q$, with $v_i \in \mathcal{D}^{1,2}(\R^N)$.

\begin{theorem}\label{thm: main 2 general}
In the previous setting, let $N=4$, and let $a^1,\dots,a^q \in \Omega$ with $a^h \neq a^k$ for $h \neq k$. Let us suppose that the $q$ systems \eqref{sist1intro} have solutions $(c_{l_{h-1}+1},\dots,c_{l_h})$, $h=1,\dots,q$, with $c_i>0$ for every $i$. Let us suppose also that each of the $q$ linearized systems \eqref{sist3intro} has a $(N+1)$-dimensional set of solutions. Then there exists $\bar \eps, \bar \beta>0$ such that, if $\eps \in (0,\bar \eps)$ and $-\infty<\beta_{ij}< \bar \beta$ for every $(i,j) \in \mathcal{K}$, then problem \eqref{system} has a solution $(u_{1,\eps},\dots,u_{m,\eps})$, where $u_{1,\eps},\dots, u_{l_1,\eps}$ are concentrating around $a^1$, $u_{l_1+1,\eps},\dots, u_{l_2,\eps}$ are concentrating around $a^2$, \dots, and $u_{l_{q-1}+1,\eps},\dots, u_{l_q,\eps}$ are concentrating around $a^q$ as $\eps \to 0$.

To be precise, we have that for any $h=1,\dots,q$
\[
u_{l_{h-1}+1,\eps} = c_{l_{h-1}+1} P_\eps U_{\d_h,\x_h} + \phi_{l_{h-1}+1,\eps}, \quad  \dots \quad  u_{l_h,\eps} = c_{l_h} P_\eps U_{\d_h,\x_h} + \phi_{l_h,\eps}, 
\]
where $\d_h = d_h \sqrt{\eps}$, $\xi_h = a^h + \d_h \tau_h$ for suitably chosen $d_h>0$ and $\tau_h \in \R^N$, and there exists $C>0$ (independent of $\eps$) such that
\[
\|\phi_{i,\eps}\|_{H_0^1(\Omega_\eps)} \le C \eps^{\frac{N-2}2} \qquad i=1,\dots,m.
\]
Moreover, there exists $\tilde \beta \le \bar \beta$ such that if $-\infty<\beta_{ij}< \tilde \beta$ for every $i \neq j$, then $u_{i,\eps}>0$ in $\Omega_\eps$ for every $i$.

In dimension $N=3$, the existence result holds, without any additional assumptions on the parameters $\beta_{ij}$ with $(i,j) \in \mathcal{K}$. Also, in this case we always obtain positive solutions.
\end{theorem}

Theorem \ref{thm: main 2} can be obtained applying Theorem \ref{thm: main 2 general} for the $3$ components system in $\R^4$, considering the $2$-decomposition of $3$ given by $\mf{l}=(0,2,3)$. In this perspective, one has only to verify that the assumptions of Theorem \ref{thm: main 2 general} are satisfied, i.e. that the system
\[
\mu_1 c_1^2 + \beta_{12} c_2^2 = 1, \quad \mu_2 c_2^2 + \beta_{12} c_1^2= 1,
\]
has a solution with $c_1,c_2>0$, and that the linearized problem
\[
-\Delta v_i = (3 \mu_i c_i^2 + \beta_{ij} c_j^2) U^2 v_i + 2\beta_{ij} c_i c_j  U^2 v_j \quad \text{in $\R^4$}, i=1,2, \ j \neq i
\]
has a $5$-dimensional set of solutions.

Nevertheless, for the sake of simplicity and brevity, we directly write down the proof of Theorem \ref{thm: main 2}, without deriving it as corollary of Theorem \ref{thm: main 2}, and we omit the proof of Theorem \ref{thm: main 2 general}. The passage from the particular situation described in Theorem \ref{thm: main 2} to the general setting considered in Theorem \ref{thm: main 2 general} creates additional difficulties only from the technical and notational points of view, and once that the main idea are understood, the interested reader can easily fill the details.

Rather, one point which we want to stress is that the assumptions of Theorem \ref{thm: main 2 general} only regard the data of the problem, and in general are not too difficult to check, as shown in the following proposition. 

Given the number of components $m$, the number of groups $1<q<m$, a $q$-decomposition $\mf{l} = (l_0,\dots,l_q)$ of $m$, and the coupling parameters $\beta_{ij}$, $i,j=1,\dots,m$, let us consider the $q$ matrices
\[
\mathcal{B}_h := \left(\beta_{ij}\right)_{(i,j) \in I_h \times I_h} \quad h =1,\dots,q.
\]
Notice that if $l_{h}-l_{h-1} = 1$ for some $i$, then $I_h$ is a singleton and the corresponding $\mathcal{B}_h$ is simply given by the real number $\beta_{ii}$, $i \in I_h$. 

Conditions on $\mathcal{B}_{h}$ ensuring the existence of  solutions $(c_{l_{h-1}+1},\dots,c_{l_h})$, $h=1,\dots,q$, for the algebraic problem \eqref{sist2intro} are given in \cite[Section 2]{Ba} or \cite[Section 4]{Sir}. 

So, let us suppose that such $(c_{l_{h-1}+1},\dots,c_{l_h})$, $h=1,\dots,q$, do exist. 

\begin{proposition}\label{prop: on hp 2}
Let $N=4$. In the previous setting, and using the notation of Theorem \ref{thm: main 2 general}, let us suppose that for any $h=1,\dots,q$ we have
\[
\text{the matrix $\mathcal{B}_h$ is invertible and has only positive elements}.
\]
Then each of the $h$ linearized systems \eqref{sist3intro} has a $(N+1)$-dimensional set of solutions.
\end{proposition}

Combining then the results in \cite{Ba,Sir} with Proposition \ref{prop: on hp 2}, Theorem \ref{thm: main 2 general} permits to prove existence of a bunch of solutions, with arbitrary number $q$ of groups (each of them having an arbitrary number of components) concentrating around $q$ different points.

\medskip

The proofs of Theorems \ref{thm: main 1} and \ref{thm: main 2} rest upon a finite-dimensional reduction, and will be the object of the next sections. Before proceeding, we conclude the introduction with some comments.

\begin{remark}
We point out that, both in Theorems \ref{thm: main 1}, \ref{thm: main 2} and \ref{thm: main 2 general}, we can deal with systems with mixed cooperation and competition (that is, we can deal with systems where some $\beta_{ij}$ is positive, and some other is negative). This is particularly interesting since the mixed coupling critical case was completely open, and also for subcritical problems has been investigated only in few recent contributions (see \cite{BySaWa, SaWa1, SaWa2, So, SoTa}).
\end{remark}

\begin{remark}
Regarding the positivity in dimension $N=4$, it is natural to think that in Theorems \ref{thm: main 2} and \ref{thm: main 2 general} we have a positive solution without any additional assumptions on $\beta_{ij}$. Indeed, any component $u_{i,\eps}$ is a superposition of positive function and small perturbation term. If the positive part of $\beta_{ij}$ is small for every $i \neq j$, then a short rigorous proof of the positivity can be given arguing as in \cite{PisTav}. If on the other hand some $\beta_{ij}$ is allowed to be large, such proof does not work and one is forced to approach the problem with finer (and much longer) techniques, such as careful $L^\infty$-estimates on the error $\phi_{i,\eps}$ (see for instance \cite[Section 8]{IacVai} for related computations). We decided to not insist on this point for the sake of brevity.
% On the other hand, if some $\beta_{ij}$ are allowed to be large, the positivity of the solution does not follow easily by an application of the maximum principle.
%
%
% becomes tricky and lengthy. 
%we believe that, in , all the time we have existence, we immediately have positivity by the particular shape of the solution (which is a superposition of a positive function and a small perturbation term), without any additional assumptions on $\beta_{ij}$. 
%We decided to assume the condition $-\infty < \beta_{ij} < \tilde \beta$ since, firstly, it does not represent a real restriction in Theorem \ref{thm: main 1} (where we already supposed that coupling parameters $\beta_{ij}$ are either negative, or small), and, secondly, it permits to prove the fact that $u_{i,\eps}>0$ in $\Omega_\eps$ for every $i$ briefly. 
%On the other hand, 
%We do not want to insist too much on this point for the sake of brevity, and we only mention that it should be possible to derive some uniform $L^\infty$ bounds for the remainder terms $\phi_{i,\eps}$, finally leading to the positivity of each components (see for instance \cite[Section 8]{IacVai} for related computations). This kind of arguments has the disadvantage of being quite long.
\end{remark}

\begin{remark}\label{rem: N=4,3}
As it emerges from Theorem \ref{thm: main 1}, in $\R^4$ we have to suppose extra-conditions on $\beta_{ij}$ in order to have existence of fully nontrivial positive solutions, with respect to the $3$-dimensional problem. This is somehow natural, and is related to the particular shape of the cubic system (in dimension $N=4$ we have $p=3$) which makes possible to prove non-existence results (for instance, one can repeat word by word the proof of Theorem 3-($iii$) in \cite{Sir} to rule out the existence of positive solutions under some assumptions on the parameters). This difference between the dimension $N=4$ and $N \neq 4$ was already observed in \cite{ChenZou2}.

The extra-condition $-\infty< \beta_{ij} < \bar \beta$ with $\bar \beta$ small enough was already considered in \cite{PisTav}, and plays the same role in our proof and the one in \cite{PisTav}; it enters in the analysis of the linearization of \eqref{system} and in the proof of the positivity of the solutions. It is interesting that this assumption is not needed in dimension $3$. 

Regarding the higher dimensional case $N \ge 5$, the main obstruction for our method is represented by the lack of regularity of the interaction term $|u_j|^{\frac{p+1}{2}} |u_i|^{\frac{p-3}{2}} u_i$, which is not of class $\C^1$. This creates several additional difficulties in dealing with a linearization of \eqref{system}, difficulties which were already observed in \cite{GuLiWe, PisTav}. 
\end{remark}

\begin{remark}
Both Theorems \ref{thm: main 1} and \ref{thm: main 2} regard the case when at least two components concentrate around different points, and hence it is natural to wonder what happens if we search for solutions to \eqref{system} with all the components concentrating around the same point. Under appropriate assumptions on $\beta_{ij}$, this case is actually much simpler to deal with, since one can obtain a solution with all the components proportional among each other, reducing system \eqref{system} to the scalar problem \eqref{critical single}. For instance, for a system with $2$ components in a domain $\Omega \subset \R^4$, it is not difficult to check that if \eqref{critical single} has a solution $w$ (positive or sign-changing, concentrating around one or more points), and if 
\begin{equation}\label{BaCor}
\text{either $-\sqrt{\mu_1 \mu_2 } <\beta_{12} < \min \{\mu_1,\mu_2\}$, or $\beta_{12}>\max\{\mu_1,\mu_2\}$},
\end{equation}
then \eqref{system} in $\Omega$ has a solution of type $(u_1,u_2) = (c_1 w, c_2 w)$, with $c_1,c_2$ as in Theorem \ref{thm: main 2}. 
In particular, we have:
\begin{itemize}
\item \emph{Bahri-Coron-type result:} if $\Omega \subset \R^4$ has non-trivial topology and \eqref{BaCor} is in force, then the critical system \eqref{system} has a positive solutions.
\item \emph{Solutions concentrating around the same point in domains with one shrinking hole:} if $a \in \Omega \subset \R^4$, $\Omega_\eps= \Omega \setminus B_\eps(a)$, and \eqref{BaCor} holds true, then \eqref{system} has a family of positive solutions with $u_{1,\eps}$ and $u_{2,\eps}$ concentrating around the same point with the same speed.
\item \emph{Multipeak solutions and multiplicity results}, extending \cite{ClMuPi, ClWe, GeMuPi, LiYaYa, MuPi, Rey} as in the previous points.
\end{itemize}
Similar results in higher or lower dimension can be obtained using the existence of $(k_0,l_0)$ as in \cite[Theorem 1.1]{ChenZou2}. In the same way, in order to deal with more than $2$ components one can use the results in \cite[Section 2]{Ba} or \cite[Section 4]{Sir}.
%These considerations motivated us to focus on the case when at least two components concentrate around different points, which cannot be treat using the results for the scalar equation.
\end{remark}

\begin{remark}
In our results, we focus on domains $\Omega_\eps$ with radially symmetric holes $B_{r_i \eps}(a_i)$. It is not difficult to check that we can treat also the case of non-symmetric holes. In such a situation, we could repeat our argument essentially word by word, simply replacing Lemma \ref{lem a.1} with Lemma 1.1 in \cite{MuPiErr}; this introduces some technical complications which we preferred to avoid.
\end{remark}

\subsection*{Structure of the paper}
The proof of Theorem \ref{thm: main 1} is contained in Sections \ref{sec: scheme}-\ref{sec: reduced problem}, while the proof of Theorem \ref{thm: main 2} is given in \ref{sec: thm 2}. Both are based upon the Lyapunov-Schmidt finite dimensional reduction method, which in the context of systems was already adopted in \cite{PisTav} to deal with the Brezis-Nirenberg-type problem. The proof of Proposition \ref{prop: on hp 2} is given in Section \ref{sec: prop on hp 2}.

In Section \ref{sec: scheme} we set up the reduction scheme, splitting system \eqref{system} in two new systems of $m$ equations, one of them living in finite dimension. The infinite dimensional problem is then treated in Section \ref{sec: eq ort} via a fixed point argument, while the finitely-reduced problem is the object of Section \ref{sec: reduced problem}. 

Since Theorem \ref{thm: main 2} shares the same structure and many passages with that of Theorem \ref{thm: main 1}, we put particular emphasis on the main differences. 

As usual when dealing with a perturbation approach, many proofs contains very long computations. In order to keep the presentation as smooth as possible, we collect them in several appendixes, so that the reader can easily understand the main strategy behind each proof, and check the details in a second time.

\subsection*{Notation and preliminary results}
We recall that with $U_{\d,\x}$ and $P_\eps U_{\d,\x}$ we denote the standard bubble and its projection in $H_0^1(\Omega_\eps)$, defined in \eqref{def bubble} and \eqref{def projection PU}. We shall use many times the fact that $0 \le P_\eps U_{\d,\x} \le U_{\d,\x}$, which is a simple consequence of the maximum principle. 

Since $\{U_{\d,\x}\}$ is the set of all the solutions to \eqref{critical single}, it is easy to check that any solution to 
\[
\begin{cases}
-\Delta U = \mu U^p & \text{in $\R^N$} \\
U >0 & \text{in $\R^N$} \\
U \in \mathcal{D}^{1,2}(\R^N)
\end{cases} 
\] 
with $\mu>0$ is given by $\mu^{-\frac{1}{p-1}} U_{\d,\x}$, for some $\d >0$ and $\x \in \R^N$.

Coming back to problem \eqref{critical single}, we recall some properties of the linearized equation
\begin{equation}\label{linear eq}
-\Delta \phi = p U_{\d,\x}^{p-1} \phi \quad \text{in $\R^N$}, \quad \phi \in \mathcal{D}^{1,2}(\R^N).
\end{equation}
It is clear that there are $N+1$ solutions given by 
\begin{equation}\label{def psi}
\begin{split}
\psi_{\d,\x}^0 &:= \frac{\pa U_{\d,\x}}{\pa \d} = \alpha_N \left(\frac{N-2}{2}\right) \delta^{\frac{N-4}{2}} \frac{|x-\xi|^2 -\delta^2}{\left( \d^2 + |x-\x|^2 \right)^{\frac{N}{2}}} \\
\psi_{\d,\x}^h&:=   \frac{\pa U_{\d,\x}}{\pa \x_\ell} = \alpha_N (N-2) \d^{\frac{N-2}{2}} \frac{x_\ell-\xi_\ell}{\left( \d^2 + |x-\x|^2 \right)^{\frac{N}{2}}}, \quad h=1,\dots,N.
\end{split}
\end{equation}
By \cite[Lemma A.1]{BiaEgn}, these functions span the set of solutions to \eqref{linear eq}.

We consider the projections $P_\eps \psi_{\d,\xi}^h$ of $\psi_{\d,\xi}^h$ ($h=0,\dots,N$) into $H_0^1(\Omega_\eps)$, i.e.
\[
\begin{cases}
-\Delta (P_\eps \psi_{\d,\xi}^h) = -\Delta \psi_{\d,\x}^h = p U_{\d,\x}^{p-1} \psi_{\d,\x}^h & \text{in $\Omega_\eps$} \\
P_\eps \psi_{\d,\x}^h  = 0 & \text{on $\pa \Omega_\eps$}.
\end{cases} 
\] 
Notice that the derivative commutes with the projection $P_\eps$, in the sense that $P_\eps \psi_{\d,\xi}^h = \pa_h (P_\eps U_{\d,\xi})$, where $\pa_h$ denotes the partial derivative with respect to $\xi_h$ if $h=1,\dots,N$, and $\pa_0$ denotes the partial derivative with respect to $\d$.

We denote by $G(x,y)$ the Green function of $-\Delta$ with Dirichlet boundary condition in $\Omega$, that is the function satisfying, for fixed $y \in \Omega$,
\[
\begin{cases}
-\Delta_x G(\cdot\,,y) = \delta_y & \text{in $\Omega$} \\
G(\cdot\,,y) = 0 &\text{on $\pa \Omega$},
\end{cases}
\]
where $\delta_y$ is the Dirac delta centered in $y$. It is well known that the Green function can be decomposed as
\[
G(x,y) = \frac{1}{N (N-2) \omega_N |x-y|^{N-2}} + H(x,y),
\]
where $\omega_N$ is the volume of the unit ball in $\R^N$, and $H(x,y)$ is the \emph{regular part} of the Green function, defined for fixed $y \in \Omega$ as the solution to
\[
\begin{cases}
-\Delta_x H(\cdot\,,y) = 0 & \text{in $\Omega$} \\
H(\cdot\,,y) = \frac{1}{N (N-2) \omega_N |\cdot \, -y|^{N-2}} &\text{on $\pa \Omega$},
\end{cases}
\]

Finally, we denote the standard inner product and norms in $H_0^1(\Omega)$ by 
\[
\langle u, v \rangle_{H_0^1(\Omega)}:= \int_{\Omega} \nabla u \cdot \nabla v, \quad \|u\|_{H_0^1(\Omega)}:= \left( \langle u,u \rangle\right)^\frac12,
\]
and the $L^q$-norm ($q \ge 1$) by $|\cdot|_{L^q(\Omega)}$. When there is no possibility of misunderstanding, we shall often adopt the simplified notations $\|\cdot\|$ and $|\cdot|_q$, for the sake of brevity.

In the rest of the paper $C$ always denotes a positive constant which can depend on the dimension $N$, on the data $\mu_i$ and $\beta_{ij}$, but not on $\eps$. In two steps we will need to point out that a constant $C$ does not depend on $\beta_{ij}$; in such cases we will explicitly write it. The exact value of $C$ can change from line to line.

In many cases, to estimate some quantity involving $P_\eps U_{\d,\x}$ or $P_\eps \psi^i_{\d,\x}$, it will be necessary to approximate the projections with the original functions, carefully controlling the difference. These kind of results are mainly collected in Appendix \ref{app tools}.

\section{Proof of Theorem \ref{thm: main 1}: the reduction scheme}\label{sec: scheme}

Being interested in positive solutions, instead of problem \eqref{system} we consider
\begin{equation}\label{system 2}
\begin{cases}
-\Delta u_i = \mu_i f(u_i) + \sum_{j \neq i} \beta_{ij} |u_j|^{\frac{p+1}{2}} |u_i|^{\frac{p-3}{2}} u_i & \text{in $\Omega_\eps$} \\
u_i = 0 & \text{on $\pa \Omega_\eps$},
\end{cases} \qquad i =1,\dots,m,
\end{equation}
where $f(s) := |s|^{p-1} s^+$. We shall see that this replacement makes possible to prove positivity of the solutions in a very simple way, under the assumptions described by the main theorems. 

Let $i:H_0^1(\Omega_\eps) \to L^{2^*}(\Omega_\eps)$ be the canonical Sobolev embedding. We consider the adjoint operator $i^*: L^{\frac{2N}{N+2}}(\Omega_\eps) \to H_0^1(\Omega_\eps)$, characterized by
\[
i^*(u) = v \quad \iff \quad \begin{cases} -\Delta v= u & \text{in $\Omega$ (in weak sense)} \\ v \in H_0^1(\Omega) \end{cases}
\]
It is well known that $i^*$ is a continuous operator, and using it we can rewrite \eqref{system 2} as
\begin{equation}\label{pb 1}
u_i = i^*\left( \mu_i f(u_i) + \sum_{j \neq i} \beta_{ij} |u_j|^{\frac{p+1}{2}} |u_i|^{\frac{p-3}{2}} u_i\right).
\end{equation}
We search for solutions as perturbation of bubbles centered at different points: let $\eta \in (0,1)$ be small, and let 
\begin{equation}\label{def X}
X_{\eta}:=\left\{ (\mf{d},\bs{\tau}) \in \R^m \times (\R^N)^m: \ \eta < d_i< \eta^{-1}, \ |\tau_i|<\eta^{-1} \right\}.
\end{equation}
Our ansatz is that
\begin{equation}\label{17ott}
u_i = \mu_i^{-\frac{1}{p-1}} P_\eps U_{\d_i,\x_i} + \phi_i,
\end{equation}
where for some $(\mf{d},\bs{\tau}) = (d_1,\dots,d_m, \tau_1,\dots,\tau_m) \in X_\eta$ we have
\begin{equation}\label{asympt expansion}
\d_i := d_i \sqrt{\eps}, \quad \xi_i:= a_i + d_i \sqrt{\eps} \tau_i.
\end{equation}

We stress that the quantity $\eta$ will always be fixed and small, while the unknowns will be $d_i>0$, $\tau_i \in \R^N$, and $\phi_i \in H_0^1(\Omega_\eps)$. 

\begin{remark}
Once that $\eta>0$ is fixed, we observe that 
\[
\d_i = O(\eps^\frac12)\quad \text{and} \quad \eps = O(\d_i^2) \quad \forall i=1,\dots,m
\]
as $\eps \to 0$.
\end{remark}

Plugging ansatz \eqref{17ott} into \eqref{pb 1}, our problem is transformed in the research of $d_i$, $\tau_i$, and $\phi_i$ such that
\begin{multline}\label{pb ansatz}
\mu_i^{-\frac{1}{p-1}} P_\eps U_{\d_i,\x_i} + \phi_i = i^*\Bigg[ \mu_i f(\mu_i^{-\frac{1}{p-1}} P_\eps U_{\d_i,\x_i} + \phi_i)  \\
 + \sum_{j \neq i} \beta_{ij} |\mu_j^{-\frac{1}{p-1}} P_\eps U_{\d_j,\x_j} + \phi_j|^{\frac{p+1}{2}} |\mu_i^{-\frac{1}{p-1}} P_\eps U_{\d_i,\x_i} + \phi_i|^{\frac{p-3}{2}} (\mu_i^{-\frac{1}{p-1}} P_\eps U_{\d_i,\x_i} + \phi_i)\Bigg]
\end{multline}
for $i=1,\dots,m$, with each equality which takes place in $H_0^1(\Omega_\eps)$. To proceed, the idea is then to split the space into two orthogonal subspaces, one of them having finite dimension. To be precise, for $\eps>0$, $d_1,\dots, d_m>0$, and $\tau_1,\dots,\tau_m \in \R^N$, we define
\[
K_i = K_{d_i,\tau_i,\eps}:= \textrm{span}\left\{P_\eps \psi_{\d_i,\xi_i}^h: \ h=0,\dots,N\right\}, \quad \mf{K}_{\mf{d},\bs{\tau},\eps}:= K_1 \times \cdots \times K_m
\]
(recall that $\d_i$ and $\xi_i$ are determined by $d_i, \tau_i$ and $\eps$ through the ansatz \eqref{asympt expansion}).
%and 
%\[
%K_{\mf{\d},\mf{\xi},\eps}:= K_1 \times \cdots \times K_m, \qquad K_{\mf{\d},\mf{\xi},\eps}^{\perp}:= K_1^\perp  \times \cdots \times K_m^\perp.
%\]
Notice that $K_{\mf{d},\bs{\tau},\eps}^{\perp} = K_1^\perp \times \cdots \times K_m^\perp$. 

If $\Pi_i = \Pi_{\d_i,\xi_i,\eps}$ (resp. $\Pi_i^{\perp} = \Pi_{\d_i,\xi_i,\eps}^{\perp}$) denotes the orthogonal projection $H_0^1(\Omega_\eps) \to K_i$ (resp. $H_0^1(\Omega_\eps) \to K_i^{\perp}$), then \eqref{pb ansatz} can be further rewritten as a system of $2m$ equations
\begin{multline}\label{pb K}
\Pi_i(\mu_i^{-\frac{1}{p-1}} P_\eps U_{\d_i,\x_i} + \phi_i )= (\Pi_i \circ i^*)\Bigg[ \mu_i f(\mu_i^{-\frac{1}{p-1}} P_\eps U_{\d_i,\x_i} + \phi_i) \\
 + \sum_{j \neq i} \beta_{ij} |\mu_j^{-\frac{1}{p-1}} P_\eps U_{\d_j,\x_j} + \phi_j|^{\frac{p+1}{2}} |\mu_i^{-\frac{1}{p-1}} P_\eps U_{\d_i,\x_i} + \phi_i|^{\frac{p-3}{2}} (\mu_i^{-\frac{1}{p-1}} P_\eps U_{\d_i,\x_i} + \phi_i)\Bigg],
\end{multline}
and
\begin{multline}\label{pb K ort}
\Pi_i^{\perp}(\mu_i^{-\frac{1}{p-1}} P_\eps U_{\d_i,\x_i} + \phi_i )= (\Pi_i^{\perp} \circ i^*)\Bigg[ \mu_i f(\mu_i^{-\frac{1}{p-1}} P_\eps U_{\d_i,\x_i} + \phi_i)  \\
 + \sum_{j \neq i} \beta_{ij} |\mu_j^{-\frac{1}{p-1}} P_\eps U_{\d_j,\x_j} + \phi_j|^{\frac{p+1}{2}} |\mu_i^{-\frac{1}{p-1}} P_\eps U_{\d_i,\x_i} + \phi_i|^{\frac{p-3}{2}} (\mu_i^{-\frac{1}{p-1}} P_\eps U_{\d_i,\x_i} + \phi_i)\Bigg],
\end{multline}                                                       
$i=1,\dots,m$.

The proof of Theorem \ref{thm: main 1} consists of two main steps: first, for fixed $\eps$, $d_1,\dots,d_m$, and $\tau_1,\dots,\tau_m$ we solve the $m$-equations system \eqref{pb K ort}, finding $(\phi_1^{\mf{d},\bs{\tau},\eps},\dots,\phi_m^{\mf{d},\bs{\tau},\eps}) \in \mf{K}_{\mf{d},\bs{\tau},\eps}^\perp$. Plugging this choice of $\bs{\phi}$ into \eqref{pb K}, we obtain a finite dimensional problem in the unknowns $\mf{d}$ and $\bs{\tau}$, which will be solved in a second step for any $\eps>0$ sufficiently small.

\section{The equations in $\mf{K}_{\mf{d},\bs{\tau},\eps}^\perp$}\label{sec: eq ort}

In this section we study the solvability of \eqref{pb K ort}. In a functional analytic perspective, \eqref{pb K ort} reads
\begin{equation}\label{L=N+R}
L_{\mf{d},\bs{\tau},\eps}^i(\bs{\phi}) = N_{\mf{d},\bs{\tau},\eps}^i(\bs{\phi}) + R_{\mf{d},\bs{\tau},\eps}^i,
\end{equation}
where $L$ stays for the linear part
\begin{equation}\label{def L}
\begin{split}
L_{\mf{d},\bs{\tau},\eps}^i(\bs{\phi}) & = \Pi_i^{\perp} \Bigg\{ \phi_i -i^*\Bigg[ \mu_i f'(\mu_i^{-\frac{1}{p-1}} P_\eps U_{\d_i,\xi_i}) \phi_i   \\
&  + \frac{p-1}{2} \sum_{j \neq i} \beta_{ij} (\mu_j^{-\frac{1}{p-1}} P_\eps U_{\d_j,\xi_j})^{\frac{p+1}{2}}   
(\mu_i^{-\frac{1}{p-1}} P_\eps U_{\d_i,\xi_i})^{\frac{p-3}{2}}  \phi_i \\
&  + \frac{p+1}{2} \sum_{j \neq i} \beta_{ij} (\mu_j^{-\frac{1}{p-1}} P_\eps U_{\d_j,\xi_j})^{\frac{p-1}{2}}   
(\mu_i^{-\frac{1}{p-1}} P_\eps U_{\d_i,\xi_i})^{\frac{p-1}{2}}  \phi_j \Bigg]\Bigg\},
\end{split}
\end{equation} 
$N$ stays for the nonlinear part
\begin{equation}\label{def N}
\begin{split}
N_{\mf{d},\bs{\tau},\eps}^i &(\bs{\phi}) = \\
& \Pi_i^{\perp} \circ i^* \Bigg[ \mu_i f(\mu_i^{-\frac{1}{p-1}} P_\eps U_{\d_i,\xi_i}+\phi_i)-\mu_i f(\mu_i^{-\frac{1}{p-1}} P_\eps U_{\d_i,\xi_i}) - \mu_i f'(\mu_i^{-\frac{1}{p-1}} P_\eps U_{\d_i,\xi_i}) \phi_i  \\
& + \sum_{j \neq i} \beta_{ij} |\mu_j^{-\frac{1}{p-1}} P_\eps U_{\d_j,\x_j} + \phi_j|^{\frac{p+1}{2}} |\mu_i^{-\frac{1}{p-1}} P_\eps U_{\d_i,\x_i} + \phi_i|^{\frac{p-3}{2}} (\mu_i^{-\frac{1}{p-1}} P_\eps U_{\d_i,\x_i} + \phi_i) \\
& - \sum_{j \neq i} \beta_{ij} (\mu_j^{-\frac{1}{p-1}} P_\eps U_{\d_j,\x_j})^{\frac{p+1}{2}}(\mu_i^{-\frac{1}{p-1}} P_\eps U_{\d_i,\x_i})^{\frac{p-1}{2}} \\
& - \frac{p-1}{2} \sum_{j \neq i} \beta_{ij} (\mu_j^{-\frac{1}{p-1}} P_\eps U_{\d_j,\xi_j})^{\frac{p+1}{2}}   
(\mu_i^{-\frac{1}{p-1}} P_\eps U_{\d_i,\xi_i})^{\frac{p-3}{2}}  \phi_i \\
& - \frac{p+1}{2} \sum_{j \neq i} \beta_{ij} (\mu_j^{-\frac{1}{p-1}} P_\eps U_{\d_j,\xi_j})^{\frac{p-1}{2}}   
(\mu_i^{-\frac{1}{p-1}} P_\eps U_{\d_i,\xi_i})^{\frac{p-1}{2}}  \phi_j \Bigg],
\end{split}
\end{equation}
and $R$ is the remainder term
\begin{equation}\label{def R}
\begin{split}
R_{\mf{d},\bs{\tau},\eps}^i & = \Pi_i^{\perp} \Bigg\{ -\mu_i^{-\frac{1}{p-1}} P_\eps U_{\d_i,\x_i}   \\
& \hphantom{= \Pi_i^{\perp} \Bigg\{ \ }+ i^*\Bigg[ \mu_i f(\mu_i^{-\frac{1}{p-1}} P_\eps U_{\d_i,\xi_i}) + \sum_{j \neq i} \beta_{ij} (\mu_j^{-\frac{1}{p-1}} P_\eps U_{\d_j,\x_j})^{\frac{p+1}{2}}(\mu_i^{-\frac{1}{p-1}} P_\eps U_{\d_i,\x_i})^{\frac{p-1}{2}} \Bigg]\Bigg\} \\
& = \Pi_i^{\perp} \circ i^* \Bigg[ \mu_i^{-\frac{1}{p-1}} \left( P_\eps U_{\d_i,\x_i}^p- U_{\d_i,\x_i}^p\right)  +  \sum_{j \neq i} \beta_{ij} (\mu_j^{-\frac{1}{p-1}} P_\eps U_{\d_j,\x_j})^{\frac{p+1}{2}}(\mu_i^{-\frac{1}{p-1}} P_\eps U_{\d_i,\x_i})^{\frac{p-1}{2}} \Bigg],
\end{split}
\end{equation}
where the last equality is a consequence of the definitions of $i^*$ and of $f$.

For future convenience, we also define 
\[
\mf{L}_{\mf{d},\bs{\tau},\eps} := (L_{\mf{d},\bs{\tau},\eps}^1, \dots,L_{\mf{d},\bs{\tau},\eps}^m): \mf{K}_{\mf{d},\bs{\tau},\eps}^\perp \to \mf{K}_{\mf{d},\bs{\tau},\eps}^\perp,
\]
and $\mf{R}_{\mf{d},\bs{\tau},\eps}$ and $\mf{N}_{\mf{d},\bs{\tau},\eps} $ in an analogue way.

%Given $\eta>0$ small, let us set
%\begin{equation}\label{def X}
%X_{\eta}:=\left\{ (\mf{d},\bs{\tau}) \in \R^m \times (\R^N)^m: \ \eta < d_i< \eta^{-1}, \ |\tau_i|<\eta^{-1} \right\}.
%\end{equation}

The main result of this section is the following:
\begin{proposition}\label{prop: pb K ort}
Let $N=4$. For every $\eta>0$ small enough there exists $\bar \beta, \eps_0>0$ small, and $C>0$, such that if $\eps \in (0,\eps_0)$, and $-\infty <\beta_{ij}<\bar \beta$ for every $i \neq j$, then for any $(\mf{d},\bs{\tau}) \in X_\eta$ (see \eqref{def X}) there exists a unique function $\bs{\phi}^{\mf{d},\bs{\tau},\eps} \in K_{\mf{d},\bs{\tau},\eps}^\perp$ solving the equation
\begin{equation}\label{nonlin eq}
\mf{L}_{\mf{d},\bs{\tau},\eps}(\bs{\phi}) = \mf{R}_{\mf{d},\bs{\tau},\eps} + \mf{N}_{\mf{d},\bs{\tau},\eps}(\bs{\phi})
\end{equation}
and satisfying
\[
\|\bs{\phi}^{\mf{d},\bs{\tau},\eps}\|_{H_0^1(\Omega_\eps)} \le C \eps^{\frac{N-2}2}.
\] 
Furthermore, the map $(\eps,\mf{d},\bs{\tau}) \mapsto \bs{\phi}^{\mf{d},\bs{\tau},\eps}$ is of class $\mathcal{C}^1$, and
\[
\|\nabla_{(\mf{d},\bs{\tau})} \bs{\phi}^{\mf{d},\bs{\tau},\eps}\|_{H_0^1(\Omega_\eps)} \le C \eps^{\frac{N-2}2}. 
\] 
If $N=3$, the same conclusion holds without any restriction on $\beta_{ij}$.
\end{proposition}

The proof of the proposition takes the rest of this section, and is divided into several intermediate lemmas. 

\subsection{Study of the linear part}

As a first step, it is important to understand the solvability of the linear problem associated to \eqref{L=N+R}, i.e. 
\begin{equation}\label{linear pb}
L_{\mf{d},\bs{\tau},\eps}^i(\bs{\phi}) = f_i, \quad \text{with} \quad f_i \in K_{d_i,\tau_i,\eps}^\perp.
\end{equation}

\begin{lemma}\label{lem: linear part}
Let $N=4$. For every $\eta>0$ small enough there exists $\bar \beta, \eps_0>0$ small, and $C>0$, such that if $\eps \in (0,\eps_0)$, and $-\infty <\beta_{ij}<\bar \beta$ for every $i \neq j$, then
\begin{equation}\label{17ott1}
\|\mf{L}_{\mf{d},\bs{\tau},\eps}(\bs{\phi})\|_{H_0^1(\Omega_\eps)} \ge C \|\bs{\phi}\|_{H_0^1(\Omega_\eps)}  \qquad \forall \bs{\phi} \in H_0^1(\Omega_\eps,\R^m)
\end{equation}
for every $(\mf{d}, \bs{\tau}) \in X_\eta$. Moreover, $\mf{L}_{\mf{d},\bs{\tau},\eps}$ is invertible in $\mf{K}_{\mf{d},\bs{\tau},\eps}^\perp$, with continuous inverse. \\
If $N=3$, the same conclusion holds true without restrictions on $\beta_{ij}$.
\end{lemma}

\begin{proof}
The long proof proceed by contradiction. Let us suppose that there exist sequences \[
\{\eps_n\} \subset \R^+, \ \eps_n \to 0, \  \{(\mf{d}_n,\boldsymbol{\tau}_n)\} \subset X_\eta,\  \{\bs{\phi}_n\} \subset K_{1,n}^\perp \times \cdots \times K_{m,n}^{\perp}
\]
such that
\[
\| \bs{\phi}_n\|_{H_0^1(\Omega_{\eps_n})}=1 \quad \text{and} \quad \|\mf{L}_{n}(\bs{\phi}_n)\|_{H_0^1(\Omega_{\eps_n})} \to 0
\]
as $n \to \infty$, where we wrote $K_{i,n}:= K_{d_{i,n},\tau_{i,n},\eps_n}$ and $\mf{L}_n:= \mf{L}_{\mf{d}_n, \bs{\tau}_n, \eps_n}$ for short. In the same spirit in this proof we write $P_n:= P_{\eps_n}$, $U_{i,n}:= U_{\d_{i,n},\xi_{i,n}}$, $ \psi_{i,n}^h:= \psi_{\d_{i,n},\xi_{i,n}}^h$, and $\Omega_n:= \Omega_{\eps_n}$.

Let $\mf{h}_n:= \mf{L}_{n}(\bs{\phi}_n)$. Then, observing that 
\[
\mu_i f'(\mu_i^{-\frac{1}{p-1}} P_\eps U_{\d,\xi}) = p (P_\eps U_{\d,\x})^{p-1} \qquad \forall \eps,\d>0,\ \xi \in \R^N,
\] 
we have by definition of $\mf{L}_n$
\begin{equation}\label{eq phi}
\begin{split}
 \phi_{i,n}   & =  i^*\Bigg[ p(P_n U_{i,n})^{p-1} \phi_{i,n}   + \frac{p-1}{2} \sum_{j \neq i} \beta_{ij} (\mu_j^{-\frac{1}{p-1}} P_n U_{j,n})^{\frac{p+1}{2}}   
(\mu_i^{-\frac{1}{p-1}} P_n U_{i,n})^{\frac{p-3}{2}}  \phi_{i,n} \\
& \hphantom{= i^*\Bigg[ \ } + \frac{p+1}{2} \sum_{j \neq i} \beta_{ij} (\mu_j^{-\frac{1}{p-1}} P_n U_{j,n})^{\frac{p-1}{2}}   
(\mu_i^{-\frac{1}{p-1}} P_n U_{i,n})^{\frac{p-1}{2}}  \phi_{j,n} \Bigg] +h_{i,n} - w_{i,n}
\end{split}
\end{equation}
for some $w_{i,n} \in K_{i,n}$. 

\smallskip

\noindent \textbf{Step 1)} We show that $\|w_{i,n}\|_{H_0^1(\Omega_n)} \to 0$ as $n \to \infty$. 

Since $w_{i,n} \in K_{i,n}$, there exists constants $c_{i,n}^k$ such that
\[
w_{i,n} = \sum_{k=0}^N c_{i,n}^k P_n \psi^k_{i,n}.
\]
Let us multiply equation \eqref{eq phi} with $\delta_{i,n}^2 w_{i,n}$: taking into account that $\phi_{i,n}, h_{i,n} \in K_{i,n}^\perp$, we deduce that
\begin{equation}\label{eq phi test} 
\begin{split}
\underbrace{\d_{i,n}^2 \|w_{i,n}\|_{H_0^1(\Omega_n)}^2}_{=:(I)} 
% \sum_{k=0}^N  c_{i,n}^l  \underbrace{\delta_{i,n}^2\int_{\Omega_n}\nabla (P_n \psi_{i,n}^l) \cdot  \nabla (P_n \psi_{i,n}^k)}_{=:(I)} 
= p  \underbrace{\delta_{i,n}^2 \int_{\Omega_n} (P_n U_{i,n})^{p-1} \phi_{i,n} w_{i,n}}_{=:(II)} \\
+   \frac{p-1}{2} \sum_{l=0}^N c_{i,n}^l\underbrace{ \delta_{i,n}^2 \int_{\Omega_n}  \sum_{j \neq i} \beta_{ij} (\mu_j^{-\frac{1}{p-1}} P_n U_{j,n})^{\frac{p+1}{2}}   
(\mu_i^{-\frac{1}{p-1}} P_n U_{i,n})^{\frac{p-3}{2}}  \phi_{i,n} (P_n \psi_{i,n}^l) }_{=:(III)}\\
+ \frac{p+1}{2} \sum_{l=0}^N c_{i,n}^l \underbrace{ \delta_{i,n}^2 \int_{\Omega_n}   \sum_{j \neq i} \beta_{ij} (\mu_j^{-\frac{1}{p-1}} P_n U_{j,n})^{\frac{p-1}{2}}   
(\mu_i^{-\frac{1}{p-1}} P_n U_{i,n})^{\frac{p-1}{2}}  \phi_{j,n} (P_n \psi_{i,n}^l)}_{=:(IV)}.
\end{split}
\end{equation}
The rest of the proof of step 1 consists in a careful (and very long) asymptotic expansion of the terms ($I$)-($IV$), whose details are contained in Appendix \ref{app step 1 linear}. Therein we prove that
\begin{equation}\label{estimates step 1}
\begin{split}
|(I)| & =  \sum_{l= 0}^N (c_{i,n}^l)^2 \sigma_{ll} + o(1) \sum_{l,k=0}^N c_{i,n}^l c_{i,n}^k  \\
|(II)| & = o(\d_{i,n}^2) \|w_{i,n}\|_{H_0^1(\Omega_n)} + O(\d_{i,n}^2) \sum_{l=0}^N c_{i,n}^l \\
|(III)| & = o(\d_{i,n}^2) \\
|(IV)| & = o(\d_{i,n}^2),
\end{split}
\end{equation}
as $n \to \infty$, with 
\begin{equation}\label{sigma lk}
\sigma_{lk} = \begin{cases} 0 & \text{if $l \neq k$} \\  p \alpha_N^{p+1} (N-2)^2 \int_{\R^N} \frac{y_l^2}{(1+|y|^2)^{N+2}}\,dy & \text{if $k=l \ge 1$} \\  p  \alpha_N^{p+1} \left( \frac{N-2}{2}\right)^2 \int_{\R^N} \frac{(|y|^2-1)^2}{(1+|y|^2)^{N+2}}\,dy & \text{if $k=l=0$}.
\end{cases}
\end{equation}
Using the second to fourth estimates in \eqref{estimates step 1}, equation \eqref{eq phi test} becomes
\begin{equation}\label{261}
\d_{i,n}^2 \|w_{i,n}\|_{H_0^1(\Omega_n)}^2 =   o(\d_{i,n}^2)  \|w_{i,n}\|_{H_0^1(\Omega_n)} + O(\d_{i,n}^2) \sum_{k=0}^N c_{i,n}^k.
\end{equation}
Due to the first estimate in \eqref{estimates step 1}, the previous expression yields
\[
C \left(\sum_{k=0}^N |c_{i,n}^k|\right)^2 \le o(1) \sum_{k=0}^N |c_{i,n}^k|,   
%\sum_{k=0}^N \sigma_{kk} (c_{i,n}^k)^2 \le o(\d_{i,n}) \left(\sum_{k=0}^N \sigma_{kk} (c_{i,n}^k)^2\right)^{\frac12} + o(\d_{i,n}^2) \sum_{k=0}^N |c_{i,n}^k|,
\]
so that, firstly, $\{c_{i,n}^k\}$ is a bounded sequence, for any $k$, and in turn, this implies that $c_{i,n}^k \to 0$ as $n \to \infty$, for every $k$. Hence, using this into \eqref{261}, we deduce that
\[
\|w_{i,n}\|_{H_0^1(\Omega_n)}^2 = o(1) \|w_{i,n}\|_{H_0^1(\Omega_n)} + o(1), 
\]
whence $\|w_{i,n}\|_{H_0^1(\Omega_n)} \to 0$ as $n \to \infty$.

\smallskip

\noindent \textbf{Step 2)} For a fixed $\kappa=1,\dots,m$, let us introduce
\[
\tilde \phi_{i,n}^\kappa (y):= \begin{cases} \d_{\kappa,n}^{\frac{N-2}{2}} \phi_{i,n}(\x_{\kappa,n}+\d_{\kappa,n} y) & y \in \frac{\Omega_n-\x_{\kappa,n}}{\d_{\kappa,n}}=:\tilde \Omega_{\kappa,n} \\
0 & y \in \R^N \setminus \tilde\Omega_{\kappa,n},
\end{cases}   \quad i=1,\dots,m.
\]
In a completely analogue way, we define $\tilde h_{i,n}^\kappa$ and $\tilde w_{i,n}^\kappa$, and we set $\tilde \phi_{\kappa,n}:= \tilde \phi_{\kappa,n}^\kappa$, $\tilde h_{\kappa,n}:= \tilde h_{\kappa,n}^\kappa$, $\tilde w_{\kappa,n}:= \tilde w_{\kappa,n}^\kappa$. 
%In what followsthat the index $i$ of the different components $\phi_{i,n}$ varies, while the index $\kappa$ of the parameters $\d_{\kappa,n}$ and $\xi_{\kappa,n}$ is fixed. 

In this step we show that $\tilde \phi_{\kappa,n} \wc 0$ in $\mathcal{D}^{1,2}(\R^N)$ (i.e. $\nabla \tilde \phi_{\kappa,n} \wc 0$ in $L^2(\R^N)$) as $n \to \infty$, for every $\kappa=1,\dots,m$.

At first, we observe that $\|\tilde \phi_{\kappa,n}\|_{H_0^1(\tilde \Omega_{\kappa,n})} = \|\phi_{\kappa,n}\|_{H_0^1(\Omega_n)} \le 1$, and hence up to a subsequence $\tilde \phi_{\kappa,n} \wc \tilde \phi_\kappa$ weakly in $\mathcal{D}^{1,2}(\R^N)$ for every $i$. Now we rewrite the equation \eqref{eq phi} for $\phi_{\kappa,n}$ in terms of $\tilde \phi_{\kappa,n}$: if $\psi \in \C^\infty_c(\R^N)$, we have 
\begin{equation}\label{eq phi tilde}
\begin{split}
\int_{\tilde \Omega_{\kappa,n}} \nabla \tilde \phi_{\kappa,n}  &\cdot \nabla \psi = p \d_{\kappa,n}^2 \int_{\tilde \Omega_{\kappa,n}} (\widehat{P_n U_{\kappa,n}})^{p-1}  \tilde \phi_{\kappa,n} \psi \\
& + \frac{p-1}{2} \sum_{j \neq \kappa} \d_{\kappa,n}^2 \beta_{\kappa j} \int_{\tilde \Omega_{\kappa,n}} (\mu_\kappa^{-\frac{1}{p-1}} \widehat{P_n U_{\kappa,n}})^{\frac{p-3}{2}}  (\mu_j^{-\frac{1}{p-1}} \widehat{P_n U_{j,n}})^{\frac{p+1}{2}}  \tilde \phi_{\kappa,n} \psi \\
& + \frac{p+1}{2} \sum_{j \neq \kappa} \d_{\kappa,n}^2 \beta_{\kappa j} \int_{\tilde \Omega_{\kappa,n}} (\mu_\kappa^{-\frac{1}{p+1}} \widehat{P_n U_{\kappa,n}})^{\frac{p-1}{2}} (\mu_j^{-\frac{1}{p-1}} \widehat{P_n U_{j,n}})^{\frac{p-1}{2}}\tilde \phi_{j,n} \psi \\
& + \int_{\tilde \Omega_{\kappa,n}} \nabla \tilde \phi_{\kappa,n} \cdot \nabla (\tilde h_{\kappa,n} - \tilde w_{\kappa,n}),
\end{split}
\end{equation}
where $\widehat{P_n U_{j,n}}(y):= P_n U_{j,n}(\x_{\kappa,n} +\d_{\kappa,n} y)$ for all $j=1,\dots,m$. In Appendix \ref{app est step 2}, we show that \eqref{eq phi tilde} yields
\begin{equation}\label{estimate step 2}
\int_{\tilde{\Omega}_{\kappa,n}} \nabla \tilde \phi_{\kappa,n}  \cdot \nabla \psi = p \int_{\R^N} U_{1,0}^{p-1} \tilde \phi_{\kappa} \psi + o(1) \qquad \text{for every $\psi \in \C^\infty_c(\R^N)$}
\end{equation}
which in turn, by weak convergence, implies that 
\begin{equation}\label{eq phi lim}
-\Delta \tilde \phi_\kappa = p U_{1,0}^{p-1} \tilde \phi_\kappa, \qquad \tilde \phi_i \in \D^{1,2}(\R^N).
\end{equation}
Since our final goal consists in proving that $\tilde \phi_\kappa=0$, due to the previous equation it will suffices to show that $\tilde \phi_\kappa$ is orthogonal, in $\D^{1,2}(\R^N)$, to the $N+1$ partial derivatives $\psi_{1,0}^\ell$ ($\ell=0,\dots, N$) of $U_{1,0}$. Indeed, we already know that these partial derivatives span the sets of the solutions to \eqref{eq phi lim}. The orthogonality condition comes from the fact that $\phi_{\kappa,n} \in K_{\kappa,n}^\perp$ for every $i$ and $n$: indeed, for any $\ell=1,\dots,N$ (the case $\ell=0$ is analogue), we have
\begin{align*}
0 & = \d_{\kappa,n} \int_{\Omega_{n}} \nabla \phi_{\kappa,n} \cdot \nabla(P_n \psi_{\kappa,n}^\ell) = p \d_{\kappa,n} \int_{\Omega_{n}} p U_{\kappa,n}^{p-1} \psi_{\kappa,n}^\ell \phi_{\kappa,n} \\
& = \d_{\kappa,n} \int_{\Omega_{n}} p (N-2) \alpha_N^p \left( \frac{\d_{\kappa,n}}{\d_{\kappa,n}^2 + |x-\xi_{\kappa,n}|^2}\right)^2 \d_{\kappa,n}^{\frac{N-2}{2}} \frac{x_\ell-\xi_{\kappa,n}^\ell}{\left( \d_{\kappa,n}^2 + |x-\x_{\kappa,n}|^2 \right)^{\frac{N}{2}}} \, dx\\
& = \int_{\tilde \Omega_{\kappa,n}} p(N-2) \alpha_N^p \frac{1}{(1+|y|^2)^2} \frac{y_\ell}{(1+|y|^2)^{\frac{N}{2}}}\tilde \phi_{\kappa,n} \,dy \\
& = \int_{\tilde \Omega_{\kappa,n}} p U_{1,0}^{p-1} \psi_{1,0}^\ell \tilde \phi_{\kappa,n},
\end{align*}
and by weak convergence we deduce that for every $\ell$ and $i$
\[
0= \int_{\R^N} p U_{1,0}^{p-1} \psi_{1,0}^\ell \tilde \phi_{\kappa} = \int_{\R^N} \nabla \psi_{1,0}^\ell \cdot \nabla \tilde \phi_\kappa,
\]
as desired.

\smallskip

\noindent \textbf{Step 3)} We prove that $\|\phi_{i,n}\|_{H_0^1(\Omega_n)} \to 0$ as $n \to \infty$ for every $i$. This is in contradiction with the fact that $\|\bs{\phi}_{n}\|_{H_0^1(\Omega_n)}=1$, and completes the proof of \eqref{17ott1}. Let us test \eqref{eq phi} with $\phi_{i,n}$: recalling that $\{\phi_{i,n}\}$ is bounded in $H_0^1(\Omega_n)$ and that $w_{i,n}, h_{i,n} \to 0$ strongly, we deduce that
\begin{equation}\label{eq phi phi}
\begin{split}
\|\phi_{i,n}\|_{H_0^1(\Omega_n)}^2 = o(1) + \underbrace{p \int_{\Omega_n} (P_n U_{i,n})^{p-1} \phi_{i,n}^2}_{=:(I)} \\
+ \frac{p-1}{2} \sum_{j \neq i} \underbrace{\beta_{ij} \int_{\Omega_{n}} (\mu_i^{-\frac{1}{p-1}} P_n U_{i,n})^{\frac{p-3}{2}}  (\mu_j^{-\frac{1}{p-1}} P_n U_{j,n})^{\frac{p+1}{2}}   \phi_{i,n}^2}_{=:(II)}  \\
 + \frac{p+1}{2} \sum_{j \neq i} \underbrace{ \beta_{ij} \int_{\Omega_{n}} (\mu_i^{-\frac{1}{p-1}} P_n U_{i,n})^{\frac{p-1}{2}} (\mu_j^{-\frac{1}{p-1}} P_n U_{j,n})^{\frac{p-1}{2}}\phi_{i,n} \phi_{j,n} }_{=:(III)}
\end{split}
\end{equation}
We have to estimate the right hand side. At first, recalling as usual that $0 \le P_n U_{i,n} \le U_{i,n}$, we have
\begin{equation}\label{3061}
|(I)| \le p\int_{\Omega_n} 
U_{i,n}^{p-1} \phi_{i,n}^2  = C p \int_{\tilde \Omega_{i,n}} U_{1,0}^{p-1} \tilde \phi_{i,n}^2 \to 0
\end{equation}
as $n \to \infty$, since $\tilde \phi_{i,n}^2 \wc 0$ in $L^{\frac{N}{N-2}}(\R^N)$ by step 2, and $U_{1,0}^{p-1} \in L^{\frac{N}{2}}(\R^N)$. 

The second term on the right hand side in \eqref{eq phi phi} can be estimated discussing several possibilities: if $\beta_{ij} <0$, we simply observe that $(II)  \le 0$. Otherwise, if $\beta_{ij}>0$,
\begin{align*}
(II)
& \le C  \beta_{ij} \int_{\Omega_n} U_{i,n}^\frac{p-3}{2} U_{j,n}^\frac{p+1}{2} \phi_{i,n}^2 \\
& \le C \beta_{ij} \left( \int_{\Omega_n} U_{i,n}^\frac{(4-N)N}{2(N-2)} U_{j,n}^\frac{N^2}{2(N-2)} \right)^\frac{2}{N} |\phi_{i,n}|_{2^*}^2.
\end{align*}
In case $N=3$, the last term reads
\begin{align*}
C \beta_{ij} \left( \int_{\Omega_n} U_{i,n}^\frac{3}{2} U_{j,n}^\frac{9}{2} \right)^\frac{2}{3} |\phi_{i,n}|_{2^*}^2 \le C \beta_{ij}\|\phi_{i,n}\|_{H_0^1(\Omega_n)}^2 \left(\d_{i,n}^\frac34 \delta_{j,n}^\frac34\right)^\frac23 = o(1)
\end{align*}
(no matter how large $\beta_{ij}$ is), where we used Lemmas \ref{lem a.4} and \ref{lem a.5}. If on the other hand $N=4$, we have
\begin{align*}
C  \beta_{ij} \left( \int_{\Omega_n} U_{j,n}^4 \right)^\frac{1}{2} |\phi_{i,n}|_{2^*}^2 &\le C  \beta_{ij}\|\phi_{i,n}\|_{H_0^1(\Omega_n)}^2,
\end{align*}
% \\
%& \le \frac12 \|\phi_{i,n}\|_{H_0^1(\Omega_n)}^2
%\end{align*}
where $C$ does not depend on $\beta_{ij}$. We conclude that there exists $\bar \beta>0$ sufficiently small such that, if $\beta_{ij} < \bar \beta$ for all $i \neq j$, then
\begin{equation}\label{N=4 beta piccolo}
|(II)| \le C   \beta_{ij}\|\phi_{i,n}\|_{H_0^1(\Omega_n)}^2 \le \frac12 \|\phi_{i,n}\|_{H_0^1(\Omega_n)}^2
\end{equation}
To sum up, in any case we can conclude that for every $i \neq j$ and for any $n$ large enough
\begin{equation}\label{3062}
\frac{p-1}{2} \sum_{j \neq i} (II) 
\le \frac12 \|\phi_{i,n}\|_{H_0^1(\Omega_n)}^2 \quad \begin{cases} \text{for every $\beta_{ij}$, if $N=3$} \\ \text{provided $\beta_{ij}< \bar \beta$, if $N=4$}.
\end{cases}
\end{equation}

It remains to consider the last term in \eqref{eq phi phi}: by the H\"older inequality and Lemmas \ref{lem a.4} and \ref{lem a.5}
\begin{equation}\label{3063}
\begin{split}
|(III)| &\le \beta_{ij} \int_{\Omega_n} U_{i,n}^\frac{p-1}{2} U_{j,n}^\frac{p-1}{2} |\phi_{i,n}| | \phi_{j,n} | \\
& \le \beta_{ij} |\phi_{i,n}|_{2^*} |\phi_{j,n}|_{2^*} \left( \int_{\Omega_n} U_{i,n}^\frac{N}{N-2} U_{j,n}^{\frac{N}{N-2}} \right)^\frac2N \\
& \le C\beta_{ij} \|\phi_{i,n}\|_{H_0^1(\Omega_n)} \|\phi_{j,n}\|_{H_0^1(\Omega_n)} \d_{i,n} \d_{j,n} \left( |\log{\dn}| + |\log{\delta_{j,n} }|\right)^\frac2N = o(1),
\end{split}
\end{equation}
as $n \to \infty$, where we used the boundedness of $\{\boldsymbol{\phi}_n\}$ in $H_0^1$. 

Plugging \eqref{3061}-\eqref{3063} into \eqref{eq phi phi}, we conclude that $\phi_{i,n} \to 0$ strongly in $H_0^1(\Omega_n)$, which gives the desired contradiction and completes the first part of the lemma. 

It remains still to show the invertibility of $\mf{L}_{\mf{d},\bs{\tau},\eps}$, and this is the object of the last step.

\smallskip

\noindent \textbf{Step 4)} We start recalling that the operator $i^*:L^{\frac{2N}{N+2}}(\Omega_\eps) \to H_0^1(\Omega_\eps)$ is compact. Therefore, by definition, the restriction of $\mf{L}_{\mf{d},\bs{\tau},\eps}$ to $\mf{K}_{\mf{d},\bs{\tau},\eps}^\perp$ is a compact perturbation of the identity. So far we showed that 
\begin{equation}\label{3064}
\|\mf{L}_{\mf{d},\bs{\tau},\eps}(\bs{\phi})\|_{H_0^1(\Omega_\eps, \R^m)} \ge C \|\bs{\phi}\|_{H_0^1(\Omega_\eps, \R^m)} \qquad \text{for every } \bs{\phi} \in \mf{K}_{\mf{d},\bs{\tau},\eps}^\perp,
\end{equation}
and hence $\mf{L}_{\mf{d},\bs{\tau},\eps}$ is injective. By the Fredholm alternative, it is also surjective, thus invertible, and the inverse is continuous (due to \eqref{3064}). 
\end{proof}

Now we prove the solvability of equation \eqref{nonlin eq}.

\begin{lemma}\label{lem: exist nonlin}
Let $N=4$. For every $\eta>0$ small enough there exists $\bar \beta, \eps_1>0$ small and $C>0$ such that: if $\eps \in (0,\eps_1)$, and $-\infty <\beta_{ij}<\bar \beta$ for every $i \neq j$, then for any $(\mf{d},\bs{\tau}) \in X_\eta$ there exists a unique function $\bs{\phi}^{\mf{d},\bs{\tau},\eps} \in K_{\mf{d},\bs{\tau},\eps}^\perp$ solving equation \eqref{nonlin eq}:
\[
\mf{L}_{\mf{d},\bs{\tau},\eps}(\bs{\phi}^{\mf{d},\bs{\tau},\eps}) = \mf{R}_{\mf{d},\bs{\tau},\eps} + \mf{N}_{\mf{d},\bs{\tau},\eps}(\bs{\phi}^{\mf{d},\bs{\tau},\eps}),
\]
and satysfying
\[
\|\bs{\phi}^{\mf{d},\bs{\tau},\eps}\|_{H_0^1(\Omega_\eps)} \le C \eps^\frac{N-2}2.
\] 
If $N=3$, the same conclusion holds without any restriction on $\beta_{ij}$.
\end{lemma}

\begin{proof}
Let $\eps \in (0,\eps_0)$, with $\eps_0$ given by Lemma \ref{lem: linear part}. Then, to solve equation \eqref{nonlin eq} is equivalent to find $\bs{\phi} \in \mf{K}_{\mf{d},\bs{\tau},\eps}^\perp$ such that  
\[
\bs{\phi} = \mf{L}_{\mf{d},\bs{\tau},\eps}^{-1} \left( \mf{R}_{\mf{d},\bs{\tau},\eps} + \mf{N}_{\mf{d},\bs{\tau},\eps}(\bs{\phi}) \right)=:  \mf{T}_{\mf{d},\bs{\tau},\eps}(\bs{\phi}).
\]
We aim at proving that $\mf{T}_{\mf{d},\bs{\tau},\eps}$ is a contraction inside a properly chosen region.

\smallskip

\noindent \textbf{Step 1)} $\mf{T}_{\mf{d},\bs{\tau},\eps}: Y_\eps \to Y_\eps$ for a suitable subset $Y_\eps \subset H_0^1(\Omega_\eps,\R^m)$.

Using the continuity of $\mf{L}_{\mf{d},\bs{\tau},\eps}^{-1}$, of $\Pi_i^\perp$ and of $i^*$, we have
\begin{equation}\label{473}
\| \mf{T}_{\mf{d},\bs{\tau},\eps}(\bs{\phi})\|_{H_0^1(\Omega_\eps)} 
\le C \left( | \tilde{\mf{R}}_{\mf{d},\bs{\tau},\eps}|_{L^{\frac{2N}{N+2}}(\Omega_\eps)} + |\tilde{\mf{N}}_{\mf{d},\bs{\tau},\eps}(\bs{\phi}) |_{L^{\frac{2N}{N+2}}(\Omega_\eps)}\right),
\end{equation}
where 
\begin{align*}
&\tilde{R}^i_{\mf{d},\bs{\tau},\eps} := \mu_i^{-\frac{1}{p-1}} \left( P_\eps U_{\d_i,\x_i}^p- U_{\d_i,\x_i}^p\right) +  \sum_{j \neq i} \beta_{ij} (\mu_j^{-\frac{1}{p-1}} P_\eps U_{\d_j,\x_j})^{\frac{p+1}{2}}(\mu_i^{-\frac{1}{p-1}} P_\eps U_{\d_i,\x_i})^{\frac{p-1}{2}},\\
&\tilde{N}^i_{\mf{d},\bs{\tau},\eps}  (\bs{\phi}) := \tilde{P}^i_{\mf{d},\bs{\tau},\eps}(\bs{\phi}) + \tilde{Q}^i_{\mf{d},\bs{\tau},\eps}(\bs{\phi}), \\
& \tilde{P}^i_{\mf{d},\bs{\tau},\eps}  (\bs{\phi}) := \mu_i f(\mu_i^{-\frac{1}{p-1}} P_\eps U_{\d_i,\xi_i}+\phi_i)-\mu_i f(\mu_i^{-\frac{1}{p-1}} P_\eps U_{\d_i,\xi_i}) - p(P_\eps U_{\d_i,\xi_i})^{p-1} \phi_i,
\end{align*}
and
\begin{align*}
\tilde{Q}^i_{\mf{d},\bs{\tau},\eps} (\bs{\phi}) & :=   \sum_{j \neq i} \beta_{ij} |\mu_j^{-\frac{1}{p-1}} P_\eps U_{\d_j,\x_j} + \phi_j|^{\frac{p+1}{2}} |\mu_i^{-\frac{1}{p-1}} P_\eps U_{\d_i,\x_i} + \phi_i|^{\frac{p-3}{2}} (\mu_i^{-\frac{1}{p-1}} P_\eps U_{\d_i,\x_i} + \phi_i) \\
& - \sum_{j \neq i} \beta_{ij} (\mu_j^{-\frac{1}{p-1}} P_\eps U_{\d_j,\x_j})^{\frac{p+1}{2}}(\mu_i^{-\frac{1}{p-1}} P_\eps U_{\d_i,\x_i})^{\frac{p-1}{2}} \\
& - \frac{p-1}{2} \sum_{j \neq i} \beta_{ij} (\mu_j^{-\frac{1}{p-1}} P_\eps U_{\d_j,\xi_j})^{\frac{p+1}{2}}   
(\mu_i^{-\frac{1}{p-1}} P_\eps U_{\d_i,\xi_i})^{\frac{p-3}{2}}  \phi_i \\
& - \frac{p+1}{2} \sum_{j \neq i} \beta_{ij} (\mu_j^{-\frac{1}{p-1}} P_\eps U_{\d_j,\xi_j})^{\frac{p-1}{2}}   
(\mu_i^{-\frac{1}{p-1}} P_\eps U_{\d_i,\xi_i})^{\frac{p-1}{2}}  \phi_j,
\end{align*}
with $i=1,\dots,m$. In Appendix \ref{app estimate nonlinear}, we prove that
\begin{equation}\label{estimate nonlinear}
\begin{split}
& |\tilde{R}^i_{\mf{d},\bs{\tau},\eps}|_{L^\frac{2N}{N+2}(\Omega_\eps)} \le C \eps^\frac{N-2}2 \\
& |\tilde{N}^i_{\mf{d},\bs{\tau},\eps}(\bs{\phi})|_{L^\frac{2N}{N+2}(\Omega_\eps)} \le C\|\bs{\phi}\|^2_{H_0^1(\Omega_\eps)}
%& |\tilde{P}^i_{\mf{d},\bs{\tau},\eps}(\bs{\phi})|_{L^\frac{2N}{N+2}(\Omega_\eps)} \le C\|\phi_i\|^2_{H_0^1(\Omega_\eps)} \\
%& |\tilde{Q}^i_{\mf{d},\bs{\tau},\eps}(\bs{\phi})|_{L^\frac{2N}{N+2}(\Omega_\eps)} \le C\|\bs{\phi}\|^2_{H_0^1(\Omega_\eps)}. 
\end{split}
\end{equation}
Therefore, by equation \eqref{473} there exist $C_1,C_2>0$ such that 
\begin{align*}
\| \mf{T}_{\mf{d},\bs{\tau},\eps}(\bs{\phi})\|_{H_0^1(\Omega_\eps)} & \le  C_1 \eps^\frac{N-2}2 + C_2  \|\bs{\phi}\|_{H_0^1(\Omega_\eps)}^2 
\end{align*}
Let $\bar C>C_1$ arbitrarily chosen, and let 
\[
Y_{\eps}:=  \left\{ \bs{\phi} \in H_0^1(\Omega_\eps,\R^m): \|\bs{\phi}\|_{H_0^1(\Omega_\eps)} \le \bar C \eps^\frac{N-2}2\right\}.
\]
Then there exists $\eps_1 \in (0,\eps_0]$ sufficiently small such that 
\[
\| \mf{T}_{\mf{d},\bs{\tau},\eps}(\bs{\phi})\|_{H_0^1(\Omega_\eps)}  \le C_1 \eps^\frac{N-2}2 + C_2 {\bar{C}}^2 \eps^{N-2} \le \bar C \eps^\frac{N-2}2
\]
for every $\eps \in (0,\eps_1)$ and $\bs{\phi} \in H_0^1(\Omega_\eps,\R^m)$, that is, $\mf{T}_{\mf{d},\bs{\tau},\eps}: Y_{\eps} \to Y_{\eps}$.

\noindent \textbf{Step 2)} $\mf{T}_{\mf{d},\bs{\tau},\eps}$ is a contraction in $Y_\eps$.

Notice that
\begin{equation}\label{673}
\| \mf{T}_{\mf{d},\bs{\tau},\eps}(\bs{\phi}^1)- \mf{T}_{\mf{d},\bs{\tau},\eps}(\bs{\phi}^2)\|_{H_0^1(\Omega_\eps)}  \le C | \tilde{\mf{N}}_{\mf{d},\bs{\tau},\eps}(\bs{\phi}^1)- \tilde{\mf{N}}_{\mf{d},\bs{\tau},\eps}(\bs{\phi}^2)|_{L^{\frac{2N}{N+2}}(\Omega_\eps)}.
\end{equation}
Recalling that $\tilde N^i_{\mf{d},\bs{\tau},\eps}(\bs{\phi}) =\tilde P^i_{\mf{d},\bs{\tau},\eps}(\bs{\phi}) + \tilde Q^i_{\mf{d},\bs{\tau},\eps}(\bs{\phi})$, we compute with a Taylor expansion (see Lemma \ref{lem a.3})
\[
| \tilde P^i_{\mf{d},\bs{\tau},\eps}(\bs{\phi}^1)- \tilde P^i_{\mf{d},\bs{\tau},\eps}(\bs{\phi}^2)|_{L^{\frac{2N}{N+2}}(\Omega_\eps)} \\
\le C \left| (P_\eps U_{\d_i,\x_i} + |\phi_i^1|+ |\phi_i^2|)^{p-2} ( |\phi_i^1|+ |\phi_i^2|)  |\phi_i^1-\phi_i^2| \right|_{L^{\frac{2N}{N+2}}(\Omega_\eps)}.
\]
Therefore, by the H\"older and the Sobolev inequalities
\begin{equation}\label{571}
\begin{split}
| \tilde P^i_{\mf{d},\bs{\tau},\eps}(\bs{\phi}^1) - \tilde P^i_{\mf{d},\bs{\tau},\eps}(\bs{\phi}^2)|_{L^{\frac{2N}{N+2}}(\Omega_\eps)}  &\le  ( C  + |\phi_i^1|_{2^*}^{p-2}+ |\phi_i^2|_{2^*}^{p-2}) \left( |\phi_i^1|_{2^*}+ |\phi_i^2|_{2^*}\right) |\phi_i^1-\phi_i^2|_{2^*} \\
& \le C \|\bs{\phi}\|_{H_0^1(\Omega_\eps)} \|\phi_i^1-\phi_i^2\|_{H_0^1(\Omega_\eps)} \le C\eps^\frac{N-2}2 \|\phi_i^1-\phi_i^2\|_{H_0^1(\Omega_\eps)} 
\end{split}
\end{equation}
for every $\bs{\phi}^1, \bs{\phi}^2 \in Y_\eps$, $i=1,\dots,m$.
Regarding $\tilde Q^i_{\bs{\d},\bs{\tau},\eps}$, in Appendix \ref{app estimate nonlinear} we show that
\begin{equation}\label{674}
 |  \tilde Q^i_{\mf{d},\bs{\tau},\eps}(\bs{\phi}^1) - \tilde Q^i_{\mf{d},\bs{\tau},\eps}(\bs{\phi}^2)|_{L^{\frac{2N}{N+2}}(\Omega_\eps)}  \le C \eps^\frac{N-2}2 \|\bs{\phi}^1-\bs{\phi}^2\|_{H_0^1(\Omega_\eps)}
\end{equation} 
or every $\bs{\phi}^1, \bs{\phi}^2 \in Y_\eps$, $i=1,\dots,m$. Collecting together \eqref{673}, \eqref{571} and \eqref{674}, we deduce that there exists $C>0$ such that
\[
\| \mf{T}_{\mf{d},\bs{\tau},\eps}(\bs{\phi}^1)- \mf{T}_{\mf{d},\bs{\tau},\eps}(\bs{\phi}^2)\|_{H_0^1(\Omega_\eps)}  \le C\eps^\frac{N-2}2 \|\bs{\phi}^1-\bs{\phi}^2\|_{H_0^1(\Omega_\eps)} \qquad \forall \bs{\phi}^1,\bs{\phi}^2 \in Y_{\eps}.
\]
By replacing the necessary $\eps_1$ with a smaller quantity, we see that for $\eps \in (0,\eps_1)$ the map $\mf{T}_{\mf{d},\bs{\tau},\eps}$ is a contraction in $Y_\eps$, and hence the thesis follows by the contraction mapping theorem.
\end{proof}

Lemma \ref{lem: exist nonlin} enables us to define a map
\[
A:(0,\eps_1) \times X_\eta \to \mf{K}_{\mf{d},\bs{\tau},\eps}^\perp, \quad (\eps, \mf{d},\bs{\tau}) \mapsto \bs{\phi}^{\mf{d},\bs{\tau},\eps}.
\]
To complete the proof of Proposition \ref{prop: pb K ort}, it remains to check that this map if differentiable, and to prove the desired estimate on the derivative.

\begin{lemma}\label{lem: diff map}
There exists $\eps_2>0$ small enough such that the map $A$ is of class $\mathcal{C}^1$ in $(0,\eps_2) \times X_\eta$.
\end{lemma}

\begin{proof}
We apply the implicit function theorem to $\mf{T}: (0,\eps_1) \times X_\eta \times \mf{K}_{\mf{d},\bs{\tau},\eps}^\perp \to \mf{K}_{\mf{d},\bs{\tau},\eps}^\perp$ defined by
\[
\mf{T}(\eps, \mf{d},\bs{\tau}, \bs{\phi}) = \mf{L}_{\mf{d},\bs{\tau},\eps}(\bs{\phi}) - \mf{R}_{\mf{d},\bs{\tau},\eps} - \mf{N}_{\mf{d},\bs{\tau},\eps}(\bs{\phi}).
\]
By Lemma \ref{lem: exist nonlin}, we know that $\mf{T}(\eps, \mf{d},\bs{\tau}, \bs{\phi}^{\mf{d},\bs{\tau},\eps}) = 0$. We shall prove that $D_{\bs{\phi}}\mf{T}(\eps, \mf{d},\bs{\tau}, \bs{\phi}^{\mf{d},\bs{\tau},\eps})$ is invertible. To this aim, by the Fredholm alternative, it is sufficient to check that $D_{\bs{\phi}}\mf{T}(\eps, \mf{d},\bs{\tau}, \bs{\phi}^{\mf{d},\bs{\tau},\eps})$ is injective, since $D_{\bs{\phi}}\mf{T}(\eps, \mf{d},\bs{\tau}, \bs{\phi}^{\mf{d},\bs{\tau},\eps})$ is a compact perturbation of the identity (due to the compactness of $i^*$). 

In the rest of the proof we often write $P$, $U_i$, $\phi_i$ instead of $P_\eps$, $U_{\d_i,\x_i}$, $\phi_i^{\mf{d},\bs{\tau},\eps}$, to ease the notation.

Notice that
\[
D_{\bs{\phi}}\mf{T}(\eps, \mf{d},\bs{\tau}, \bs{\phi}^{\mf{d},\bs{\tau},\eps})[\bs{\psi}] = \mf{L}_{\mf{d},\bs{\tau},\eps}(\bs{\psi}) - D_{\bs{\phi}}\mf{N}_{\mf{d},\bs{\tau},\eps}(\bs{\phi}^{\mf{d},\bs{\tau},\eps})[\bs{\psi}].
\]
By definition of $\mf{N}_{\mf{d},\bs{\tau},\eps}$ (see \eqref{def N})
\begin{align*}
D_{\bs{\phi}} &N_{\mf{d},\bs{\tau},\eps}^i (\bs{\phi}^{\mf{d},\bs{\tau},\eps})[\bs{\psi}] = \Pi_i^{\perp} \circ i^* \Bigg[ \mu_i \left( f'(\mu_i^{-\frac{1}{p-1}} P U_i+\phi_i)- f'(\mu_i^{-\frac{1}{p-1}} P U_i)  \right)\psi_i  \\
& + \frac{p+1}2 \sum_{j \neq i} \beta_{ij}  \bigg( |\mu_i^{-\frac{1}{p-1}} P U_i + \phi_i|^{\frac{p-3}{2}} |\mu_j^{-\frac{1}{p-1}} P U_j + \phi_j|^{\frac{p-3}{2}} 
% \hphantom{+ \frac{p+1}2 \sum_{j \neq i}} \cdot 
(\mu_i^{-\frac{1}{p-1}} P U_i + \phi_i) (\mu_j^{-\frac{1}{p-1}} P U_j + \phi_j)  \\
& \hphantom{+ \frac{p+1}2 \sum_{j \neq i} \beta_{ij}} \ \ -  (\mu_i^{-\frac{1}{p-1}} P U_i)^{\frac{p-1}{2}} (\mu_j^{-\frac{1}{p-1}} P U_j)^{\frac{p-1}{2}}   
  \bigg) \psi_j \\
&   + \frac{p-1}2 \sum_{j \neq i} \beta_{ij} \bigg(  |\mu_i^{-\frac{1}{p-1}} P U_i + \phi_i|^{\frac{p-3}{2}} |\mu_j^{-\frac{1}{p-1}} P U_j + \phi_j|^{\frac{p+1}{2}} \\ 
& \hphantom{+ \frac{p+1}2 \sum_{j \neq i} \beta_{ij}} \ \ -(\mu_i^{-\frac{1}{p-1}} P U_i)^{\frac{p-3}{2}}  (\mu_j^{-\frac{1}{p-1}} P U_j)^{\frac{p+1}{2}}   \bigg)
 \psi_i \Bigg].
\end{align*}
Then, if $N=3$ (i.e. $p=5$), by the Lagrange theorem and using the fact that $0 \le P_\eps U_{\d_i,\x_i} \le U_{\d_i,\x_i}$,
\begin{align*}
\|D_{\bs{\phi}} &N_{\mf{d},\bs{\tau},\eps}^i (\bs{\phi}^{\mf{d},\bs{\tau},\eps})[\bs{\psi}]\|_{H_0^1(\Omega_\eps)}  \le C \Bigg[ | U_i^{3} \phi_i \psi_i|_{\frac{6}{5}} 
 + \big| |\phi_i|^4 \psi_i\big|_{\frac{6}{5}} \\
 & + \sum_{j \neq i} \bigg( \big| | U_i + |\phi_i| |  \, |U_j + |\phi_j| |^2\,|\phi_i| \, |\psi_j|\big|_{\frac65} + \big| | U_i + |\phi_i| |^2  \, |U_j + |\phi_j| |\,|\phi_j|\, |\psi_j| \big|_{\frac65}  \\
 & \hphantom{+ \sum_{j \neq i} \bigg(}\ + \big| |\phi_i|  \, |U_j + |\phi_j| |^3 \, |\psi_i|\big|_{\frac65} + \big| | U_i + |\phi_i| |  \, |U_j + |\phi_j| |^2\,|\phi_j|\, |\psi_i| \big|_{\frac65} \bigg)\Bigg];
% &\hphantom{\le C \Bigg[} + \sum_{j \neq i} \bigg( | U_{\d_i,\x_i}^\frac{p-1}2 U_{\d_j,\x_j}^{\frac{p-1}2} \psi_j|_{\frac{2N}{N+2}} +    \big| U_{\d_i,\x_i}^\frac{p-1}2 |\phi_j^{\mf{d},\bs{\tau},\eps}|^{\frac{p-1}2} \psi_j\big|_{\frac{2N}{N+2}}\\
% & \hphantom{\le C \Bigg[+ \sum_{j \neq i} \bigg(} + \big| |\phi_i^{\mf{d},\bs{\tau},\eps} |^{\frac{p-1}2} U_{\d_j,\x_j}^\frac{p-1}2  \psi_j\big|_{\frac{2N}{N+2}} + \big| |\phi_i^{\mf{d},\bs{\tau},\eps} |^{\frac{p-1}2} |\phi_i^{\mf{d},\bs{\tau},\eps}|^\frac{p-1}2  \psi_j\big|_{\frac{2N}{N+2}} \\
% &  \hphantom{\le C \Bigg[+ \sum_{j \neq i} \bigg(} +  | U_{\d_i,\x_i}^\frac{p-3}2 U_{\d_j,\x_j}^{\frac{p+1}2} \psi_i|_{\frac{2N}{N+2}} +   \big| U_{\d_i,\x_i}^\frac{p-3}2 |\phi_j^{\mf{d},\bs{\tau},\eps}|^{\frac{p+1}2} \psi_i \\
% &  \hphantom{\le C \Bigg[+ \sum_{j \neq i} \bigg(}  +  \big| |\phi_i^{\mf{d},\bs{\tau},\eps} |^{\frac{p-3}2} U_{\d_j,\x_j}^\frac{p+1}2  \psi_j\big|_{\frac{2N}{N+2}} + \big| |\phi_i^{\mf{d},\bs{\tau},\eps} |^{\frac{p-3}2} |\phi_i^{\mf{d},\bs{\tau},\eps}|^\frac{p+1}2  \psi_j\big|_{\frac{2N}{N+2}} \bigg) \Bigg].
\end{align*}
similarly, if $N=4$, using the fact that $p-3=0$, we find
\begin{align*}
\|D_{\bs{\phi}} &N_{\mf{d},\bs{\tau},\eps}^i (\bs{\phi}^{\mf{d},\bs{\tau},\eps})[\bs{\psi}]\|_{H_0^1(\Omega_\eps)}  \le C \Bigg[ | U_i \phi_i \psi_i|_{\frac{4}{3}} 
 + \big| |\phi_i|^2 \psi_i\big|_{\frac{4}{3}} \\
 & + \sum_{j \neq i} \bigg( |  U_i \phi_j \psi_j |_{\frac43} + |\phi_i U_j \psi_j|_\frac43 + |\phi_i \phi_j \psi_j|_\frac43  + | U_j \phi_j \psi_i|_\frac43 + |\phi_j^2 \psi_i|_\frac43 \bigg)\Bigg].
 \end{align*}
In order to estimate the right hand side, it is not difficult to apply the H\"older and the Sobolev inequalities, as well as the estimate in Lemma \ref{lem: exist nonlin}, to deduce that
\begin{equation}\label{771}
\|D_{\bs{\phi}} N_{\mf{d},\bs{\tau},\eps}^i (\bs{\phi}^{\mf{d},\bs{\tau},\eps})[\bs{\psi}]\|_{H_0^1(\Omega_\eps)} = o(1) \|\bs{\psi}\|_{H_0^1(\Omega_\eps)}, 
\end{equation}
where $o(1) \to 0$ as $\eps \to 0$ (in particular, we use the fact that $|U_i|_{2^*} \le C$, and $\|\phi_i\| \to 0$ as $\eps \to 0$).

Using \eqref{771} and \eqref{17ott1}, we infer that if $D_{\bs{\phi}}\mf{T}(\eps, \mf{d},\bs{\tau}, \bs{\phi}^{\mf{d},\bs{\tau},\eps})[\bs{\psi}] =0$, that is
\[ 
\mf{L}_{\mf{d},\bs{\tau},\eps}(\bs{\psi}) = D_{\bs{\phi}}\mf{N}_{\mf{d},\bs{\tau},\eps}(\bs{\phi}^{\mf{d},\bs{\tau},\eps})[\bs{\psi}],
\]
then
\[
C \|\bs{\psi}\|_{H_0^1(\Omega_\eps)} \le o(1) \|\bs{\psi}\|_{H_0^1(\Omega_\eps)},
\]
which finally implies that $\bs{\psi}=0$. This means that $D_{\bs{\phi}}\mf{T}(\eps, \mf{d},\bs{\tau}, \bs{\phi}^{\mf{d},\bs{\tau},\eps})$ is injective for $\eps$ small enough and, as observed, this suffices to complete the proof.
\end{proof}

The following lemma completes the proof of Proposition \ref{prop: pb K ort}.

\begin{lemma}\label{lem: der error}
There exists $\bar \eps>0$ small enough and a constant $C>0$ such that
\[
\| \nabla_{(\bs{\d},\bs{\xi})} \bs{\phi}^{\mf{d},\bs{\tau},\eps} \|_{H_0^1(\Omega_\eps)} \le C \eps^{\frac{N-3}{2}}
\quad \iff \quad 
\| \nabla_{(\mf{d},\bs{\tau})} \bs{\phi}^{\mf{d},\bs{\tau},\eps} \|_{H_0^1(\Omega_\eps)} \le C \eps^{\frac{N-2}{2}}
\]
for every $\eps \in (0,\bar \eps)$ and $(\mf{d},\bs{\tau}) \in X_\eta$.
\end{lemma}

\begin{proof}
The equivalence of the two inequalities follows by the chain rule and by the ansatz \eqref{asympt expansion}. Let $s_{i,h}=\d_i$ if $h=0$, and $s_{i,h} = \xi_{i,h}$ if $h=1,\dots,N$, with $i=1,\dots,m$. We differentiate the equation
\begin{equation}\label{eq per phi 2}
\mf{L}_{\mf{d},\bs{\tau},\eps}(\bs{\phi}^{\mf{d},\bs{\tau},\eps}) = \mf{R}_{\mf{d},\bs{\tau},\eps} - \mf{N}_{\mf{d},\bs{\tau},\eps}(\bs{\phi}^{\mf{d},\bs{\tau},\eps}),
\end{equation}
with respect to a variable $s_{i,h}$, and we obtain
\begin{multline*}
\pa_{s_{i,h}} \mf{L}_{\mf{d},\bs{\tau},\eps}(\bs{\phi}^{\mf{d},\bs{\tau},\eps}) + \mf{L}_{\mf{d},\bs{\tau},\eps}(\pa_{s_{i,h}}\bs{\phi}^{\mf{d},\bs{\tau},\eps})\\
 = \pa_{s_{i,h}} \mf{R}_{\mf{d},\bs{\tau},\eps} - \pa_{s_{i,h}} \mf{N}_{\mf{d},\bs{\tau},\eps} (\bs{\phi}^{\mf{d},\bs{\tau},\eps}) - D_{\bs{\phi}} \mf{N}_{\mf{d},\bs{\tau},\eps} (\bs{\phi}^{\mf{d},\bs{\tau},\eps}) [\pa_{s_{i,h}} \bs{\phi}^{\mf{d},\bs{\tau},\eps}].
\end{multline*}
We claim that
\begin{equation}\label{stima derivative error}
\| \pa_{s_{i,h}} \mf{L}_{\mf{d},\bs{\tau},\eps}(\bs{\phi}^{\mf{d},\bs{\tau},\eps}) \|_{H_0^1(\Omega_\eps)} + \| \pa_{s_{i,h}} \mf{R}_{\mf{d},\bs{\tau},\eps} \|_{H_0^1(\Omega_\eps)}+ \| \pa_{s_{i,h}} \mf{N}_{\mf{d},\bs{\tau},\eps}(\bs{\phi}^{\mf{d},\bs{\tau},\eps}) \|_{H_0^1(\Omega_\eps)}  \le C \eps^\frac{N-3}2.
\end{equation}
With this estimate in our hands, the thesis can be easily proved. Indeed, by \eqref{771} and the inequality in Lemma \ref{lem: linear part}, we deduce that 
\begin{align*}
C \| \pa_{s_{i,h}}\bs{\phi}^{\mf{d},\bs{\tau},\eps}\| &\le \| \mf{L}_{\mf{d},\bs{\tau},\eps}(\pa_{s_{i,h}}\bs{\phi}^{\mf{d},\bs{\tau},\eps}) - D_{\bs{\phi}} \mf{N}_{\mf{d},\bs{\tau},\eps} (\bs{\phi}^{\mf{d},\bs{\tau},\eps}) [\pa_{s_{i,h}} \bs{\phi}^{\mf{d},\bs{\tau},\eps}]\| \\
& \le \|\pa_{s_{i,h}} \mf{L}_{\mf{d},\bs{\tau},\eps}(\bs{\phi}^{\mf{d},\bs{\tau},\eps})\| + \|  \pa_{s_{i,h}} \mf{R}_{\mf{d},\bs{\tau},\eps}\| + \|  \pa_{s_{i,h}} \mf{N}_{\mf{d},\bs{\tau},\eps} (\bs{\phi}^{\mf{d},\bs{\tau},\eps})\| \le C \eps^\frac{N-3}2.
\end{align*}

The validity of \eqref{stima derivative error} can be checked by direct computations, and the details are presented in Appendix \ref{app derivative error}.
\end{proof}

\section{The reduced problem}\label{sec: reduced problem}

In this section we solve equation \eqref{pb K} with $\bs{\phi}= \bs{\phi}^{\mf{d},\bs{\tau},\eps}$.

Let $J_\eps: H_0^1(\Omega_\eps,\R^m) \to \R$ be defined by
\begin{equation}\label{def action}
J_\eps(u_1,\dots,u_m) = \int_{\Omega_\eps} \frac12 \sum_{i=1}^m |\nabla u_i|^2 - \mu_i F(u_i) - \frac2{p+1}\sum_{ 1 \le i < j \le m} \beta_{ij} |u_i|^{\frac{p+1}2} |u_j|^\frac{p+1}2,
\end{equation}
where $F: \R \to \R$, $t \mapsto (t^+)^{p+1}/(p+1)$ is the primitive of $f$. Critical points of $J_\eps$ are solution to \eqref{system 2}, and hence solutions to \eqref{system} (here we use the fact that $\beta_{ij}=\beta_{ji}$).

Let us introduce the reduced functional $\tilde J_\eps: X_\eta \to \R$, 
\[
\tilde J_\eps(\mf{d},\bs{\tau}):= J_\eps\left( \mu_1^{-\frac1{p-1}} P_\eps U_{\d_1,\x_1} + \phi_1^{\mf{d},\bs{\tau},\eps}, \dots,  \mu_m^{-\frac1{p-1}} P_\eps U_{\d_m,\x_m} + \phi_m^{\mf{d},\bs{\tau},\eps}\right).
\] 
In order to simplify the notation, from now on we often write 
\[
V_i^{\mf{d},\bs{\tau},\eps}:= \mu_i^{-\frac1{p-1}} P_\eps U_{\d_i,\x_i} + \phi_i^{\mf{d},\bs{\tau},\eps}.
\]

\begin{lemma}\label{lem: the reduction works}
There exists $\bar \eps>0$ sufficiently small such that if $(\mf{d},\bs{\tau})$ is a critical point of $\tilde J_\eps$, and $\eps \in (0, \bar \eps)$, then 
\[
\left( \mu_1^{-\frac1{p-1}} P_\eps U_{\d_1,\x_1} + \phi_1^{\mf{d},\bs{\tau},\eps}, \dots,  \mu_m^{-\frac1{p-1}} P_\eps U_{\d_m,\x_m} + \phi_m^{\mf{d},\bs{\tau},\eps}\right)
\]
is a solution to \eqref{pb K}, and hence a solution to \eqref{system 2}.
\end{lemma}

\begin{proof}
We start observing that, denoting by $\langle \cdot, \cdot \rangle$ the scalar product in $H_0^1(\Omega_\eps)$, we have
\begin{align*}
\frac{\pa \tilde J_\eps}{\pa d_i}&(\mf{d},\bs{\tau})  =\sqrt{\eps} \frac{\pa \tilde J_\eps}{\pa \d_i}(\mf{d},\bs{\tau}) \\
& = \sqrt{\eps} d J_\eps\left (V_1^{\mf{d},\bs{\tau},\eps}, \dots, V_m^{\mf{d},\bs{\tau},\eps} \right)\left[\left( \frac{\pa \phi_1^{\mf{d},\bs{\tau},\eps}}{\pa \delta_i}, \dots, \frac{\pa \phi_m^{\mf{d},\bs{\tau},\eps}}{\pa \delta_i}\right)\right] \\
& + \sqrt{\eps} \pa_i J_\eps\left (V_1^{\mf{d},\bs{\tau},\eps}, \dots, V_m^{\mf{d},\bs{\tau},\eps} \right)\left[ \mu_i^{-\frac{1}{p-1}} P_\eps \psi^0_{\delta_i,\xi_i} \right] \\
& = \sqrt{\eps} \sum_{k =1}^m \left\langle V_k^{\mf{d},\bs{\tau},\eps} -  i^*\left[ \mu_k f(V_k^{\mf{d},\bs{\tau},\eps})  + \sum_{j \neq k} \beta_{kj} |V_j^{\mf{d},\bs{\tau},\eps}|^{\frac{p+1}{2}} |V_k^{\mf{d},\bs{\tau},\eps}|^{\frac{p-3}{2}} V_k^{\mf{d},\bs{\tau},\eps}\right], \frac{\pa \phi_k^{\mf{d},\bs{\tau},\eps}}{\pa \delta_i}\right\rangle \\
& + \sqrt{\eps} \left\langle V_i^{\mf{d},\bs{\tau},\eps} -  i^*\left[ \mu_i f(V_i^{\mf{d},\bs{\tau},\eps})  + \sum_{j \neq i} \beta_{ij} |V_j^{\mf{d},\bs{\tau},\eps}|^{\frac{p+1}{2}} |V_i^{\mf{d},\bs{\tau},\eps}|^{\frac{p-3}{2}} V_i^{\mf{d},\bs{\tau},\eps}\right], \mu_i^{-\frac{1}{p-1}} P_\eps \psi^0_{\delta_i,\xi_i}\right\rangle
\end{align*}
By Proposition \ref{prop: pb K ort}, we know that the projection of $\nabla J_\eps(V_1^{\mf{d},\bs{\tau},\eps},\dots,V_m^{\mf{d},\bs{\tau},\eps})$ on $\bs{K}^\perp_{\mf{d},\bs{\tau},\eps}$ is $0$. This means that the terms in left position inside the brackets are linear combination of the partial derivatives $P_\eps \psi^\ell_{\delta_k,\xi_k}$, and hence the previous chain of equalities can be continued in the following way:
\begin{align*}
\frac{\pa \tilde J_\eps}{\pa d_i}&(\mf{d},\bs{\tau})  = \sqrt{\eps} \sum_{k=1}^m \left\langle \sum_{\ell=0}^N c_{k,\ell}^\eps P_\eps \psi^{\ell}_{\d_k,\x_k},  \frac{\pa \phi_k^{\mf{d},\bs{\tau},\eps}}{\pa \delta_i}\right\rangle + \sqrt{\eps}\left\langle \sum_{\ell=0}^N c_{i,\ell}^\eps P_\eps \psi^{\ell}_{\d_k,\x_k}, \mu_i^{-\frac{1}{p-1}} P_\eps \psi^0_{\delta_i,\xi_i}\right\rangle.
\end{align*}
Now, let $(\mf{d},\bs{\tau})$ be a critical point for $\tilde J_\eps$. We have then (multiplying by $\sqrt{\eps}$)
\begin{equation}\label{481}
\eps \sum_{k=1}^m \left\langle \sum_{\ell=0}^N c_{k,\ell}^\eps P_\eps \psi^{\ell}_{\d_k,\x_k},  \frac{\pa \phi_k^{\mf{d},\bs{\tau},\eps}}{\pa \delta_i}\right\rangle + \eps \left\langle \sum_{\ell=0}^N c_{i,\ell}^\eps P_\eps \psi^{\ell}_{\d_k,\x_k}, \mu_i^{-\frac{1}{p-1}} P_\eps \psi^0_{\delta_i,\xi_i}\right\rangle = 0.
\end{equation}
In the same way, if we compute the derivatives with respect to $\tau_{i,h}$ ($h=1,\dots,m$) and we evaluate them in a critical point, we obtain
\begin{equation}\label{482}
\eps \sum_{k=1}^m \left\langle \sum_{\ell=0}^N c_{k,\ell}^\eps P_\eps \psi^{\ell}_{\d_k,\x_k},  \frac{\pa \phi_k^{\mf{d},\bs{\tau},\eps}}{\pa \xi_{i,h}}\right\rangle + \eps \left\langle \sum_{\ell=0}^N c_{i,\ell}^\eps P_\eps \psi^{\ell}_{\d_k,\x_k}, \mu_i^{-\frac{1}{p-1}} P_\eps \psi^h_{\delta_i,\xi_i}\right\rangle = 0.
\end{equation}
Letting $i$ and $h$ vary, \eqref{481} and \eqref{482} provides us a linear homogeneous system of $m(N+1)$ equations in the $m(N+1)$ unknowns $c_{i,h}^\eps$. We aim at showing that the system has only the trivial solution; this means that also the projection on $\nabla J_\eps(V_1^{\mf{d},\bs{\tau},\eps},\dots,V_m^{\mf{d},\bs{\tau},\eps})$ on $\mf{K}_{\mf{d},\bs{\tau},\eps}$ vanishes, i.e. $\nabla J_\eps(V_1^{\mf{d},\bs{\tau},\eps},\dots,V_m^{\mf{d},\bs{\tau},\eps}) = 0$ in $H_0^1(\Omega_\eps)$, and completes the proof. Thus, we consider now \eqref{481} and \eqref{482}, and we show that the matrix of the coefficients is invertible for $\eps$ small enough.

The last term in both \eqref{481} and \eqref{482} can be estimated as in step 1 in Lemma \ref{lem: linear part}: recalling that $\eps \simeq \delta_i^2$ as $\eps \to 0$, we have for $h=0,\dots, N$
\begin{equation}\label{480}
 \eps \left\langle \sum_{\ell=0}^N c_{i,\ell}^\eps P_\eps \psi^{\ell}_{\d_k,\x_k}, \mu_i^{-\frac{1}{p-1}} P_\eps \psi^h_{\delta_i,\xi_i}\right\rangle = \mu_i^{-\frac{1}{p-1}} \sigma_{hh} c_{i,h}^\eps + o(1) \sum_{\ell \neq h} c_{i,\ell}^\eps , 
 \end{equation}
as $\eps \to 0$, with $\sigma_{hh} >0$ defined by \eqref{sigma lk}. Let us now consider the terms involving the derivative of $\phi_k^{\mf{d},\bs{\tau},\eps}$. By Lemma \ref{lem: der error}, the Cauchy-Schwarz inequality, and recalling again that $\d_i^2 \|P_\eps \psi^{\ell}_{\d_i,\x_i}\|^2_{H_0^1(\Omega_\eps)} \to \sigma_{\ell \ell}$, we deduce that 
\begin{equation}\label{486}
\eps \left| \left\langle P_\eps \psi^{\ell}_{\d_i,\x_i},  \frac{\pa \phi_i^{\mf{d},\bs{\tau},\eps}}{\pa \xi_{i,h}}\right\rangle \right| \le C \delta_i^2  \left\| P_\eps \psi^{\ell}_{\d_i,\x_i}\right\|_{H_0^1(\Omega_\eps)} \left\| \frac{\pa \phi_i^{\mf{d},\bs{\tau},\eps}}{\pa \xi_{i,h}} \right\|_{H_0^1(\Omega_\eps)} \le C \d_i \eps^\frac{N-3}2 \le C \eps^\frac{N-2}2,
\end{equation}
where we used also the fact that $\d_i \simeq \eps^\frac12$. The same estimate holds if we consider the derivatives with respect to $\d_{i}$, $i=1,\dots,m$.
Thus, plugging \eqref{480} and \eqref{486} inside \eqref{481} and \eqref{482}, we infer that the matrix associated to the system is strictly diagonally dominant, hence invertible, and being homogeneous has only the trivial solution. As observed, this completes the proof.
\end{proof}

Due to Lemma \ref{lem: the reduction works}, in order to complete the proof of Theorem \ref{thm: main 1} we have to find a critical point of the reduced functional $\tilde J_\eps$ in $X_\eta$. In this perspective, we need the asymptotic expansion of $J_\eps(V_1^{\mf{d},\bs{\tau},\eps},\dots,V_m^{\mf{d},\bs{\tau},\eps})$ as $\eps \to 0^+$. We introduce the real numbers 
\begin{equation}\label{def b_i}
\begin{split}
b_1 := \frac{\alpha_N^{p+1} }{N}  \int_{\R^N} \frac{dy}{(1+|y^2)^N}, \qquad  
b_2  := \frac{\alpha_N^{p+1}}2   \int_{\R^N}  \frac{dy}{(1+|y|^2)^\frac{N+2}2},
\end{split}
\end{equation}
and the functions $\Gamma: \R^N \to \R$, $\Psi: X_\eta \to \R$ defined by
\begin{equation}\label{def gamma e psi}
\begin{split}
\Gamma(x) & := \int_{\R^N}   \frac{dy}{|y+x|^{N-2} (1+|y|^2)^{\frac{N+2}2}} \\
\Psi(\mf{d},\bs{\tau}) & := \sum_{i=1}^m \mu_i^{-\frac2{p-1}} \left[ b_2 H(a_i,a_i) d_i^{N-2}+ \frac{\alpha_N^{p+1} r_i^{N-2}}2    \frac{ \Gamma(\tau_i)}{d_i^{N-2}(1+|\tau_i|^2)^{\frac{N-2}2}}\right].
\end{split}
\end{equation} 

\begin{proposition}\label{prop: exp func}
We have 
\[
\tilde J_\eps(\mf{d},\bs{\tau}) = \left(\sum_{i=1}^m \mu_i^{-\frac2{p-1}}\right) b_1 + \Psi(\mf{d},\bs{\tau})\eps^\frac{N-2}2 + R_\eps(\mf{d},\bs{\tau}),
\]
with $R_\eps = o(\eps^\frac{N-2}2)$ $\C^1$-uniformly in $X_\eta$ as $\eps \to 0$.
\end{proposition}

Here and what follows, we write that $f_\eps = o(\eps^\alpha)$ $\C^k$-uniformly in $X_\eta$ as $\eps \to 0$ if
\[
\lim_{\eps \to 0} \frac{\|f_\eps\|_{\C^k(X_\eta)}}{\eps^\alpha} = 0
\]

The proof of the proposition takes most of the rest of the section. In order to keep the notation short, in what follows we sometimes write $P$, $U_i$, $V_i$, $\phi_i$ and $\|\cdot\|$ instead of $P_\eps$, $U_{\d_i,\x_i}$, $V_i^{\mf{d},\bs{\tau},\eps}$, $\phi_i^{\mf{d},\bs{\phi},\eps}$ and $\|\cdot \|_{H_0^1(\Omega_\eps)}$, respectively.

Let $I_{\eps,\mu}: H_0^1(\Omega_\eps) \to \R$ be defined by 
\[
I_\eps(u) = \frac{1}{2}\int_{\Omega_\eps} |\nabla u|^2 - \int_{\Omega_\eps} \mu F(u), \quad I_\eps:= I_{\eps,1}.
\]
Then, by the Lagrange theorem,
\begin{equation}\label{981}
\begin{split}
\tilde J_\eps(\mf{d},\bs{\tau}) &= \sum_i I_{\eps,\mu_i}(V_i) - \frac2{p+1}\sum_{i<j}\int_{\Omega_\eps} \beta_{ij} |V_i|^\frac{p+1}2 |V_j|^\frac{p+1}2
\\
& = \sum_i  \mu_i^{-\frac2{p+1}}  I_\eps(P U_i) + \frac12\sum_i \int_{\Omega_\eps} 2 \mu_i^{-\frac1{p-1}} \nabla (P U_i) \cdot \nabla \phi_i + |\nabla \phi_i|^2 \\
& - \sum_i \int_{\Omega_\eps} \mu_i  \left( F(\mu_i^{-\frac1{p-1}}  P U_i +  \phi_i)   - F(\mu_i^{-\frac1{p-1}} P U_i)\right)
\\
&  - \frac{2}{p+1} \sum_{i<j} \beta_{ij} \int_{\Omega_\eps} |\mu_i^{-\frac{1}{p-1}} P U_i|^{\frac{p+1}2} |\mu_j^{-\frac{1}{p-1}} P U_j|^{\frac{p+1}2} 
\\
& - \sum_{i<j} \beta_{ij} \int_{\Omega_\eps} \left( |\mu_i^{-\frac{1}{p-1}} P U_i + \eta _i \phi_i|^{\frac{p+1}2} |\mu_j^{-\frac{1}{p-1}} P U_j + \eta _j \phi_j|^{\frac{p-1}2} \phi_j \right. \\
& \hphantom{- \sum_{i<j} \beta_{ij} \int_{\Omega_\eps} \Big( } \left.  + |\mu_i^{-\frac{1}{p-1}} P U_{i} + \eta_i \phi_i|^{\frac{p-1}2} |\mu_j^{-\frac{1}{p-1}} P U_j+ \eta_j \phi_j|^{\frac{p+1}2} \phi_i\right) 
\end{split}
\end{equation}
where $\eta_i, \eta_j$ denote continuous functions with values in $[0,1]$. 

We start considering the asymptotic expansion of $I_\eps(P_\eps U_{\d_i,\x_i})$.

\begin{lemma}\label{lem 5.3}
For every $i=1,\dots,m$, it results that
\[
I_\eps(P_\eps U_{\d_i,\x_i})  = b_1 + \left[ b_2 H(a_i,a_i) d_i^{N-2}+ \frac{\alpha_N^{p+1}   r_i^{N-2}}2    \frac{ \Gamma(\tau_i)}{d_i^{N-2}(1+|\tau_i|^2)^{\frac{N-2}2}} \right]\eps^\frac{N-2}2 + o(\eps^{\frac{N-2}2})
\]
as $\eps \to 0$, $\C^1$-uniformly in $(\mf{d},\bs{\tau}) \in X_\eta$, where the real numbers $b_1,b_2>0$ and the function $\Gamma: \R^N \to \R$ are defined by \eqref{def b_i} and \eqref{def gamma e psi}.
\end{lemma}

\begin{proof}
At first, using the definition of $P_\eps U_{\d_i,\x_i}$, we observe that
\begin{equation}\label{982}
\begin{split}
I_\eps (P_\eps U_{\d_i,\x_i}) & = \frac12 \int_{\Omega_\eps} P_\eps U_{\d_i,\x_i} U_{\d_i,\x_i}^p - \frac1{p+1}  \int_{\Omega_\eps} (P_\eps U_{\d_i,\x_i} )^{p+1} \\
& = \left(\frac12-\frac1{p+1}\right) \int_{\Omega_\eps} U_{\d_i,\x_i}^{p+1} + \frac12  \int_{\Omega_\eps} U_{\d_i,\x_i}^{p}( P_\eps U_{\d_i,\x_i}- U_{\d_i,\x_i}) \\
& \qquad \qquad -\frac1{p+1} \int_{\Omega_\eps} ((P_\eps U_{\d_i,\x_i})^{p+1}- U_{\d_i,\x_i}^{p+1}) \\
& = \frac1N \int_{\Omega_\eps} U_{\d_i,\x_i}^{p+1} -  \frac12 \int_{\Omega_\eps} U_{\d_i,\x_i}^{p}( P_\eps U_{\d_i,\x_i}- U_{\d_i,\x_i})  \\
& \qquad \qquad - \frac{p}2\int_{\Omega_\eps} (U_{\d_i,\x_i}+ \eta_i ( P_\eps U_{\d_i,\x_i}- U_{\d_i,\x_i}))^{p-1}( P_\eps U_{\d_i,\x_i}- U_{\d_i,\x_i})^2
\end{split}
\end{equation}
for a function $\eta_i$ with values in $[0,1]$. Now we estimate separately the three terms on the right hand side. Recalling that $\Omega_\eps= \Omega \setminus \bigcup_i B_{r_i \eps} (a_i)$, the first term gives
\begin{equation}\label{980}
\begin{split}
\int_{\Omega_\eps} U_{\d_i,\x_i}^{p+1} = \alpha_N^{p+1} \int_{\frac{\Omega_\eps-\xi_i}{\d_i}} \frac{dy}{(1+|y|^2)^N}  = b_1 + O\left( \left(\frac{\eps}{\d_i}\right)^{N}+\d_i^N\right).
\end{split}
\end{equation}
To treat the second term, we use Lemma \ref{lem a.1} ($R=R_{\eps,d_i,\tau_i}$):
\begin{equation}\label{983}
\begin{split}
\int_{\Omega_\eps} U_{\d_i,\x_i}^{p} & ( P_\eps U_{\d_i,\x_i}- U_{\d_i,\x_i}) = \alpha_N^{p} \int_{\Omega_\eps}   R_{\eps,\d_i,\x_i} U_{\d_i,\x_i}^p \\
&  -\alpha_N^{p+1} \int_{\Omega_\eps} \left( \d_i^\frac{N-2}2 H(x,\xi_i) + \frac{\d_i^{-\frac{N-2}2}}{(1+|\tau|^2)^{\frac{N-2}2}} \left( \frac{r_i \eps}{|x-a|}\right)^{N-2} \right) U_{\d_i,\x_i}^p.
\end{split}
\end{equation}
Now, by dominated convergence
\begin{equation}\label{984}
\begin{split}
\int_{\Omega_\eps} \d_i^\frac{N-2}2 H(x,\xi_i) U_{\d_i,\x_i}^p & = \int_{\frac{\Omega_\eps-\xi_i}{\d_i}} \d_i^{N-2} \frac{H(\xi_i+ \d_i y,\xi_i)}{(1+|y|^2)^{\frac{N+2}2}}\,dy \\
& = b_2 H(a_i,a_i) \d_i^{N-2}  + o(\d_i^{N-2}),
\end{split}
\end{equation}
and
\begin{equation}\label{985}
\begin{split}
\int_{\Omega_\eps} \frac{\d_i^{-\frac{N-2}2}}{(1+|\tau_i|^2)^{\frac{N-2}2}} & \left( \frac{r_i \eps}{|x-a|}\right)^{N-2}  U_{\d_i,\x_i}^p \\
& = \frac{1}{(1+|\tau_i|^2)^{\frac{N-2}2}}  \left( \frac{r_i \eps}{\d_i}\right)^{N-2} \int_{\frac{\Omega_\eps-\xi_i}{\d_i}}  \frac{dy}{|y+\tau_i|^{N-2} (1+|y|^2)^{\frac{N+2}2}} \\
& =  \left( \frac{\eps}{\d_i}\right)^{N-2} \frac{r_i^{N-2} \Gamma(\tau_i)}{(1+|\tau_i|^2)^{\frac{N-2}2}} + o \left(\left( \frac{\eps}{\d_i}\right)^{N-2}\right),
\end{split}
\end{equation}
as $\eps \to 0$. Also, by Lemma \ref{lem a.1} the term with $R_{\eps,\d_i,\x_i}$ is of lower order with respect to those in \eqref{984} and \eqref{985}, and can be absorbed in the small $o$ therein.

Since $\d_i \simeq \eps^\frac12$, it remains to estimate only the last term in \eqref{982}, and this can be done using \cite[Lemma 3.2]{GeMuPi}: as $0 \le P_\eps U_{\d_i,\x_i} \le U_{\d_i,\x_i}$, we have
\begin{align*}
\Bigg| \int_{\Omega_\eps} (U_{\d_i,\x_i}+  \eta_i ( P_\eps U_{\d_i,\x_i}- U_{\d_i,\x_i}))^{p-1} &( P_\eps U_{\d_i,\x_i}- U_{\d_i,\x_i})^2\Bigg|  \\
& \le C   \int_{\Omega_\eps} U_{\d_i,\x_i}^{p-1} ( P_\eps U_{\d_i,\x_i}- U_{\d_i,\x_i})^2 = o(\eps^\frac{N-2}2)
\end{align*}
Collecting all the previous computations, we infer that the expansion in the thesis holds $\C^0$ uniformly as $\eps \to 0$. The estimates for the derivatives can be obtained in a similar way.
\end{proof}
%
%
%
% can be obtained as in the proof of Lemma 4.1 in {MuPi}, see also \cite{MuPiErr} for the necessary corrections: using the last formulae on \cite[pages 208--209]{MuPiErr}, together with \cite[formulae (4.6) and (4.12)]{MuPi}, \cite[Lemma 2.3]{GeMuPi}, and the ansatz \eqref{asympt expansion}, we deduce that
%\begin{align*}
%I_\eps(P_\eps U_{\d,\x}) & = \frac{1}{N} \int_{\Omega_\eps} U_{\d,\x}^{p+1} -\frac12 \int_{\Omega_\eps} (P_\eps U_{\d,\x}-U_{\d,\x} ) U_{\d,\x}^p \\
%& - \frac{p}2\int_{\Omega_\eps} (U_{\d,\x} + \eta (V_{\d,\x}-U_{\d,\x}))^{p-1} (P_\eps U_{\d,\x}-U_{\d,\x} )^2 \\
%& = \frac{\alpha_N^{p+1} }{N} \left( \int_{\R^N} \frac{dy}{(1+|y^2)^N} +O(\d^N) \right) \\
%&+\frac{\alpha_N^{p+1}}2 \d^{N-2} H(0,0) \left( \int_{\R^N}  \frac{dy}{(1+|y|^2)^\frac{N+2}2} +o(1) \right) \\
%&  +\frac{\alpha_N^{p+1}}2 \d^{N-2} \left( \int_{\R^N}  \frac{dy}{(1+|\tau|^2)^{\frac{N-2}{2}} |y-\tau|^{N-2}(1+|y|^2)^\frac{N+2}2} +o(1) \right)\\
%& + \begin{cases} O(\d^4 |\log \d|) & \text{if $N=4$} \\ O(\d^2) & \text{if $N=3$} \end{cases}.
%\end{align*}
%Notice that the last term is, both for $N=3$ and $4$, of lower order with respect to $\d^{N-2}$. 

Coming back to \eqref{981}, we now show that the other terms are perturbation of $I_\eps(P_\eps U_{\d_i,\x_i})$.

\begin{lemma}\label{lem 5.4}
Let $R_{1,\eps}: X_\eta \to \R$ be defined by
\[
R_{1,\eps}(\mf{d},\bs{\tau}):= \tilde J_\eps(\mf{d},\bs{\tau}) - \sum_{i=1}^m \mu_i^{-\frac2{p-1}} I_\eps(P_\eps U_{\d_i,\x_i}).
\]
Then $R_{1,\eps} = o(\eps^\frac{N-2}2)$ $\C^1$-uniformly in $X_\eta$ as $\eps \to 0$.
\end{lemma}
\begin{proof}
Using the definitions of $P U_i$ and of $F$, we have
\begin{equation}\label{991}
\begin{split}
\Bigg[\frac12 \int_{\Omega_\eps} 2 \mu_i^{-\frac1{p-1}} & \nabla P U_i \cdot \nabla \phi_i - \int_{\Omega_\eps} \mu_i  \left( F(\mu_i^{-\frac1{p-1}}  P U_i +  \phi_i)   - F(\mu_i^{-\frac1{p-1}} P U_i)\right)\Bigg] \\
& = \int_{\Omega_\eps} \mu_i^{-\frac1{p-1}}(U_i^p - P U_i^p) \phi_i \\
& -  \int_{\Omega_\eps} \mu_i  \left( F(\mu_i^{-\frac1{p-1}}  P U_i +  \phi_i)   - F(\mu_i^{-\frac1{p-1}} P U_i) - F'(\mu_i^{-\frac1{p-1}} P U_i) \phi_i \right).
\end{split}
\end{equation}
The first term on the right hand side can be controlled using Corollary \ref{corol a.2}, Lemmas \ref{lem a.3}, \ref{lem a.4} and \ref{lem a.11}, the fact that $0 \le P U_i \le U_i$, and the H\"older and the Sobolev inequalities:
\begin{align*}
\Bigg| \int_{\Omega_\eps}  & (U_i^p - P U_i^p) \phi_i \Bigg| \le C \int_{\Omega_\eps} U_i^{p-1} |P U_i - U_i| |\phi_i| \\
& \le C \left( \int_{\Omega_\eps} U_i^\frac{8N}{(N-2)(N+2)} \left( \d_i^\frac{N-2}2 + \frac{\d_i^{\frac32N-3}}{|x-a_i|^{N-2}}\right) \right)^\frac{N+2}{2N} \|\phi_i\| \\
& \le  C \d_i^\frac{N-2}2 \left( \int_{\Omega_\eps} U_i^\frac{8N}{(N-2)(N+2)} \right)^\frac{N+2}{2N} \|\phi_i\| + C \d_i^{\frac32N-3} \left( \int_{\Omega_\eps} \frac{U_i^\frac{8N}{(N-2)(N+2)} }{|x-a_i|^\frac{2N(N-2)}{N+2} }\right)^\frac{N+2}{2N} \|\phi_i\|  \\
& \le C \d_i^\frac{N-2}2 \left(\d_i^\frac{N(N-2)}{N+2}\right)^\frac{N+2}{2N} \|\phi_i\|+  C \d_i^{\frac32N-3} \left(\d_i^\frac{N(2-N)}{N+2}\right)^\frac{N+2}{2N} \|\phi_i\| \le C \d_i^{N-2} \|\phi_i\| = o(\eps^\frac{N-2}2)
\end{align*}
uniformly in $X_\eta$ as $ \eps \to 0$, where the last equality follows by Proposition \ref{prop: pb K ort}. Regarding the second term on the right hand side in \eqref{991}, by the Lagrange theorem there exists a function $\eta_i$ with values in $[0,1]$ such that\begin{align*}
\int_{\Omega_\eps}   &  \left| F(\mu_i^{-\frac1{p-1}}  P U_i +  \phi_i)   - F(\mu_i^{-\frac1{p-1}} P U_i) - F'(\mu_i^{-\frac1{p-1}} P U_i) \phi_i \right| \\
& = \int_{\Omega_\eps} |F''(\mu_i^{-\frac1{p-1}} P U_i + \eta_i \phi_i) | \phi_i^2 \le C \int_{\Omega_\eps} \left( (P U_i)^{p-1} \phi_i^2 + |\phi_i|^{p+1} \right) \\
& \le C \int_{\Omega_\eps} U_i^{p-1} \phi_i^2 + C |\phi_i|_{2^*}^{2^*} \le C |U_i|_{2^*}^{p-1} |\phi_i|_{2^*}^2 + C \|\phi_i\|^{2*} \\
& \le C \|\phi_i\|^2 + \|\phi_i\|^{2^*} = o(\eps^\frac{N-2}2)
\end{align*}
uniformly in $X_\eta$ as $ \eps \to 0$. To sum up, the left hand side in \eqref{991} is $o(\eps^\frac{N-2}2)$ uniformly in $X_\eta$ as $\eps \to 0$.

To estimate the other terms on the right hand side in \eqref{981}, we use Lemmas \ref{lem a.4} and \ref{lem a.5}, the fact that $0 \le P_\eps U_{\d,\xi} \le U_{\d,\x}$, the ansatz \eqref{asympt expansion}, and the estimate in Proposition \ref{prop: pb K ort}: 
\[
\left| \int_{\Omega_\eps} P U_i^{\frac{p+1}2} P U_j^{\frac{p+1}2} \right|  \le \int_{\Omega_\eps} U_i^{\frac{N}{N-2}} U_j^{\frac{N}{N-2}}  \le C \d_i^\frac{N}2 \d_j^\frac{N}2 \left( |\log \d_i| + |\log \d_j| \right) = o(\eps^\frac{N-2}2), 
\]
\begin{align*}
\int_{\Omega_\eps} |\mu_i^{-\frac1{p-1}} & P U_i+ \eta_i \phi_i|^\frac{p+1}2 |\mu_j^{-\frac1{p-1}} P U_j+ \eta_j \phi_j|^\frac{p-1}2 |\phi_j| \\
& \le C \int_{\Omega_\eps} \left( U_i^\frac{p+1}2 U_j^\frac{p-1}2 |\phi_j| + |\phi_i|^\frac{p+1}2 U_j^\frac{p-1}2 |\phi_j| + U_i^\frac{p-1}2 |\phi_j|^\frac{p+1}2 + |\phi_i|^\frac{p+1}2|\phi_j|^\frac{p+1}2 \right) \\
&  \le C \left( \int_{\Omega_\eps} U_i^\frac{2N^2}{(N+2)(N-2)} U_j^\frac{4N}{(N+2)(N-2)} \right)^\frac{N+2}{2N} |\phi_j|_{2^*} + C |\phi_i|_{2^*}^\frac{p+1}2 |U_j|_{2^*}^\frac{p-1}2 |\phi_j|_{2^*} \\
& \qquad \qquad +  C|U_i|_{2^*}^\frac{p-1}2 |\phi_j|_{2^*}^\frac{p+1}2 + C + |\phi_i|_{2^*}^\frac{p+1}2|\phi_j|_{2^*}^\frac{p+1}2 \\
& \le C \d_i \d_j \|\phi_j\| + C  \|\phi_i\|^\frac{p+1}2\|\phi_j\| + C \|\phi_j\|^\frac{p+1}2 + C \|\phi_i\|^\frac{p+1}2 \|\phi_j\|^\frac{p+1}2 = o(\eps^\frac{N-2}2),
\end{align*}
and analogously
\begin{align*}
\int_{\Omega_\eps} |\mu_i^{-\frac1{p-1}} & P U_i+ \eta_i \phi_i|^\frac{p-1}2 |\mu_j^{-\frac1{p-1}} P U_j+ \eta_j \phi_j|^\frac{p+1}2 |\phi_i| = o(\eps^\frac{N-2}2)
\end{align*}
uniformly in $X_\eta$ as $ \eps \to 0$.

Collecting together the previous estimates, we deduce that $R_{1,\eps} = o(\eps^\frac{N-2}2)$ $\C^0$-uniformly as $\eps \to 0$. The estimates on the derivatives can be obtained in a similar way.
\end{proof}

Proposition \ref{prop: exp func} follows from Lemmas \ref{lem 5.3} and \ref{lem 5.4} straightforwardly, and allows us to complete the proof of Theorem \ref{thm: main 1}: we show that for every $\eps>0$ small enough the function $\tilde J_\eps$ has a critical point in $X_\eta$. A crucial lemma is the following.

\begin{lemma}\label{lem: ex e non-deg}
The function $\Psi$ defined in \eqref{def gamma e psi} has a non-degenrate critical point $(\tilde{\mf{d}}, 0) \in X_\eta$, provided $\eta>0$ was chosen small enough at the beginning.
\end{lemma}
\begin{proof}
We have
\begin{align*}
\pa_{d_i} \Psi(\mf{d},\bs{\tau}) &= (N-2) \mu_i^{-\frac2{p-1}}  b_2 H(a_i,a_i) d_i^{N-3} - (N-2) \frac{\alpha_N^{p+1} \mu_i^{-\frac2{p-1}} r_i^{N-2}}2 \frac{\Gamma(\tau_i)}{d_i^{N-1} (1+|\tau_i|^2)^\frac{N-2}2} \\
\pa_{\tau_{i,h}} \Psi(\mf{d},\bs{\tau}) &= \frac{\alpha_N^{p+1} \mu_i^{-\frac2{p-1}}r_i^{N-2}}{2 d_i^{N-2} } \left( \frac{\pa_{\tau_{i,h}}\Gamma(\tau_i)}{(1+|\tau_i|^2)^\frac{N-2}2} - \frac{(N-2) \Gamma(\tau_i) \tau_{i,h}}{ (1+|\tau_i|^2)^\frac{N}2} \right).
\end{align*}
As proved in \cite[Lemma 4.1]{GeMuPi}, the function $\Gamma$ has a non-degenerate maximum in $0$ (the Hessian matrix is diagonal with all negative entries), and hence we deduce that for every $\mf{d}$
\[
\nabla_{\bs{\tau}} \Psi(\mf{d},0) = 0.
\]
Thus, to find a critical point of $\Psi$, it is sufficient to find $\mf{d}$ with $\eta < d_i < \eta^{-1}$ for every $i$ such that 
\[
\nabla_{\mf{d}} \Psi(\mf{d},0) = 0.
\]
The existence of such critical point follows straightforwardly by the fact that 
\[
\Psi(\mf{d},0) = \sum_i g_i(d_i), \quad \text{with} \quad g_i(d_i) :=\tilde b_{i,1} d_i^{N-2} + \frac{\tilde b_{i,2}}{d_i^{N-2}}
\]
for suitable positive constants $\tilde b_{i,1}$, $\tilde b_{i,2}$, so that $\mf{d} \mapsto \Psi(\mf{d},0)$ admits a global minimum $\tilde{\mf{d}}$ in the set $\{\eta < d_i< \eta^{-1}\}$ (at least for $\eta>0$ small enough). 

Now we show that $(\tilde{\mf{d}},0)$ is non-degenerate. The Hessian matrix $D^2 \Psi(\tilde{\mf{d}},0)$ can be divided in blocks in the following way:
\[
D^2 \Psi(\tilde{\mf{d}},0) = \left(  \begin{array}{lr}   D^2_{\mf{d}} \Psi(\tilde{\mf{d}},0) & \left(\pa_{\tau_{i,h}, d_j} \Psi (\mf{d},0)\right) \\ \left(\pa_{\tau_{i,h}, d_j} \Psi (\mf{d},0)\right) & D^2_{\bs{\tau}} \Psi(\tilde{\mf{d}},0)
  \end{array} \right)
\]
Recalling that $0$ is a non-degenerate maximum for the function $\Gamma$, it is not difficult to check by direct computations that the matrix $D^2_{\bs{\tau}} \Psi(\tilde{\mf{d}},0)$ is diagonal and negative definite. Moreover, for every $i,j =1,\dots,m$ and $h=1,\dots,N$ we have $\pa_{\tau_{i,h}, d_j} \Psi (\tilde{\mf{d}},0) = 0$, and hence to prove the non-degeneracy of $(\tilde{\mf{d}},0)$ it remains only to analyze the block $D^2_{\mf{d}} \Psi(\tilde{\mf{d}},0)$. It is clear that this is another diagonal matrix, with
\[
\pa_{d_i , d_i} \Psi(\tilde{\mf{d}},0) =\frac{2}{\mu_i}  b_2 H(a_i,a_i) + 6 \frac{\alpha_4^4 r_i^2}{2 \mu_i} \frac{\Gamma(0)}{\tilde d_i^4}  \qquad \text{if $N=4$}, 
\]
and
\[
\pa_{d_i , d_i} \Psi(\tilde{\mf{d}},0) = 2 \frac{\alpha_3^6 r_i}{2 \mu_i^\frac12} \frac{\Gamma(0)}{\tilde d_i^3} \qquad \text{if $N=3$}.
\]
In any case, $D^2_{\mf{d}} \Psi(\tilde{\mf{d}},0)$ is diagonal and positive definite, and hence we deduce that $(\tilde{\mf{d}},0)$ is a non-degenerate saddle-point for $\Psi$.
\end{proof}

\begin{proof}[Conclusion of the proof of Theorem \ref{thm: main 1}]
Thanks to Lemma \ref{lem: the reduction works}, we prove that for $\eps>0$ small the function $\tilde J_\eps$ has a critical point in $X_\eta$. By Proposition \ref{prop: exp func}, this amounts to find a solution of the algebraic system
\[
\begin{cases}
\nabla_{\mf{d}} \Psi(\mf{d},\bs{\tau}) \eps^\frac{N-2}2 + \nabla_{\mf{d}} R_\eps (\mf{d},\bs{\tau})  = 0 \\
\nabla_{\bs{\tau}} \Psi(\mf{d},\bs{\tau}) \eps^\frac{N-2}2 + \nabla_{\bs{\tau}}  R_\eps(\mf{d},\bs{\tau}) = 0,
\end{cases}
\]
with $\nabla_{\mf{d}} R_\eps (\mf{d},\bs{\tau}), \nabla_{\bs{\tau}} R_\eps(\mf{d},\bs{\tau}) = o (\eps^\frac{N-2}2)$ as $\eps \to 0$, $\C^0$-uniformly in $X_\eta$. The previous system can be rewritten as 
\[
\begin{cases}
\nabla_{\mf{d}} \Psi(\mf{d},\bs{\tau})  + g_1(\mf{d},\mf{\tau},\eps)  = 0 \\
\nabla_{\bs{\tau}} \Psi(\mf{d},\bs{\tau}) +g_2(\mf{d},\mf{\tau},\eps) = 0,
\end{cases}
\]
with $g_1,g_2 = o(1)$ as $\eps \to 0$, $\C^0$-uniformly in $(\mf{d},\bs{\tau}) \in X_\eta$. Now, let us define the two maps
\[
\Lambda_0( \mf{d},\bs{\tau} ) := \left( \begin{array}{c}\nabla_{\mf{d}} \Psi(\mf{d},\bs{\tau}) \\ \nabla_{\bs{\tau}} \Psi(\mf{d},\bs{\tau}) \end{array} \right), \quad \quad \Lambda_{1,\eps}( \mf{d},\bs{\tau} ) := \left( \begin{array}{c}\nabla_{\mf{d}} \Psi(\mf{d},\bs{\tau}) + g_1(\mf{d},\bs{\tau},\eps) \\ \nabla_{\bs{\tau}} \Psi(\mf{d},\bs{\tau}) + g_2(\mf{d},\bs{\tau},\eps) \end{array} \right).
\]
The zeros of $\Lambda_{1,\eps}$ are critical point of the reduced functional $\tilde J_\eps$, and hence gives solutions to \eqref{system}. By Lemma \ref{lem: ex e non-deg}, we know that $\Lambda_0$ has a zero in $(\tilde{\mf{d}},0)$, and by non-degeneracy there exists a neighbourhood $U \subset X_\eta$ of $(\tilde{\mf{d}},0)$ such that $(\tilde{\mf{d}},0)$ is the unique $0$ of $\Lambda_0$ in $\overline{U}$, and $\textrm{deg}(\Lambda_0, U,0) \neq 0$. 

Let $H_\eps:[0,1] \times U \to \R^m \times \R^{Nm}$ be defined by 
\[
H_\eps(t,\mf{d},\bs{\tau}):= \Lambda_0(\mf{d},\bs{\tau}) + t \left( \begin{array}{c} g_1(\mf{d},\bs{\tau},\eps) \\ g_2(\mf{d},\bs{\tau},\eps) \end{array} \right).
\]
This is an homotopy between $\Lambda_0$ and $\Lambda_{1,\eps}$, and since $g_1,g_2 \to 0$ as $\eps \to 0$ uniformly in $X_\eta$, is such that $H_\eps(t,\mf{d},\bs{\tau}) \neq 0$ for every $t \in [0,1]$, $(\mf{d},\bs{\tau}) \in \pa U$, at least for $\eps>0$ small enough. As a consequence, by the homotopy-invariance property of the topological degree, we conclude that
\[
\textrm{deg}(\Lambda_0, U,0) = \textrm{deg}(\Lambda_{1,\eps}, H_\eps(1,U),0) \neq 0,
\]
and hence $\Lambda_{1,\eps}$ has a zero in $H(1,U)$. That is, $\tilde J_\eps$ has a critical point, as desired. This completes the proof of the existence of a solution $(u_{1,\eps},\dots,u_{m,\eps})$ to system \eqref{system 2}.

It remains to show that, if necessary replacing $\bar \beta$ with a smaller quantity, $u_{i,\eps}>0$ in $\Omega_\eps$ for every $i$, so that in particular $(u_{1,\eps},\dots,u_{m,\eps})$ solves \eqref{system}. 

We start from the case $N=4$, in which case we have $\beta_{ij} < \bar \beta$ for every $i \neq j$. This case can be treated exactly as in \cite[Conclusion of the proof of Theorem 1.1]{PisTav}. \footnote{Some careful is needed, since to apply the argument in \cite{PisTav} it is required that in the bound $\|\bs{\phi}\| \le C \eps^\frac{N-2}2$, given in Proposition \ref{prop: pb K ort}, the constant $C>0$ is independent on the particular choice of $\beta_{ij}$ with $-\infty < \beta_{ij} < \bar \beta$. Going through the proof of Proposition \ref{prop: pb K ort}, it is not difficult to check that this is possible.}

Regarding the case $N=3$, the positivity of the solutions (without any assumption on $\beta_{ij}$) can be obtained replacing system \eqref{system 2} with
\begin{equation}\label{system 3}
\begin{cases}
-\Delta u_i = \mu_i f(u_i) + \sum_{j \neq i} \beta_{ij} |u_j|^{\frac{p+1}{2}} |u_i|^{\frac{p-3}{2}} u_i^+ & \text{in $\Omega_\eps$} \\
u_i = 0 & \text{on $\pa \Omega_\eps$},
\end{cases} \qquad i =1,\dots,m.
\end{equation}
The key fact is that for $N=3$ the interaction term on the right hand side is of type $F(x,y) = |x| x^+ |y|^3$, which is of class $\mathcal{C}^1$, so that the proof we used to deal with system \eqref{system 2} in dimensions $N=3,4$ can be used word by word to produce a solution to \eqref{system 3} with $u_{i,\eps} \not \equiv 0$ for every $i$. This immediately implies (by the classical maximum principle) that $u_{i,\eps} >0$ in $\Omega_\eps$. 
\end{proof}

\begin{remark}
We stress that the strategy to prove the positivity of $u_{i,\eps}$ in dimension $N=3$ does not work in dimension $N=4$, since for the problem in $\R^4$ the interaction term is of type $F(x,y) = x^+ y^2$, which is not $\mathcal{C}^1$. We needed the smoothness of $F$ to prove Lemma \ref{lem a.6}, which is the key ingredient to estimate the nonlinear part for the equations in $\mf{K}_{\mf{d},\bs{\tau},\eps}^\perp$.
\end{remark}

\section{Proof of Theorem \ref{thm: main 2}}\label{sec: thm 2}

We start with some preliminaries about the shape of the approximate solution.

\subsection{Some preliminaries}
Let us consider the $2$ components system
\begin{equation}\label{2 comp}
\begin{cases}
-\Delta v_1= \mu_1 v_1^3 + \beta_{12} v_1 v_2^2 & \text{in $\R^4$} \\
-\Delta v_2= \mu_2 v_2^3 + \beta_{12} v_1^2 v_1^2 & \text{in $\R^4$} \\
v_1,v_2 >0 & \text{in $\R^4$} \\
v_1,v_2 \in \mathcal{D}^{1,2}(\R^4),
\end{cases}
\end{equation}
and let us search for solutions of the form $(v_1,v_2) = (c_1 U, c_2 U)$ (where $U=U_{1,0}$ is a standard bubble in $\R^4$ centered in $0$ and with $\d=1$), with $c_1,c_2 >0$. This ansatz leads to the algebraic system
\begin{equation}\label{syst c_i}
\mu_1 c_1^2 + \beta_{12} c_2^2 = 1, \qquad \beta_{12} c_1^2 + \mu_2 c_2^2 = 1,
\end{equation}
which admits the solution
\begin{equation}\label{def c_i}
c_1^2:= \frac{\beta_{12}-\mu_2}{\beta^2 - \mu_1 \mu_2},  \quad c_2^2:= \frac{\beta_{12}-\mu_1}{\beta^2 - \mu_1 \mu_2}
\end{equation}
if either $-\sqrt{\mu_1 \mu_2} <\beta_{12} < \min\{\mu_1,\mu_2\}$, or $\beta_{12}> \max\{\mu_1,\mu_2\}$.

Let us consider now the linearization of \eqref{2 comp} in $(c_1 U, c_2 U)$, namely the linear system
\begin{equation}\label{linearization 2}
\left\{ \begin{aligned}  -\Delta v_1&=(3\mu_1c_1^2+\beta_{12} c_2^2)U^2 v_1+2\beta_{12} c_1 c_2 U^2 v_2= U^2( \alpha_{11} v_1+ \alpha_{12} v_2)\quad \hbox{in}\ \mathbb R^4\\
   -\Delta v_2&=2\beta_{12} c_1 c_2 U^2 v_1+(3\mu_2c_2^2+\beta_{12} c_1^2)U^2 v_2=U^2(\alpha_{21} v_1+ \alpha_{22} v_2)\quad \hbox{in}\ \mathbb R^4.\\
 \end{aligned}\right.
 \end{equation}
 where
\begin{equation}\label{def alpha}
 \alpha_{11} :=3\mu_1 c_1^2+\beta_{12} c_2^2,\quad \alpha_{12} = \alpha_{21} := 2\beta_{12} c_1 c_2, \quad  \alpha_{22}:=3\mu_2c_2^2+\beta_{12} c_1^2.
\end{equation}
We introduce the $2 \times 2$ matrix $M:= (\alpha_{ij})_{i,j=1,2}$, with eigenvalues
\begin{equation}\label{def lambda i}
\lambda_2={\alpha_{11}+\alpha_{22}+\sqrt{(\alpha_{11}-\alpha_{22})^2+4\alpha_{12}^2}\over2}, \quad  \lambda_2={\alpha_{11}+\alpha_{22}-\sqrt{(\alpha_{11}-\alpha_{22})^2+4\alpha_{12}^2}\over2}.
\end{equation}
Using \eqref{syst c_i}, it is not difficult to check by direct computations that $\lambda_1 = 3$.

We consider now the eigenvalue problem:
\[
-\Delta v = \nu U^2 v, \qquad v \in \mathcal{D}^{1,2}(\R^4),
\]
It is well known (see \cite[Lemma A.1]{BiaEgn}) that there exists a sequence of positive eigenvalues $\{\nu_k\}$ with 
\[
1=\nu_1< 3=\nu_2 < \nu_3 < \dots < \nu_k< \nu_{k+1} < \dots,
\] 
$\nu_k \to +\infty$.

\begin{lemma}\label{lem: 5 ott 1}
Let $(e_1,e_2)$ be a non-trivial eigenvector of the matrix $M$ associated with $\lambda_2$. If $\lambda_2 \neq \nu_k$ for every $k$, then the set of solutions to the linear system \eqref{linearization 2}  is
 $5-$dimensional, and is generated by
 $$(e_2 ,-e_1 )\psi_{1,0}^h \qquad  h=0,1,\dots, 4$$
 (where $\psi_{1,0}^h$ have been defined in \eqref{def psi}).
 \end{lemma}
 
This lemma is a particular case of the forthcoming Lemma \ref{prop: base thm 2 general}, to which we refer for the proof. Notice that in Lemma \ref{prop: base thm 2} we express the generators of the set of solutions to \eqref{linearization 2} as $\mathfrak e_1 \psi_{1,0}^h$, with $\mathfrak e_1$ eigenvector associated to the eigenvalue $\lambda_1=3$, while here we use an eigenvector associate to $\lambda_2$. This is possible since, being eigenvectors associated to different eigenvalues orthogonal, $(e_2,-e_1)$ is indeed an eigenvector for $\lambda_1=3$.
% \begin{remark}
% Since both $c_1$ and $c_2$ are strictly positive, we have $e_1,e_2 \neq 0$ for any eigenvector associated with $\lambda_1$.
% \end{remark}
% \begin{proof}
% Let $(v_1,v_2)$ be a solution to \eqref{linearization 2}. We multiply the first equation in \eqref{linearization 2} by $e_1$,  the second equation by 
% $e_2$, and we add the results:
%\begin{equation}\label{4}-\Delta(e_1 v_1+e_2 v_2)=\underbrace{(\alpha_{11} e_1+\alpha_{12} e_2)}_{=\lambda_1 e_1}U^2 v_1+
%\underbrace{(\alpha_{12} e_1+\alpha_{22} e_2)}_{=\lambda_1 e_2}U^2 v_2.\end{equation}
%Letting $v:= e_1 v_1 + e_2 v_2$, we have
%$$-\Delta v=\lambda_1 U^2 v.$$
%If $\lambda_1$ is different from the sequence of the eigenvalues $\nu_k$, then we infer that $v=0$, i.e. 
%$$v_2=-\frac{e_1}{e_2} v_1={\alpha_{22}-\lambda_1\over \alpha_{12}}v_1.$$
%Then we come back the the first equation in \eqref{linearization 2} and we obtain
%$$ -\Delta v_1=U^2\underbrace{(\alpha_{11}+\alpha_{22}-\lambda_1)}_{=\lambda_2 = 3} v_1=3 U^2 v_1;$$
%due to \cite[Lemma A.1]{BiaEgn}, this implies that $v_1$ is the linear combination of $\psi_{1,0}^h$, and since $v_2 = -e_1 v_1 /e_2$, the desired conclusion follows. \end{proof}

\begin{proposition}\label{prop: base thm 2}
Under the assumptions of theorem \ref{thm: main 2}, the set of solutions to the linear system \eqref{linearization 2}  is
 $5-$dimensional, and is generated by
 $$(e_2 ,-e_1 )\psi_{1,0}^h \qquad  h=0,1,\dots, 4$$
 (where $\psi_{1,0}^h$ have been defined in \eqref{def psi}).
\end{proposition}

\begin{proof}
Notice at first that $\lambda_2<\lambda_1 = 3$. Then, by Lemma \ref{lem: 5 ott 1}, to complete the proof is sufficient to show that $\lambda_2$ is different from the eigenvalues $\nu_k$. Since $\lambda_2<3$ and $1=\nu_1<\nu_2=3$, we have to check that under the assumptions of Theorem \ref{thm: main 2} it results $\lambda_2 \neq 1$. Using the definition of $\alpha_{ij}$ and the one of $c_1,c_2$, it is not difficult to check that 
\[
\lambda_{1,2} = \frac{6-2\beta (c_1^2+c_2^2) \pm 2\beta (c_1^2+c_2^2)}2.
\]
If $\lambda_2=1$, i.e.
\[
1 = 3-2\beta_{12} \frac{2\beta_{12}-\mu_1-\mu_2}{\beta^2 -\mu_1 \mu_2},
\]
then it is not difficult to infer that either $\beta_{12} = \mu_1$, or $\beta_{12}=\mu_2$. Since in Theorem \ref{thm: main 2} we suppose that either $\beta_{12} <\min\{\mu_1,\mu_2\}$, or $\beta_{12} >\max\{\mu_1,\mu_2\}$, we have $\lambda_2 \neq 1$, and the thesis follows. 
\end{proof}

\subsection{The reduction scheme}

Once again, we search for solutions of \eqref{system 2}, which can be rewritten as in \eqref{pb 1}:
\[
u_i = i^*\left( \mu_i f(u_i) + \sum_{j \neq i} \beta_{ij} |u_j|^{\frac{p+1}{2}} |u_i|^{\frac{p-3}{2}} u_i\right).
\]
Let $\eta \in (0,1)$ be small, and let 
\begin{equation}\label{def X 2}
X_{\eta}:=\left\{ (\mf{d},\bs{\tau}) = (d_1,d_3, \tau_1,\tau_3)\in \R^2 \times (\R^4)^2: \ \eta < d_1,d_3< \eta^{-1}, \ |\tau_1|, |\tau_3|<\eta^{-1} \right\}.
\end{equation}
Our ansatz is that
\[
u_1 = c_1 P_\eps U_{\d_1,\x_1} + \phi_1, \quad u_2 = c_2 P_\eps U_{\d_1,\x_1} + \phi_2, \quad u_3 = c_3 P_\eps U_{\d_3,\x_3} + \phi_3
\]
where for some $(\mf{d},\bs{\tau}) = (d_1,d_3, \tau_1,\tau_3) \in X_\eta$ we have
\[
\d_i := d_i \sqrt{\eps}, \quad \xi_i:= a_i + d_i \sqrt{\eps} \tau_i,  \qquad i=1,3.
\]
Even though we have to deal with only $2$ parameters $d_1$, $d_3$ and two vectors $\tau_1$, $\tau_3$, in order to simplify the notation it is convenient to introduce $d_2= d_1$ and $\tau_2=\tau_1$. Analogously, we often write $\d_2=\d_1$ and $\xi_2=\xi_1$.

Plugging the previous ansatz into \eqref{pb 1}, our problem is transformed in the research of $d_i$, $\tau_i$ and $\phi_i$ such that \eqref{pb ansatz} is satisfied (with $c_i$ instead of $\mu_i^{-\frac1{p-1}}$) for $i=1,\dots,m$, with each equality which takes place in $H_0^1(\Omega_\eps)$. To proceed, the idea is again to split the space into two orthogonal subspaces, one of them having finite dimension. But in doing this we take into account that $u_1$ and $u_2$ are concentrating around the same point. Then we define
\begin{align*}
K_1 = K_{d_1,\tau_1,\eps} &:= \textrm{span}\left\{(e_2,-e_1) P_\eps \psi_{\d_1,\xi_1}^h: \ h=0,\dots,N\right\} \subset H_0^1(\Omega_\eps,\R^2), \\
K_3 =  K_{d_3,\tau_3,\eps} &:= \textrm{span}\left\{P_\eps \psi_{\d_3,\xi_3}^h: \ h=0,\dots,N\right\} \subset H_0^1(\Omega_\eps,\R) \\
\mf{K}_{\mf{d},\bs{\tau},\eps} &:= K_1 \times K_3,
\end{align*}
where $(e_1,e_2)$ is an eigenvector with norm $1$ of the matrix $M$ associated with $\lambda_2$ (defined in \eqref{def lambda i}). Notice that $K_{\mf{d},\bs{\tau},\eps}^{\perp} = K_1^\perp \times K_3^\perp$. 

If, for $i=1,3$, the symbol $\Pi_i = \Pi_{\d_i,\xi_i,\eps}$ (resp. $\Pi_i^{\perp} = \Pi_{\d_i,\xi_i,\eps}^{\perp}$) denotes the orthogonal projection $H_0^1(\Omega_\eps) \to K_i$ (resp. $H_0^1(\Omega_\eps) \to K_i^{\perp}$), then \eqref{pb ansatz} can be further rewritten as a system of $4$ equations. We have \eqref{pb K} and \eqref{pb K ort} for $i=3$ (with $\mu_j^{-\frac1{p-1}}$ replaced by $c_j$), together with
\begin{equation}\label{pb K 2}
\begin{split}
\Pi_1 & (c_1  P_\eps U_{\d_i,\x_i}  + \phi_1, c_2 P_\eps U_{\d_i,\x_i} + \phi_2)\\
&= \Pi_1 \Bigg( i^*\Bigg[ \mu_1 f(c_1 P_\eps U_{\d_1,\x_1} + \phi_1) + \sum_{j \neq 1} \beta_{1j} |c_j P_\eps U_{\d_j,\x_j} + \phi_j|^2 (c_1 P_\eps U_{\d_1,\x_1} + \phi_1)\Bigg], \\
& \hphantom{\ =\Pi_1 \Bigg( }  i^*\Bigg[ \mu_2 f(c_2 P_\eps U_{\d_1,\x_1} + \phi_2)  + \sum_{j \neq 2} \beta_{2j} |c_j P_\eps U_{\d_j,\x_j} + \phi_j|^2 (c_2 P_\eps U_{\d_1,\x_1} + \phi_2)\Bigg]\Bigg),
\end{split}
\end{equation}
and
\begin{equation}\label{pb K ort 2}
\begin{split}
\Pi_1^{\perp} & (c_1  P_\eps U_{\d_i,\x_i}  + \phi_1, c_2 P_\eps U_{\d_i,\x_i} + \phi_2)\\
&= \Pi_1^{\perp} \Bigg( i^*\Bigg[ \mu_1 f(c_1 P_\eps U_{\d_1,\x_1} + \phi_1) + \sum_{j \neq 1} \beta_{1j} |c_j P_\eps U_{\d_j,\x_j} + \phi_j|^2 (c_1 P_\eps U_{\d_1,\x_1} + \phi_1)\Bigg], \\
& \hphantom{\ =\Pi_1^\perp \Bigg( }  i^*\Bigg[ \mu_2 f(c_2 P_\eps U_{\d_1,\x_1} + \phi_2)  + \sum_{j \neq 2} \beta_{2j} |c_j P_\eps U_{\d_j,\x_j} + \phi_j|^2 (c_2 P_\eps U_{\d_1,\x_1} + \phi_2)\Bigg]\Bigg),
\end{split}
\end{equation}

\subsection{The equations in $\mf{K}_{\mf{d},\bs{\tau},\eps}^\perp$}
The equations in $\mf{K}_{\mf{d},\bs{\tau},\eps}^\perp$ can be still expressed in the form
\begin{equation}\label{L=N+R 2}
L_{\mf{d},\bs{\tau},\eps}^i(\bs{\phi}) = N_{\mf{d},\bs{\tau},\eps}^i(\bs{\phi}) + R_{\mf{d},\bs{\tau},\eps}^i \qquad i=1,3,
\end{equation}
where for $i=3$ the linear and nonlinear part and the remainder term are defined as in \eqref{def L}, \eqref{def N}, \eqref{def R}, while for $i=1$ we have
\begin{equation}\label{def L 2}
\begin{split}
& L_{\mf{d},\bs{\tau},\eps}^1(\bs{\phi})  =\\
& \Pi_1^{\perp} \Bigg( \phi_1 -i^*\Bigg[ \mu_1 f'(c_1 P_\eps U_{\d_1,\xi_1}) \phi_1  + \sum_{j \neq 1} \beta_{1j} (c_j P_\eps U_{\d_j,\xi_j})^2  \phi_1 
 %\\
%& \hphantom{= \Pi_1^{\perp} \Bigg( \phi_1 - i^*\Bigg[ \ } 
%\\
%&  \hphantom{= \Pi_1^{\perp} \Bigg( \ }  
+ 2 \sum_{j \neq 1} \beta_{1j} (c_j P_\eps U_{\d_j,\xi_j})
(c_1 P_\eps U_{\d_1,\xi_1}) \phi_j \Bigg], \\
&  \hphantom{\Pi_1^{\perp} \Bigg(  } 
\phi_2 -i^*\Bigg[ \mu_2 f'(c_2 P_\eps U_{\d_1,\xi_1}) \phi_2  + \sum_{j \neq 2} \beta_{2j} (c_j P_\eps U_{\d_j,\xi_j})^2  \phi_2 
% \\
%& \hphantom{= \Pi_1^{\perp} \Bigg( \phi_1 - i^*\Bigg[ \ } 
+ 2 \sum_{j \neq 2} \beta_{2j} (c_j P_\eps U_{\d_j,\xi_j})
(c_2 P_\eps U_{\d_1,\xi_1}) \phi_j \Bigg] \Bigg),
\end{split}
\end{equation}

\begin{equation}\label{def N 2}
\begin{split}
N_{\mf{d},\bs{\tau},\eps}^1 &(\bs{\phi}) = \\
& \Pi_1^{\perp}\Bigg(  i^* \Bigg[ \mu_1 f(c_1 P_\eps U_{\d_1,\xi_1}+\phi_1)-\mu_1 f(c_1P_\eps U_{\d_1,\xi_1}) - \mu_1 f'(c_1P_\eps U_{\d_1,\xi_1}) \phi_1  \\
& \hphantom{ \Pi_1^{\perp}\Bigg( i^* \Bigg[ \ } + \sum_{j \neq 1} \beta_{1j} |c_j P_\eps U_{\d_j,\x_j} + \phi_j|^2 (c_1 P_\eps U_{\d_1,\x_1} + \phi_1)  - \beta_{1j} (c_j P_\eps U_{\d_j,\x_j})^2 (c_1 P_\eps U_{\d_1,\x_1}) \\
& \hphantom{ \Pi_1^{\perp}\Bigg( i^* \Bigg[ \ }  - \sum_{j \neq 1} \beta_{1j} (c_j P_\eps U_{\d_j,\xi_j})^2 \phi_1  - 2  \beta_{1j} (c_j P_\eps U_{\d_j,\xi_j})(c_1 P_\eps U_{\d_1,\xi_1}) \phi_j \Bigg] , \\
& \hphantom{\Pi_1^{\perp}\Bigg(  \ } i^* \Bigg[ \mu_2 f(c_2 P_\eps U_{\d_1,\xi_1}+\phi_2)-\mu_2 f(c_2P_\eps U_{\d_1,\xi_1}) - \mu_2 f'(c_2P_\eps U_{\d_1,\xi_1}) \phi_2  \\
& \hphantom{ \Pi_1^{\perp}\Bigg( i^* \Bigg[ \ } + \sum_{j \neq 2} \beta_{2j} |c_j P_\eps U_{\d_j,\x_j} + \phi_j|^2 (c_2 P_\eps U_{\d_1,\x_1} + \phi_2)  - \beta_{2j} (c_j P_\eps U_{\d_j,\x_j})^2 (c_2 P_\eps U_{\d_1,\x_1}) \\
& \hphantom{ \Pi_1^{\perp}\Bigg( i^* \Bigg[ \ }  - \sum_{j \neq 2} \beta_{2j} (c_j P_\eps U_{\d_j,\xi_j})^2 \phi_2  - 2 \beta_{2j} (c_j P_\eps U_{\d_j,\xi_j})(c_2 P_\eps U_{\d_1,\xi_1}) \phi_j \Bigg]\Bigg),
\end{split}
\end{equation}
and 
\begin{equation}\label{def R 2}
\begin{split}
R_{\mf{d},\bs{\tau},\eps}^1  
%&= \Pi_1^{\perp} \Bigg(  -c_1 P_\eps U_{\d_1,\x_1}   + i^*\Bigg[ \mu_1 f(c_1 P_\eps U_{\d_1,\xi_1}) + \sum_{j \neq 1} \beta_{1j} (c_j P_\eps U_{\d_j,\x_j})^2(c_1 P_\eps U_{\d_1,\x_1}) \Bigg], \\
%& \hphantom{= \Pi_1^{\perp} \Bigg( \ }-c_2 P_\eps U_{\d_1,\x_1}   + i^*\Bigg[ \mu_2 f(c_2 P_\eps U_{\d_1,\xi_1}) + \sum_{j \neq 2} \beta_{2j} (c_j P_\eps U_{\d_j,\x_j})^2(c_2 P_\eps U_{\d_1,\x_1}) \Bigg] \Bigg) \\
& = \Pi_1^{\perp} \Bigg( i^*\Bigg[ c_1 (P_\eps U_{\d_1,\x_1}-U_{\d_1,\x_1}) + \beta_{13} (c_3 P_\eps U_{\d_3,\x_3})^2 (c_1 P_\eps U_{\d_1,\x_1}) \Bigg], \\
& \hphantom{= \Pi_1^{\perp} \Bigg( \ } i^*\Bigg[ c_2 (P_\eps U_{\d_1,\x_1}-U_{\d_1,\x_1}) + \beta_{23} (c_3 P_\eps U_{\d_3,\x_3})^2 (c_2 P_\eps U_{\d_1,\x_1}) \Bigg]\Bigg),
\end{split}
\end{equation} 
where we used the definition of $i^*$ and the equations \eqref{syst c_i} defining $c_1,c_2$.

We define 
\[
\mf{L}_{\mf{d},\bs{\tau},\eps} := (L_{\mf{d},\bs{\tau},\eps}^1,L_{\mf{d},\bs{\tau},\eps}^3): \mf{K}_{\mf{d},\bs{\tau},\eps}^\perp \to \mf{K}_{\mf{d},\bs{\tau},\eps}^\perp,
\]
and $\mf{R}_{\mf{d},\bs{\tau},\eps}$ and $\mf{N}_{\mf{d},\bs{\tau},\eps} $ in an analogue way.

%Given $\eta>0$ small, let us set
%\begin{equation}\label{def X}
%X_{\eta}:=\left\{ (\mf{d},\bs{\tau}) \in \R^m \times (\R^N)^m: \ \eta < d_i< \eta^{-1}, \ |\tau_i|<\eta^{-1} \right\}.
%\end{equation}

The main result of this subsection is the counterpart of Proposition \ref{prop: pb K ort} in the present setting.

\begin{proposition}\label{prop: pb K ort 2}
For every $\eta>0$ small enough there exists $\bar \beta, \bar \eps>0$ small, and $C>0$, such that if $\eps \in (0,\bar \eps)$, and $-\infty <\beta_{13},\beta_{23}<\bar \beta$, then for any $(\mf{d},\bs{\tau}) \in X_\eta$ (see \eqref{def X 2}) there exists a unique function $\bs{\phi}^{\mf{d},\bs{\tau},\eps} \in K_{\mf{d},\bs{\tau},\eps}^\perp$ solving the equation
\[
\mf{L}_{\mf{d},\bs{\tau},\eps}(\bs{\phi}) = \mf{R}_{\mf{d},\bs{\tau},\eps} + \mf{N}_{\mf{d},\bs{\tau},\eps}(\bs{\phi})
\]
and satisfying
\[
\|\bs{\phi}^{\mf{d},\bs{\tau},\eps}\|_{H_0^1(\Omega_\eps)} \le C \eps.
\] 
Furthermore, the map $(\eps,\mf{d},\bs{\tau}) \mapsto \bs{\phi}^{\mf{d},\bs{\tau},\eps}$ is of class $\mathcal{C}^1$, and
\[
\|\nabla_{(\mf{d},\bs{\tau})} \bs{\phi}^{\mf{d},\bs{\tau},\eps}\|_{H_0^1(\Omega_\eps)} \le C \eps. 
\] 
\end{proposition}

For the proof, we start studying the linear part.

\begin{lemma}\label{lem: linear part 2}
For every $\eta>0$ small enough there exists $\bar \beta, \eps_0>0$ small, and $C>0$, such that if $\eps \in (0,\eps_0)$, and $-\infty <\beta_{13},\beta_{23}<\bar \beta$, then
\[
\|\mf{L}_{\mf{d},\bs{\tau},\eps}(\bs{\phi})\|_{H_0^1(\Omega_\eps)} \ge C \|\bs{\phi}\|_{H_0^1(\Omega_\eps)}  \qquad \forall \bs{\phi} \in H_0^1(\Omega_\eps,\R^3)
\]
for every $(\mf{d}, \bs{\tau}) \in X_\eta$. Moreover, $\mf{L}_{\mf{d},\bs{\tau},\eps}$ is invertible in $\mf{K}_{\mf{d},\bs{\tau},\eps}^\perp$, with continuous inverse.
\end{lemma}

\begin{proof}
The proof follows exactly the same sketch of that of Lemma \ref{lem: linear part}. We suppose by contradiction that there exist sequences 
\[
\{\eps_n\} \subset \R^+, \ \eps_n \to 0, \  \{(\mf{d}_n,\boldsymbol{\tau}_n)\} \subset X_\eta,\  \{\bs{\phi}_n= ((\phi_{1,n},\phi_{2,n}),\phi_{3,n})\} \subset K_{1,n}^\perp \times K_{3,n}^{\perp}
\]
such that
\[
\| \bs{\phi}_n\|_{H_0^1(\Omega_{\eps_n})}=1 \quad \text{and} \quad \|\mf{L}_{n}(\bs{\phi}_n)\|_{H_0^1(\Omega_{\eps_n})} \to 0
\]
as $n \to \infty$, where we adopted the same short notation as in Lemma \ref{lem: linear part}.

Let $\mf{h}_n:= \mf{L}_{n}(\bs{\phi}_n)$. Then we have three equations for $\phi_{1,n}$, $\phi_{2,n}$, $\phi_{3,n}$, completely analogue to \eqref{eq phi} (with $c_i$ instead of $\mu_i^{-\frac1{p-1}}$), but with $(w_{1,n},w_{2,n}) \in K_{1,n}$, $w_{3,n} \in K_{3,n}$.

\smallskip

\noindent \textbf{Step 1)} $\|w_{i,n}\|_{H_0^1(\Omega_n)} \to 0$ as $n \to \infty$. \\
The argument used in Lemma \ref{lem: linear part} immediately implies that $\|w_{3,n}\|_{H_0^1(\Omega_n)} \to 0$. Regarding $w_{1,n}$ and $w_{2,n}$, we multiply the first equation with $\d_{1,n}^2 w_{1,n}$, the second equation with $\d_{1,n}^2 w_{2,n}$, and we sum the results. Since $(\phi_{1,n},\phi_{2,n}), (h_{1,n},h_{2,n}) \in K_{1,n}^\perp$, and recalling the definition of $\alpha_{ij}$ given in \eqref{def alpha}, we obtain
\begin{equation}\label{4.8 2}
\begin{split}
\d_{1,n}^2 &( \|w_{1,n}\|^2 + \|w_{2,n}\|^2 ) = \int_{\Omega_n} c_3^2 (P_n U_{3,n})^2 (\beta_{13} \phi_{1,n} w_{1,n}+ \beta_{23} \phi_{2,n} w_{2,n} ) \\
&  + \int_{\Omega_n} c_3 (P_n U_{1,n}) (P_n U_{3,n}) (c_1 \beta_{13} \phi_{3,n} w_{1,n}+ c_2 \beta_{23} \phi_{3,n} w_{2,n} ) \\
& + \d_{1,n}^2 \int_{\Omega_n} (P_n U_{1,n})^2 (\alpha_{11} \phi_{1,n} w_{1,n} + \alpha_{12} (\phi_{2,n} w_{1,n} + \phi_{1,n} w_{2,n}) + \alpha_{22} \phi_{2,n} w_{2,n}). 
\end{split}
\end{equation}
In this equation, the first and the second integral on the right hand side can be treated as terms ($III$) and ($IV$) in \eqref{eq phi test}, and together give $o(\d_{1,n}^2)$ as $n \to \infty$. 

The integral on the left hand side can be treated developed as term ($I$) in \eqref{eq phi test}, in the following way: since $(w_{1,n},w_{2,n}) \in K_{1,n}$, there exists constants $c_{1,n}^h$, $h=0,\dots,4$ such that
\begin{equation}\label{6 ott 1}
(w_{1,n},w_{2,n}) = \sum_{h=0}^4 c_{1,n}^h (e_2, -e_1) P_n \psi_{1,n}^h;
\end{equation}
then, for $\sigma_{lk}$ defined in \eqref{sigma lk}, we have 
%\begin{equation}\label{6101}
\[
\begin{split}
\d_{1,n}^2 ( \|w_{1,n}\|^2 + \|w_{2,n}\|^2 ) &= \sum_{l,k=0}^4 c_{1,n}^l c_{1,n}^k \underbrace{(e_1^2 + e_2^2)}_{=1} (\sigma_{lk}+ o(1)) \\
& = \sum_{h=0}^4 (c_{1,n}^h)^2 \sigma_{hh} + o(1) \sum_{h,k=0}^4 c_{1,n}^h c_{1,n}^k .
\end{split}
\]
%\end{equation}
%The second an the third integrals in \eqref{4.8 2} can be treated as terms ($III$) and ($IV$) in \eqref{eq phi test}, respectively, and together give a $o(\d_{1,n}^2)$. It remains to analyze the last integral, which plays the role of term ($II$) in \eqref{eq phi test}. Using \eqref{6 ott 1}, we rewrite it as
%\begin{equation}\label{6 ott 2}
%\d_{1,n}^2 \sum_{h=0}^4 c_{1,n}^h   \int_{\Omega_n} (P_n U_{1,n})^2 (P_n \psi_{1,n}^h)  (\alpha_{11} e_2 \phi_{1,n}  + \alpha_{12} ( e_2\phi_{2,n}  - e_1 \phi_{1,n} ) - \alpha_{22} e_1 \phi_{2,n})
%\end{equation}

It remains then to analyze the last integral on the right hand side in \eqref{4.8 2}, and in what follows we prove that it gives 
\[
o(\d_{1,n}^2) \left( \|w_{1,n}\| + \|w_{2,n}\| \right) + O(\d_{i,n}^2) \sum_{l=0}^N c_{i,n}^l.
\]
Now, since $(\phi_{1,n},\phi_{2,n}) \in K_{1,n}^\perp$, for every $h$
\begin{align*}
0 &= \int_{\Omega_n} e_2 \nabla (P_n \psi_{1,n}^h) \cdot \nabla \phi_{1,n} - e_1 \nabla (P_n \psi_{1,n}^h) \cdot \nabla \phi_{2,n} = 3 \int_{\Omega_n} U_{1,n}^2 \psi_{1,n}^h (e_2 \phi_{1,n} - e_1 \phi_{2,n}).
\end{align*}   
Therefore, for every $h$
\begin{align*}
\int_{\Omega_n} & U_{1,n}^2 \psi_{1,n}^h (\alpha_{11} e_2 \phi_{1,n}  + \alpha_{12} ( e_2\phi_{2,n}  - e_1 \phi_{1,n} ) - \alpha_{22} e_1 \phi_{2,n}) \\
& = \int_{\Omega_n} \alpha_{11} U_{1,n}^2 \psi_{1,n}^h (e_2 \phi_{1,n} \pm e_1 \phi_{2,n}) + \int_{\Omega_n} \alpha_{12} U_{1,n}^2 \psi_{1,n}^h (e_2 \phi_{2,n}-e_1 \phi_{1,n} ) \\
& \hphantom{= \int_{\Omega_n} \ }  - \int_{\Omega_n} \alpha_{22} U_{1,n}^2 \psi_{1,n}^h (e_1 \phi_{2,n} \pm e_2 \phi_{1,n}) \\
& = \int_{\Omega_n} (\alpha_{11} e_1 +\alpha_{12} e_2) U_{1,n}^2 \psi_{1,n}^h \phi_{2,n} - \int_{\Omega_n} (\alpha_{12} e_1 +\alpha_{22} e_2) U_{1,n}^2 \psi_{1,n}^h \phi_{1,n} \\
& = \lambda_2 \int_{\Omega_n} U_{1,n}^2 \psi_{1,n}^h (e_1 \phi_{2,n}-e_2 \phi_{1,n}) = 0,
\end{align*}
where we used the definition of $\alpha_{ij}$, see \eqref{def alpha}. In turn, using \eqref{6 ott 1}, we have
\begin{align*}
\d_{1,n}^2 &\int_{\Omega_n} (P_n U_{1,n})^2 (\alpha_{11} \phi_{1,n} w_{1,n} + \alpha_{12} (\phi_{2,n} w_{1,n} + \phi_{1,n} w_{2,n}) + \alpha_{22} \phi_{2,n} w_{2,n}) \\
& = \d_{1,n}^2 \int_{\Omega_n} ((P_n U_{1,n})^2 - U_{1,n})^2 (\alpha_{11} \phi_{1,n} w_{1,n} + \alpha_{12} (\phi_{2,n} w_{1,n} + \phi_{1,n} w_{2,n}) + \alpha_{22} \phi_{2,n} w_{2,n}) \\
&+ \d_{1,n}^2\sum_{h=0}^4 c_{1,n}^h \int_{\Omega_n}  U_{1,n}^2 (P_n \psi_{1,n}^h -\psi_{1,n}^h)  (\alpha_{11} e_2 \phi_{1,n}  + \alpha_{12} ( e_2\phi_{2,n}  - e_1 \phi_{1,n} ) - \alpha_{22} e_1 \phi_{2,n}),
\end{align*} 
and hence we can proceed exactly as in \eqref{265}, deducing that 
\begin{multline*}
\d_{1,n}^2 \int_{\Omega_n} (P_n U_{1,n})^2 (\alpha_{11} \phi_{1,n} w_{1,n} + \alpha_{12} (\phi_{2,n} w_{1,n} + \phi_{1,n} w_{2,n}) + \alpha_{22} \phi_{2,n} w_{2,n}) \\
 = o(\d_{1,n}^2) \left( \|w_{1,n}\| + \|w_{2,n}\| \right) + O(\d_{i,n}^2) \sum_{l=0}^N c_{i,n}^l.
\end{multline*}

To sum up, we proved the analogue of the estimates \eqref{estimates step 1} for the couple $(w_{1,n}, w_{2,n})$, and hence we can conclude as in Lemma \ref{lem: linear part} that $\|w_{1,n}\|,\|w_{2,n}\| \to 0$ as $n \to \infty$.

\smallskip

\noindent \textbf{Step 2)} For a fixed $\kappa=1,\dots,m$, we introduce
\[
\tilde \phi_{i,n}^\kappa(y):= \begin{cases} \d_{\kappa,n}^{\frac{N-2}{2}} \phi_{i,n}(\x_{\kappa,n}+\d_{\kappa,n} y) & y \in \frac{\Omega_n-\x_{\kappa,n}}{\d_{\kappa,n}}=:\tilde \Omega_{\kappa,n} \\
0 & y \in \R^N \setminus \tilde\Omega_{\kappa,n},
\end{cases}   \quad i=1,\dots,m.
\]
In a completely analogue way, we define $\tilde h_{i,n}^\kappa$ and $\tilde w_{i,n}^\kappa$. Proceeding as in step 2 of Lemma \ref{lem: linear part} (with minor changes), it is not difficult to check that $\tilde \phi_{1,n}^1, \tilde \phi_{2,n}^1, \tilde \phi_{3,n}^3 \wc 0$ in $\mathcal{D}^{1,2}(\R^N)$ as $n \to \infty$. The only difference with respect to Lemma \ref{lem: linear part} is that this time the weak limit $(\tilde \phi_1^1, \tilde \phi_2^1)$ solves the system \eqref{linearization 2} instead of a single equation, and hence the fact that $(\tilde \phi_1^1, \tilde \phi_2^1)= (0,0)$ comes from Proposition \ref{prop: base thm 2} and the condition $(\phi_{1,n},\phi_{2,n}) \in K_{1,n}^\perp$. 

\smallskip

\noindent \textbf{Step 3)} We prove that $\|\phi_{i,n}\|_{H_0^1(\Omega_n)} \to 0$ as $n \to \infty$ for every $i$. Regarding $\phi_{3,n}$, we can proceed as in Lemma \ref{lem: linear part}. We focus then on the other components. We test \eqref{eq phi} with $\phi_{1,n}$: recalling that $\{\phi_{1,n}\}$ is bounded in $H_0^1(\Omega_n)$ and that $w_{1,n}, h_{1,n} \to 0$ strongly, we deduce that
\begin{equation}\label{eq phi phi 2}
\begin{split}
\|\phi_{1,n}\|_{H_0^1(\Omega_n)}^2 & = o(1) + \int_{\Omega_n} (P_n U_{i,n})^2 (\alpha_{11} \phi_{1,n}^2 + \alpha_{12} \phi_{1,n} \phi_{2,n}) \\
& +  \beta_{13} \int_{\Omega_{n}}   (c_3 P_n U_{3,n})^2  \phi_{1,n}^2 + 2\beta_{23} \int_{\Omega_{n}} (c_1 P_n U_{i,n}) (c_3 P_n U_{3,n}) \phi_{1,n} \phi_{3,n}
\end{split}
\end{equation}
The second and the third integral on the right hand side can be treated as in Lemma \ref{lem: linear part}. Notice that here we have only to ask that $\beta_{13}$ is small enough (we don't need any assumption on $\beta_{12}$). 

Now, to estimate the first integral in \eqref{eq phi phi 2}, we observe that as in \eqref{3061}
\[
\int_{\Omega_n} (P_n U_{1,n})^2  \phi_{1,n}^2 \to 0,
\]
and moreover
\begin{align*}
\Bigg| \int_{\Omega_n} (P_n U_{1,n})^2 \phi_{1,n} \phi_{2,n} \Bigg| &\le C \int_{\Omega_n} U_{1,n}^2 |\phi_{1,n}| \, |\phi_{2,n}| \\
& = C \int_{\tilde \Omega_{1,n}} U_{1,0}^2 \tilde \phi_{1,n}^1 \tilde \phi_{2,n}^1 \to 0
\end{align*}
as $n \to \infty$, since $\tilde \phi_{1,n}^1, \tilde \phi_{2,n}^1 \wc 0$ in $L^{\frac{2N}{N-2}}(\R^N)$ by step 2, and $U_{1,0}^2 \in L^{\frac{N}{2}}(\R^N)$. Thus, we infer that $\|\phi_{1,n}\| \to 0$; in the same way, $\|\phi_{2,n}\| \to 0$, and we reached the desired contradiction. 

We stress that, while we have to suppose $-\infty < \beta_{13}, \beta_{23}< \bar \beta$, no assumption is needed on $\beta_{12}$.

\smallskip

\noindent{Step 4)} Invertibility of $\mf{L}_{\mf{d},\bs{\tau},\eps}$. This can be proved exactly as in Lemma \ref{lem: linear part}.
\end{proof}

The rest of the proof of Proposition \ref{prop: pb K ort 2} is now a straightforward modification of that of Proposition \ref{prop: pb K ort}, and hence is omitted.

\subsection{The reduced problem}

In this section we solve equation \eqref{pb K 2} with $\bs{\phi}= \bs{\phi}^{\mf{d},\bs{\tau},\eps}$. In what follows we use the notation 
\[
V_i^{\mf{d},\bs{\tau},\eps} := c_i P_\eps U_{\d_1,\x_1}^{\mf{d},\bs{\tau},\eps}+ \phi_i^{\mf{d},\bs{\tau},\eps} \quad i=1,2, \quad V_3^{\mf{d},\bs{\tau},\eps} := c_3 P_\eps U_{\d_3,\x_3}^{\mf{d},\bs{\tau},\eps}+ \phi_3^{\mf{d},\bs{\tau},\eps}.
\]

Let $J_\eps: H_0^1(\Omega_\eps,\R^3) \to \R$ be the action functional, defined as in \eqref{def action}. Critical points of $J_\eps$ are solution to \eqref{system 2}, and hence solutions to \eqref{system}.

Let us introduce now the reduced functional $\tilde J_\eps: X_\eta \to \R$, 
\[
\tilde J_\eps(\mf{d},\bs{\tau}):= J_\eps\left( c_1 P_\eps U_{\d_1,\x_1} + \phi_1^{\mf{d},\bs{\tau},\eps}, c_2 P_\eps U_{\d_1,\x_1} + \phi_2^{\mf{d},\bs{\tau},\eps} , c_3 P_\eps U_{\d_3,\x_3} + \phi_3^{\mf{d},\bs{\tau},\eps}\right).
\] 
The counterpart of Lemma \ref{lem: the reduction works} in the present context is given in the following statement, whose proof is omitted.
\begin{lemma}\label{lem: the reduction works 2}
There exists $\bar \eps>0$ sufficiently small such that if $(\mf{d},\bs{\tau}) = (d_1,d_3, \tau_1,\tau_3)$ is a critical point of $\tilde J_\eps$, and $\eps \in (0, \bar \eps)$, then 
\[
\left( c_1 P_\eps U_{\d_1,\x_1} + \phi_1^{\mf{d},\bs{\tau},\eps}, c_2 P_\eps U_{\d_1,\x_1} + \phi_2^{\mf{d},\bs{\tau},\eps} , c_3 P_\eps U_{\d_3,\x_3} + \phi_3^{\mf{d},\bs{\tau},\eps}\right)
\]
is a solution to \eqref{pb K}, and hence a solution to \eqref{system 2}.
\end{lemma}

Due to Lemma \ref{lem: the reduction works 2}, we have to find a critical point of the reduced functional $\tilde J_\eps$ in $X_\eta$, and to this aim we derive the asymptotic expansion of $J_\eps(V_1^{\mf{d},\bs{\tau},\eps},V_2^{\mf{d},\bs{\tau},\eps},V_3^{\mf{d},\bs{\tau},\eps})$ as $\eps \to 0^+$. We recall the definition of $b_1,b_2$ and $\Gamma$, given in \eqref{def b_i} and \eqref{def gamma e psi}, but we modify the definition of $\Psi$ in the following way:
\begin{equation}\label{def gamma e psi 2}
\begin{split}
\Psi(\mf{d},\bs{\tau}) & := (c_1^2 + c_2^2) \left[ b_2 H(a_1,a_1) d_1^{N-2}+ \frac{\alpha_N^{p+1} r_1^{N-2}}2    \frac{ \Gamma(\tau_1)}{d_1^{N-2}(1+|\tau_1|^2)^{\frac{N-2}2}}\right]\\
& + c_3^2 \left[ b_2 H(a_1,a_1) d_1^{N-2}+ \frac{\alpha_N^{p+1} r_1^{N-2}}2    \frac{ \Gamma(\tau_1)}{d_1^{N-2}(1+|\tau_1|^2)^{\frac{N-2}2}}\right].
\end{split}
\end{equation}

\begin{proposition}\label{prop: exp func 2}
We have 
\[
\tilde J_\eps(\mf{d},\bs{\tau}) = \left(\sum_{i=1}^3 c_i^2 \right) b_1 + \Psi(\mf{d},\bs{\tau})\eps^\frac{N-2}2 + R_\eps(\mf{d},\bs{\tau}),
\]
with $R_\eps = o(\eps^\frac{N-2}2)$ $\C^1$-uniformly in $X_\eta$ as $\eps \to 0$.
\end{proposition}
\begin{proof}
In this proof we omit the dependence of the functions on $\mf{d}$, $\bs{\tau}$, $\eps$ for simplicity, and we write $U_1 = U_2 = U_{\d_1,\x_1}$, $U_3 = U_{\d_3,\x_3}$.

We proceed as in Proposition \ref{prop: exp func}, noting that by \eqref{syst c_i}
\[
\begin{split}
\tilde J_\eps(\mf{d},\bs{\tau}) &= \sum_i I_\eps(V_i) - \frac2{p+1}\sum_{i<j}\int_{\Omega_\eps} \beta_{ij} |V_i|^\frac{p+1}2 |V_j|^\frac{p+1}2
\\
& = \sum_i    I_\eps(c_i P U_i) + \frac12\sum_i \int_{\Omega_\eps} 2 c_i \nabla (P U_i) \cdot \nabla \phi_i + |\nabla \phi_i|^2 \\
& - \sum_i \int_{\Omega_\eps} \mu_i  \left( F( c_i  P U_i +  \phi_i)   - F(c_i P U_i)\right)
\\
&  - \frac{2}{p+1} \sum_{i<j} \int_{\Omega_\eps} |c_iP U_i|^2 |c_j P U_j|^2
\\
& - \sum_{i<j} \int_{\Omega_\eps} \left( |c_i P U_i + \eta _i \phi_i|^2 |c_j P U_j + \eta _j \phi_j| \phi_j \right. \\
& \hphantom{- \sum_{i<j} \int_{\Omega_\eps} \Big( } \left.  + |c_i P U_{i} + \eta_i \phi_i| |c_j P U_j+ \eta_j \phi_j|^2 \phi_i\right) \\
& = (c_1^2 + c_2^2) I_\eps(P U_1) + c_3^2 I_\eps(P U_3) + R_{1,\eps}
\end{split}
\]
(where $\eta_i, \eta_j$ denote continuous functions with values in $[0,1]$). Then we can slightly modify Lemmas \ref{lem 5.3} and \ref{lem 5.4} to obtain the desired result.
\end{proof}

Once that Proposition \ref{prop: exp func 2} is established, the conclusion of the proof of Theorem \ref{thm: main 2} is analogue to the one of Theorem \ref{thm: main 1}.

\section{Proof of Proposition \ref{prop: on hp 2}}\label{sec: prop on hp 2}

Let us consider a system
\begin{equation}\label{sist1}
-\Delta U_i=\sum\limits_{j =1}^k \beta_{ij} U_i U_j^2 \qquad \hbox{in}\ \mathbb R^4,  \ i=1,\dots,k,
\end{equation}
and let us suppose that it has a solution $U_i=c_i U,$ with $c_i>0$ for $i =1,\dots,k$, where $U=U_{1,0}$ is a standard bubble in $\R^4$ centered in $0$ and with $\d=1$. That is,
\[
 \sum\limits_{j=1,\dots,k}\beta_{ij}  c_j^2=1, \quad  i=1,\dots,k.
\]
We can linearize system \eqref{sist1} in $\mathcal{D}^{1,2}(\R^N)$ around the solutions $(c_1 U,\dots,c_k U)$, obtaining
\begin{equation}\label{sist3}
-\Delta v_i  =\left[ \underbrace{ \left(3\beta_{ii}c_i^2+\sum\limits_{j\not=i}\beta_{ij}  c_j^2 \right) }_{= 1+ 2 \beta_{ii} c_i^2   }v_i+2 \sum\limits_{j\not=i}\beta_{ij}c_i c_j v_j\right]U^2  \qquad \text{in }\mathbb R^4,
\end{equation}
with $v_i \in \mathcal{D}^{1,2}(\R^4)$.

It is clear that Proposition \ref{prop: on hp 2} follows if we prove that the condition
\[
\text{the matrix $(\beta_{ij})_{i,j=1,\dots,k}$ is invertible and has only positive entries}
\]
implies that \eqref{sist3} has a $5$-dimensional set of solutions.

To this aim, we observe that system \eqref{sist3} can be rewritten as
\[
-\Delta {\bf v}=U^2 {\mathcal M} {\bf v}\ \hbox{in}\ \mathbb R^4, \quad {\bf v}:=(v_1,\dots,v_k) \in \mathcal{D}^{1,2}(\R^4,\R^k),
\]
with ${\mathcal M}:=\ {\mathcal Id}+2\ {\mathcal C}$ and
\[
\mathcal{C}:=\left(\begin{matrix}   \beta_{11}c_1^2& \beta_{12}c_1c_2&\dots& \beta_{1k}c_1c_k\\
\beta_{12}c_1c_2&   \beta_{22}c_2^2&  \dots&\beta_{2k}c_2c_k\\
 \vdots& \vdots&  \ddots&\vdots\\
  \beta_{1k}c_1c_k&  \beta_{2k}c_2c_k &  \dots&  \beta_{kk}c_k^2\\
  \end{matrix}\right).
\]
Let $\Lambda$ be an eigenvalue of $\mathcal M$ and $e$ an associated eigenfunction, i.e.
\[
\mathcal M e=\Lambda e.
\]
It is useful to point out that $\Lambda_\ell$ is an eigenvalue of $\mathcal M$ if and only if $\Theta_\ell:= (\Lambda_\ell-1)/ 2$ is an eigenvalue of the matrix $\mathcal C$. 
It is immediate to check that $\Theta=1$ is an eigenvalue of $\mathcal C$ whose eigenvector is $(c_1,\dots,c_k).$ We set $\Theta_1=1$, which implies $\Lambda_1=3.$

Let us consider the eigenvalue problem
\[
-\Delta v = \nu U^2 v, \qquad v \in \mathcal{D}^{1,2}(\R^4).
\]
It is well known (see \cite[Lemma A.1]{BiaEgn}) that there exists a sequence of positive eigenvalues $\{\nu_m\}_{m\in\mathbb N}$ with 
\[
1=\nu_1< 3=\nu_2 < \nu_3 < \dots < \nu_m< \nu_{m+1} < \dots \ \quad \hbox{and}\quad \nu_m \to +\infty.
\]
The role of these eigenvalues when dealing with \eqref{sist3} is clarified by the following statement.

\begin{lemma}\label{prop: base thm 2 general}
Assume that, for any $\Lambda_2,\dots,\Lambda_k$ of $\mathcal{M}$ do not coincide with any of the eigenvalues $\{\nu_m: m \in \N\}$.
%\begin{equation}\label{eigenvalues}\Lambda_\ell \not\in\{\nu_1, \dots,\nu_m,\dots\}\ \hbox{for any}\ \ell=2,\dots,k.
%\end{equation}
 Then the set of solutions to the linear system \eqref{sist3}  is
 $5-$dimensional, and is generated by
 $$\psi_{1,0}^h\mathfrak e_1  \qquad  h=0,1,\dots, 4$$
  where $\mathfrak e_1\in \mathbb R^k$ is an eigenvector associated with $\Lambda_1=3$ and the functions $\psi_{1,0}^h$ have been defined in \eqref{def psi}.
\end{lemma}

\begin{proof}
Let $\Lambda_\ell$ be an eigenvalue of the matrix $\mathcal M$ and let $\mathfrak e_\ell\in\mathbb R^k$ an associated eigenvector.
We multiply \eqref{sist3} by $\mathfrak e_\ell$ and taking into account the symmetry of the matrix $\mathcal M$ we get
\[
-\Delta (\mathfrak e_\ell\cdot{\bf v})=\Lambda_\ell\ U^2 (\mathfrak e_\ell\cdot {\bf v})\ \hbox{in}\ \mathbb R^4.
\]
Since $\Lambda_\ell \neq \nu_m$ for every $m$, we deduce that
$$\mathfrak  e_\ell\cdot{\bf v}=0\ \hbox{for any}\ \ell=2,\dots,k,$$
which implies (by the orthogonality of eigenvectors associated to different eigenvalues) that
$${\bf v}=\psi(x) \mathfrak e_1\ \hbox{for some function $\psi$ such that}\
-\Delta \psi=3 U^2 \psi\ \hbox{in}\ \mathbb R^4.$$
The claim follows then by \cite[Lemma A.1]{BiaEgn}.
\end{proof}

\begin{proof}[Conclusion of the proof of Proposition \ref{prop: on hp 2}]
As observed above, we have to prove that if $(\beta_{ij})$ is invertible and has positive entries, then the set of solutions to \eqref{sist3} is $5$-dimensional. By Lemma \ref{prop: base thm 2 general}, this amounts to show that if $(\beta_{ij})$ is invertible and has positive entries, then the eigenvalues $\Lambda_2,\dots,\Lambda_k$ of $\mathcal{M}$ are different from $\nu_1=1$, $\nu_2 = 3$, $\nu_m >3$. 

Let us argue in terms of the matrix $\mathcal{C}$. By assumption, $\mathcal C$ has positive entries. Therefore by Perron-Frobenius Theorem we deduce that the eigenvalue $\Theta_1=1$, which is associated to the eigenvector of positive elements $(c_1,\dots,c_k)$, is simple, and any other eigenvalue $\Theta_\ell$ satisfies $|\Theta_\ell|<1$.
Moreover, $0$ is not an eigenvalue of the matrix $\mathcal C$, since a straightforward computation shows that
$$ \textrm{det} \ \mathcal C =-(c_1^2\cdot\dots\cdot c_k^2)\ \textrm{det}(\beta_{ij})\not=0 $$
being $(\beta_{ij})$ invertible.
Therefore, $\Lambda_1=3$ is a simple eigenvalue, and we have that both $-1<\Lambda_\ell<3$ and $\Lambda_\ell\not =1$ for any $\ell=2,\dots,k$. This completes the proof.
%
%
%A simple example where existence is ensured is given in \cite{Ba}. We choose
% $ \beta_{ii}=\mu_i $ with $0<\mu_1<\dots<\mu_k$ and $\beta_{ij}=\beta$ if $i\not=j $ with $\beta>\mu_k.$
%It is easy to check that
% $$c_i:=\left[(\mu_i-\beta)\left(1+\beta\sum\limits_{j=1}^k{1\over \mu_j -\beta}\right)\right]^{-1/2}\ \hbox{for}\  i=1,\dots,k$$
%is a solution to system \eqref{sist2}.
%In this case assumption \eqref{betaij} is trivially satisfied. Assumption \eqref{beta-inv} is satisfied for $\beta$ large or for generic choice of $\beta$ and $\mu_i$'s.
\end{proof}

\appendix 

\section{}\label{app tools}

In this section we collect several technical lemmas. We start recalling from \cite[Lemma 3.1]{GeMuPi} the point-wise estimate of the difference of $U_{\d,\x}$ and $P_\eps U_{\d,\x}$, as well as of their derivates.

\begin{lemma}\label{lem a.1}
Let $a \in \Omega$, $r>0$, and let $\eta,\eps>0$ be small. Let $d >0$ and $\tau \in \R^N$ be such that $\eta<d<\eta^{-1}$, $|\tau|<\eta^{-1}$. Finally, let $\d= d \sqrt{\eps}$ and $\xi = a+ \d \tau$. Let us define
\[
R(x):= P_\eps U_{\d,\x}(x) -U_{\d,\x}(x) + \alpha_N \d^{\frac{N-2}{2}} H(x,\x) + \alpha_N \frac{1}{\delta^{\frac{N-2}{2}} \left(1+|\tau|^2\right)^{\frac{N-2}{2}}} \left(\frac{r \eps}{|x-a|}\right)^{N-2}.
\]
Then there exists a positive constant $C>0$ depending only on $\eta$ and on $\dist(a,\pa \Omega)$ such that for any $x \in \Omega \setminus B_{r \eps}(a)$ 
\begin{align*}
|R(x)| & \le C \delta^{\frac{N-2}{2}} \left[ \frac{\eps^{N-2} (1+\eps \delta^{1-N}) }{|x-a|^{N-2}} + \delta^2 + \left(\frac{\eps}{\delta}\right)^{N-2} \right] \\
|\pa_{\tau_i} R(x)| &\le C \d^{\frac{N}{2}} \left[ \frac{\eps^{N-2} (1+\eps \delta^{-N}) }{|x-a|^{N-2}} + \delta^2 + \frac{\eps^{N-2}}{\delta^{N-1}} \right] \\
|\pa_\d R(x)| &\le C \d^{\frac{N-4}{2}} \left[ \frac{\eps^{N-2} (1+\eps \delta^{1-N}) }{|x-a|^{N-2}} + \delta^2 + \left(\frac{\eps}{\delta}\right)^{N-2} \right].
\end{align*}
\end{lemma}

\begin{corollary}\label{corol a.2}
In the previous setting, we have
\[
\pa_{\xi_h} \left( P_\eps U_{\d,\x}(x) -U_{\d,\x}(x) \right) = -\alpha_N \d^{\frac{N-2}{2}} \pa_{\xi_h} H(x,\xi) 
 +\frac{\alpha_N (N-2)  \tau_h}{\d^{\frac{N}{2}}(1+|\tau|^2)^{\frac{N}{2}}}\left( \frac{r \eps}{|x-a|}\right)^{N-2}  + \frac{1}{\d} \pa_{\tau_h} R(x),
\]
and
\begin{multline*}
\pa_{\d} \left( P_\eps U_{\d,\x}(x) -U_{\d,\x}(x) \right) = -\alpha_N \frac{N-2}{2} \d^{\frac{N-4}{2}} H(x,\xi) \\
 +  \frac{\alpha_N(N-2)}{2\d^{\frac{N}{2}}(1+|\tau|^2)^{\frac{N-2}{2}}}\left( \frac{r \eps}{|x-a|}\right)^{N-2}  + \pa_{\d} R(x),
\end{multline*}
and in particular there exists $C>0$ depending only on $\eta$ and on $\dist(a,\pa \Omega)$ such that
\begin{align*}
| P_\eps U_{\d,\x}(x) -U_{\d,\x}(x)| & \le C \d^{\frac{N-2}2} + C \frac{ \d^{\frac32 N-3} }{|x-a|^{N-2}} \\
| P_\eps \psi_i^h(x) -\psi_i^h(x)| & \le C  \d^\frac{N-2}2 + C \frac{\d^{\frac32 N-4}}{|x-a|^{N-2}} \\
|  P_\eps \psi_i^0(x) -\psi_i^0(x)| & \le C  \d^\frac{N-4}2 + C \frac{\d^{\frac32 N-4} }{|x-a|^{N-2}},
\end{align*}
for $h=1,\dots,N$, for every $\eps>0$ small enough. In particular, if we fix a compact set $K \Subset \Omega$, then one can choose a unique constant $C>0$ for any $a \in K$.
\end{corollary}
\begin{proof}
The result is a straightforward consequence of Lemma \ref{lem a.1} and of the boundedness of $H(x,\xi)$ together with its derivatives for $x \in \Omega$, $\xi \in K'$, with $K'$ compact of $\Omega$. We also used the facts that $\d \simeq \eps^{\frac12}$ and that the derivatives commute with the projection $P_\eps$.
\end{proof}

Now we collect some estimates regarding the derivates of the bubble.

\begin{lemma}\label{lem a.7}
Let $\psi^\ell_{\d,\xi}$ ($\ell=0,\dots,N$) be defined in \eqref{def psi}. Then we have
\begin{align*}
(i) & \quad   |\psi^0_{\d,\x}| \le \frac{C}{\d} U_{\d,\x} \\
(ii) & \quad |\psi^\ell_{\d,\x}| \le \frac{C}{\d} U_{\d,\x}^{\frac{N}{N-2}} |x_\ell-\x_\ell| \quad \ell=1,\dots,N 
%\\
%(iii) & \quad |\pa_\d \psi_{\d,\x}^0| \le \frac{C}{\d^2} (N-4) U_{\d,\x} + \frac{C}{\d} U_{\d,\x}^{\frac{N}{N-2}} \\
%(iv) & \quad |\pa_\ell \psi^0_{\d,\x}| =  |\pa_0 \psi^\ell_{\d,\x}| \le \frac{C}{\d^2}     U_{\d,\x}^{\frac{N}{N-2}} |x_\ell-\x_\ell|  \quad \ell=1,\dots,N \\
%(v) & \quad |\pa_h \psi^\ell_{\d,\x}| \le \frac{C}{\d} U_{\d,\x}^{\frac{N}{N-2}}  \quad h,\ell=1,\dots,N
\end{align*}
point-wisely in $\R^N$.
\end{lemma}
\begin{proof}
The thesis is a simple consequence of the explicit expressions.
\end{proof}

The following lemma is an easy consequence of Taylor expansion.

\begin{lemma}\label{lem a.3}
For every $q>1$, there exists $C>0$ such that
\[
||a+b|^q-|a|^q| \le C (|a|^{q-1}|b| +|b|^q) \qquad \text{for every $a,b \in \R$}.
\]
Let also $f(s):= (s^+)^{q}$, with $q \ge 2$. Then there exists $C>0$ such that 
\begin{align*}
| f(a+b_1)-f(a)-f'(a) b_1 & - (f(a+b_2)-f(a)-f'(a)b_2)| \\
&\le C | |a| + |b_1| + |b_2| |^{q-2} (|b_1|+|b_2|)|b_1-b_2|
\end{align*}
for every $a,b_1,b_2 \in \R$.
\end{lemma}

A slightly more involved result is contained in the following statement.

\begin{lemma}\label{lem a.6}
Let $N=3,4$, $p=(N+2)/(N-2)$, 
\[
F(x,y):= |x|^\frac{p-3}2 x |y|^\frac{p+1}2,
\]
and
\[
\psi(h,k):= F(x+h  ,y+k)  -F(x,y)  - \pa_1 F(x,y) h - \pa_2 F(x,y) k.
\]
Then $F,\psi \in \C^{1,1}(\R^2)$, and
\begin{align*}
|\psi &(h_1,k_1) -\psi(h_2,k_2)|  \\
&\le C\left[ \frac{p-3}2 \big| |x| + |h_1| + |h_2| \big|^\frac{p-5}2 \big||y|+|k_1|+|k_2|\big|^\frac{p+1}2 \big||h_1| + |h_2|\big| \big|h_1-h_2\big| \right. \\
& \hphantom{C [ p-3} + \big| |x| + |h_1| + |h_2| \big|^\frac{p-3}2 \big| |y|  + |k_1| + |k_2| \big|^\frac{p-1}2 \big| |k_1|+|k_2| \big| \big|h_1-h_2\big|  \\
& \hphantom{C [ p-3} + \big| |x| + |h_1| + |h_2| \big|^\frac{p-3}2 \big| |y|  + |k_1| + |k_2| \big|^\frac{p-1}2 \big| |h_1|+|h_2| \big| \big|k_1-k_2\big|  \\
& \hphantom{C [ p-3}   \left. + \big| |x| + |h_1| + |h_2| \big|^{\frac{p-1}2} \big| |y| + |k_1| + |k_2| \big|^\frac{p-3}2  \big||k_1| + |k_2| \big| \big|k_1-k_2\big| \right]
\end{align*}
for every $x,y,h_1,k_1,h_2,k_2 \in \R$.
\end{lemma}
\begin{proof}
By the Lagrange intermediate value theorem
\[
|\psi(h_1,k_1) -\psi(h_2,k_2)| =  |\nabla \psi (\bar h,\bar k) \cdot (h_1-h_2,k_1-k_2)|,
\]
for some $(\bar h, \bar k)$ on the segment joining $(h_1,k_1)$ with $(h_2,k_2)$. By definition, and using again the Lagrange theorem,
\begin{align*}
|\pa_1 \psi(\bar h,\bar k)| & = |\pa_1 F(x + \bar h, y + \bar k) - \pa_1 F(x,y)| \\
& \le |\pa_{11} F(x+ \tilde h, y+\tilde k) \bar h| + |\pa_{12} F(x+ \tilde h, y + \tilde k) \bar k|,
\end{align*}
where $(\tilde h,\tilde k)$ is on the segment joining $(0,0)$ and $(\bar h,\bar k)$. A similar estimate holds for $|\pa_2 \psi (\bar h, \bar k)|$. Observing that 
\begin{align*}
|\bar h| \le |h_1| + |h_2| \quad \text{and} \quad |\bar k| \le |k_1| + |k_2|, \\
\Longrightarrow \quad  |\tilde h| \le |h_1| + |h_2| \quad \text{and} \quad |\tilde k| \le |k_1| + |k_2|,
\end{align*}
and using the explicit expression of $F$, the thesis follows.
\end{proof}

Now we collect several estimates regarding integrals of the bubbles.

\begin{lemma}\label{lem a.4}
Let $K \subset \subset \Omega$. Then, as $\d \to 0^+$, we have
\[
\int_{\Omega} U_{\d,\x}^q = \begin{cases} O(\d^{\frac{N-2}{2}q}) & \text{if }0<q<\frac{N}{N-2}, \\
C \d^{\frac{N}{2}} |\log \d| + O(\d^{\frac{N}{2}}) & \text{if }q=\frac{N}{N-2}, \\
O(\d^{N-\frac{N-2}{2} q}) & \text{if }   \frac{N}{N-2} < q <+\infty, \  q \neq  \frac{2N}{N-2}, \\
C & q = \frac{2N}{N-2} \end{cases}
\]
uniformly in $\xi \in K$, where $C>0$ denotes a constant depending only on the dimension $N$.
\end{lemma}
\begin{proof}
We focus on the cases $q \neq N/(N-2), 2N/(N-2)$ (if equality holds, one can proceed in the same way). Since $\Omega$ is bounded, there exists $R>0$ such that $\Omega \subset B_R(\xi)$. Then
\[
\int_{\Omega} U_{\d,\xi}^q \le  \int_{B_R(\xi)} U_{\d,\xi}^q = \alpha_N^q \int_{B_{R/\d}} \d^{N- \frac{N-2}{2} q} \frac{dy}{(1+|y|^2)^{\frac{N-2}{2} q}}.
\]
Now, if $q>N/(N-2)$ the last integral can be controlled with the one over $\R^N$, which is convergent, and the thesis follows. If on the other hand $q<N/(N-2)$, we have
\begin{align*}
\int_{\Omega} U_{\d,\xi}^q & \le  C \d^{N- \frac{N-2}{2} q} \left(1+ \int_1^{R/\d} r^{N-1- (N-2)q} \,dr \right) \\
& \le C \d^{N- \frac{N-2}{2} q} \cdot \frac{C}{\d^{N- (N-2)q}} \le C  \d^{\frac{N-2}{2} q}. \qedhere
\end{align*}
\end{proof}
 
In a similar way: 
 
\begin{lemma}\label{lem a.11}
Let $r>0$, $\eta \in (0,1)$, $a \in K \subset \subset \Omega$, $|\tau|<\eta^{-1}$, $\xi= a+\d \tau$, $0\le \nu_2 <N$, $\nu_1 \ge 0$, $h=1,\dots,N$, and let $(N-2)q + \nu_2-\nu_1>N$. Then, as $\d \to 0^+$, we have
\[
\int_{\Omega \setminus B_{r \d^2}(a)}  U_{\d,\xi}^q \frac{|x_h-\xi_h|^{\nu_1}}{|x-a|^{\nu_2}}  = O(\d^{N+\nu_1-\nu_2-\frac{N-2}2q}),
\]
uniformly in $a\in K$ and $|\tau|<\eta^{-1}$. \\
If $\nu_2=0$ and $q = (N+\nu_1)/(N-2)$, then 
\[
\int_{\Omega} U_{\d,\x}^q |x-\xi|^{\nu_1}  = O(\d^{\frac{N-2}{2}q} |\log \d|)
\]
as $\d \to 0^+$, uniformly in $a\in K$ and $|\tau|<\eta^{-1}$.\\
If $\nu_1=0$, $\nu_2=N$, and $(N-2)q + \nu_2>N$, then 
\[
\int_{\Omega \setminus B_{r \d^2}(a)} \frac{U_{\d,\xi}^q }{|x-a|^{N}}  = O(\d^{-\frac{N-2}2q} |\log \d|),
\]
as $\d \to 0^+$, uniformly in $a\in K$ and $|\tau|<\eta^{-1}$.
\end{lemma}
\begin{proof}
We only prove the first inequality. Arguing as in the previous lemma, we have
\begin{align*}
\int_{\Omega \setminus B_{r \d^2}(a)}  U_{\d,\xi}^q \frac{|x_h-\xi_h|^{\nu_1}}{|x-a|^{\nu_2}} & = C \int_{\frac{\Omega \setminus B_{r \d^2}(a)-a}{\d}} \d^{N+\nu_1-\frac{N-2}2q -\nu_2} \frac{ |y-\tau|^{\nu_1}}{ (1+|y|^2)^{q\frac{N-2}2} |y|^{\nu_2} }\,dy \\ & \le  C  \d^{N+\nu_1-\frac{N-2}2q -\nu_2}  \int_{B_{R/\d} \setminus B_{r\d}} \frac{ |y-\tau|^{\nu_1}}{ (1+|y|^2)^{q\frac{N-2}2} |y|^{\nu_2} }\,dy,
\end{align*} 
whence the thesis follows.
\end{proof}

\begin{lemma}\label{lem a.5}
Let $R, q_1,q_2>0$. Then there exists $C>0$ such that
\[
\int_{\Omega} U_{\d_1,\x_1}^{q_1} U_{\d_2,\x_2}^{q_2} \le C \d_1^{\frac{N-2}{2}q_1} \d_2^{\frac{N-2}{2}q_2}
+ C \d_1^{\frac{N-2}{2} q_1} \int_{\Omega} U_{\d_2,\x_2}^{q_2} + C \d_2^{\frac{N-2}{2} q_2} \int_{\Omega} U_{\d_1,\x_1}^{q_1}  
\]
for every $\xi_1,\xi_2 \in \Omega$ such that $\dist(\xi_i,\pa \Omega) > 2R$ and $|\xi_1-\xi_2| > 2R$. 
\end{lemma}
\begin{proof}
We can split the integral over $\Omega$ using the fact that   
\[
\Omega = \big( \Omega \setminus (B_R(\xi_1) \cup B_R(\xi_2))\big)  \cup B_R(\xi_1) \cup B_R(\xi_2).
\]
Then, the thesis follows using the positivity of $U_{\d,\x}$, and the fact that $0 \le U_{\d,\x} \le C \d^{\frac{N-2}2}$ in $\R^N \setminus B_R(\x)$.
\end{proof}

Similarly:

\begin{lemma}\label{lem a.9}
Let $R, q_1,q_2>0$. Then there exists $C>0$ such that
\[
\int_{\Omega} U_{\d_1,\x_1}^{q_1} U_{\d_2,\x_2}^{q_2} \frac{|x-\xi_1|^{\nu_1}}{|x-a|^{\nu_2}} \le C \d_1^{\frac{N-2}{2}q_1} \d_2^{\frac{N-2}{2}q_2}
+ C \d_1^{\frac{N-2}{2} q_1} \int_{\Omega} U_{\d_2,\x_2}^{q_2} + C \d_2^{\frac{N-2}{2} q_2} \int_{\Omega} U_{\d_1,\x_1}^{q_1}   \frac{|x-\xi_1|^{\nu_1}}{|x-a|^{\nu_2}}
\]
for every $\xi_1,\xi_2,a \in \Omega$ such that $\dist(\xi_i,\pa \Omega) > 2R$, $|\xi_1-\xi_2| > 2R$, $|\xi_1-a|<R/2$, $\nu_1 \ge 0$ and $0 \le \nu_2 < N$. 
\end{lemma}
%\begin{proof}
%We can argue as in the proof of Lemma \ref{lem a.5}.\end{proof}

%
%\begin{lemma}\label{lem a.8}
%Let $K \subset \subset \Omega$, $\nu>0$, and let $q>(N+\nu)/(N-2)$. We have
%\[
%\int_{\Omega} U_{\d,\x}^q |x-\xi|^\nu  = O(\d^{N+\nu- \frac{N-2}{2}q})
%\]
%as $\d \to 0^+$, uniformly in $\xi \in K$. If $q = (N+\nu)/(N-2)$, then 
%\[
%\int_{\Omega} U_{\d,\x}^q |x-\xi|^\nu  = O(\d^{N+\nu- \frac{N-2}{2}q} |\log \d|)
%\]
%as $\d \to 0^+$, uniformly in $\xi \in K$.
%\end{lemma}
%\begin{proof}
%We argue as in Lemma \ref{lem a.4}:
%\begin{align*}
%\int_{\Omega} U_{\d,\xi}^q |x-\xi|^{\nu} &  \le \alpha_N^q \int_{\R^N} \d^{N- \frac{N-2}{2} q + \nu} \frac{|y|^\alpha dy}{(1+|y|^2)^{\frac{N-2}{2} q}}.
%\end{align*}
%As $q>(N+\nu)/(N-2)$, the integral is convergent, and the thesis follows. The second part of the lemma can be proved in the same way.
%\end{proof}

%\begin{lemma}
%Let $K \subset \subset \Omega$, $0< \nu <N$, and let $q>(N-\nu)/(N-2)$. Then, as $\d \to 0^+$, we have
%\[
%\int_{\Omega_\eps} U_{\d,\xi}^q \frac{dx}{|x-a|^\nu}  = O(\d^{N-\nu-\frac{N-2}2q}),
%\]
%uniformly in $a,\xi \in K$.
%\end{lemma}
%
%\begin{lemma}\label{lem a.11}
%Let $K \subset \subset \Omega$, $0\le \nu_2 <N$, $\nu_1 \ge 0$, $h=1,\dots,N$, and let $(N-2)q + \nu_2-\nu_1>N$. Then, as $\d \to 0^+$, we have
%\[
%\int_{\Omega_\eps} U_{\d,\xi}^q |x_h-\xi_h|^{\nu_1}\frac{dx}{|x-a|^{\nu_2}}  = O(\d^{N+\nu_1-\nu_2-\frac{N-2}2q}),
%\]
%uniformly in $a,\xi \in K$.
%\end{lemma}

\section{}\label{app step 1 linear}

In this appendix we provide the details for the estimates \eqref{estimates step 1}. 

\smallskip

\noindent \textbf{Estimate of ($I$):} integrating by parts we deduce that
\begin{multline}\label{eq app 1 1}
\delta_{i,n}^2\int_{\Omega_n}\nabla (P_n \psi_{i,n}^l) \cdot  \nabla (P_n \psi_{i,n}^k)  = \d_{i,n}^2 \int_{\Omega_n} p U_{i,n}^{p-1} \psi_{i,n}^k (P_n \psi_{i,n}^l) \\
 = \d_{i,n}^2  \int_{\Omega_n} p U_{i,n}^{p-1} \psi_{i,n}^k \psi_{i,n}^l + \d_{i,n}^2  \int_{\Omega_n} p U_{i,n}^{p-1} \psi_{i,n}^k (P_n \psi_{i,n}^l-U_{i,n}^l).
\end{multline}
Due to Corollary \ref{corol a.2}, it is not difficult to show that the last term tends to $0$ as $n \to \infty$. We consider here in full details the case $l,k =1,\dots,N$ (the others can be treated similarly): first, we observe that Corollary \ref{corol a.2} is applicable since $\{\xi_{i,n}\} \subset K$ for some compact set $K \Subset \Omega$, at least for $\eta>0$ small enough. Then, using also Lemmas \ref{lem a.7} and \ref{lem a.11}, and recalling that $\eps_n \simeq \d_{i,n}^2$ since $\eta$ is fixed, we deduce that
\begin{align*}
\d_{i,n}^2 \Bigg| \int_{\Omega_n}  U_{i,n}^{p-1} &  \psi_{i,n}^k (P_n \psi_{i,n}^l-U_{i,n}^l) \Bigg| \le C \d_{i,n} \int_{\Omega_n} U_{i,n}^{\frac{4+N}{N-2}} |x_k- \xi^{i,n}_k| \left( \d_{i,n}^\frac{N-2}2 + \frac{ \d_{i,n}^{\frac32N-4}}{|x-a_i|^{N-2}} \right) \\
& = C  \d_{i,n}^\frac{N}2  \int_{\Omega_n} U_{i,n}^{\frac{4+N}{N-2}} |x_k- \xi^{i,n}_k|  + C  \d_{i,n}^{\frac32N-3}  \int_{\Omega_n} U_{i,n}^{\frac{4+N}{N-2}} \frac{|x_k- \xi^{i,n}_k|}{|x-a_i|^{N-2}} \\
& \le C \d_{i,n}^{\frac{N}2 + N+1-\frac{4+N}2} + C \d_{i,n}^{\frac32N-3+N+1- \frac{4+N}2-N+2} = o(1)
\end{align*}
as $n \to \infty$.
%To apply Lemma \ref{lem a.1}, we point out that the Robin function $H(x,\xi_{i,n})$ is uniformly bounded for $x \in \Omega_n$ and $\xi_{i,n}$ detached from $\pa \Omega$, uniormly in $n$ (it is not restrictive to suppose that this condition is satisfied; indeed, it is sufficient to choose $\eta>0$ small enough from the beginning). 

Coming back to \eqref{eq app 1 1}, using the explicit expression of $U_{i,n}$ and $\psi_{i,n}^k$, it is not difficult to check that 
\[
\delta_{i,n}^2\int_{\Omega_n}\nabla (P_n \psi_{i,n}^l) \cdot  \nabla (P_n \psi_{i,n}^k)= \sigma_{lk}+o(1),
\]
where $o(1) \to 0$ as $n \to \infty$, and the values $\sigma_{lk}$ are defined in \eqref{sigma lk}.
%\[
%\sigma_{lk} = \begin{cases} 0 & \text{if $l \neq k$} \\  p \alpha_N^{p+1} (N-2)^2 \int_{\R^N} \frac{y_l^2}{(1+|y|^2)^{N+2}}\,dy & \text{if $k=l \ge 1$} \\  p  \alpha_N^{p+1} \left( \frac{N-2}{2}\right)^2 \int_{\R^N} \frac{(|y|^2-1)^2}{(1+|y|^2)^{N+2}}\,dy & \text{if $k=l=0$}.
%\end{cases}
%\] 
Therefore, as $n \to \infty$ 
\begin{equation}\label{develop w}
\begin{split}
\d_{i,n}^2 \|w_{i,n}\|_{H_0^1(\Omega_n)}^2 & = \sum_{l,k=0}^N c_{i,n}^l c_{i,n}^k \delta_{i,n}^2\int_{\Omega_n}\nabla (P_n \psi_{i,n}^l) \cdot  \nabla (P_n \psi_{i,n}^k) \\ &=\sum_{l,k=0}^N c_{i,n}^l c_{i,n}^k (\sigma_{lk} +o(1)) =  \sum_{l= 0}^N (c_{i,n}^l)^2 \sigma_{ll} + o(1) \sum_{l,k=0}^N c_{i,n}^l c_{i,n}^k.
\end{split}
\end{equation}
\textbf{Estimate of ($II$):} we aim at proving the second equality in \eqref{estimates step 1}. 
Since $\phi_{i,n} \in K_{i,n}^\perp$, we know that $\langle \phi_{i,n}, P_n \psi_{i,n}^l \rangle_{H_0^1} = 0$, i.e.
\[
0 = \int_{\Omega_n} \nabla (P_n \psi_{i,n}^l) \cdot \nabla \phi_{i,n} = p\int_{\Omega_n}  U_{i,n}^{p-1} \phi_{i,n} \psi_{i,n}^l.
\] 
As a consequence
\begin{equation}\label{265}
\begin{split}
|(II)| \le \delta_{i,n}^2\Bigg[ &\underbrace{\int_{\Omega_n} |(P_n U_{i,n})^{\frac4{N-2}} - U_{i,n}^{\frac4{N-2}}| |\phi_{i,n}||w_{i,n}|}_{=:(II.1)}  + \sum_{l=0}^N c_{i,n}^l \underbrace{\int_{\Omega_n} U_{i,n}^\frac4{N-2}|\phi_{i,n}||P_n \psi_{i,n}^l- \psi_{i,n}^l|}_{=:(II.2)} \Bigg].
 \end{split}
\end{equation}
Concerning the first integral, by the H\"older and the Sobolev inequalities
\begin{equation}\label{3051}
\begin{split}
|(II.1)| &\le |(P_n U_{i,n})^\frac4{N-2}-U_{i,n}^\frac4{N-2}|_{\frac{N}2} |\phi_{i,n}|_{2^*} |w_{i,n} |_{2^*} \\
& \le  C|(P_n U_{i,n})^\frac4{N-2}-U_{i,n}^\frac4{N-2}|_{\frac{N}2} \|\phi_{i,n} \| \|w_{i,n} \|. 
\end{split}
\end{equation}
We recall that $\|\phi_{i,n}\|=1$ for every $n$. Moreover, since $0 \le P_n U_{i,n} \le U_{i,n}$ in $\Omega_n$ by the maximum principle, applying Lemma \ref{lem a.3} we have 
\begin{align*}
|(P_n U_{i,n})^\frac4{N-2}-U_{i,n}^\frac4{N-2}|_{\frac{N}2} & \le |U_{i,n}^{\frac4{N-2}-1} (P_n U_{i,n}-U_{i,n})|_{\frac{N}2} + |(P_n U_{i,n}-U_{i,n})^\frac4{N-2}|_{\frac{N}2}  \\
& \le 2|U_{i,n}^{\frac4{N-2}-1} (P_n U_{i,n}-U_{i,n})|_{\frac{N}2} = 2 \left( \int_{\Omega_n} U_{i,n}^{\frac{(6-N)N}{2(N-2)}} |P_n U_{i,n}-U_{i,n}|^{\frac{N}2} \right)^\frac2N.
\end{align*}
Thus, using Corollary \ref{corol a.2}, the fact that $\eps_n \simeq \d_{i,n}^2$ (see ansatz \eqref{asympt expansion}), and Lemmas \ref{lem a.4} and \ref{lem a.11}, we obtain for $N=3$
\begin{align*}
|(P_n U_{i,n})^\frac4{N-2} & -U_{i,n}^\frac4{N-2}|_{\frac{N}2} \\
&\le C \d_{i,n}^{\frac{N-2}2} \left( \int_{\Omega_n} U_{i,n}^{\frac{(6-N)N}{2(N-2)}}  \right)^\frac2N 
 + C \d_{i,n}^{\frac32N-3} \left( \int_{\Omega_n} \frac{U_{i,n}^{\frac{(6-N)N}{2(N-2)}}}{|x-a_i|^{\frac{N-2}2N}}  \right)^\frac2N  = C \d_{i,n}^{N-2},
\end{align*}
and similarly for $N=4$
\[
|(P_n U_{i,n})^\frac4{N-2}  -U_{i,n}^\frac4{N-2}|_{\frac{N}2} \le C \d_{i,n}^{N-2} |\log \d_{i,n}|^\frac2N.
\] 
Thus, coming back to \eqref{3051}, we proved that $|(II.1)| = o(1) \|w_{i,n}\|_{H_0^1(\Omega_n)}$.

Let us consider now $|(II.2)|$. By the H\"older and the Sobolev inequalities
\[
|(II.2)| \le C\|\phi_{i,n}\| |U_{i,n}^{\frac4{N-2}} (P_n \psi_{i,n}^l-\psi_{i,n}^l)|_{\frac{2N}{N+2}}.
\]
The right hand side can be controlled using Corollary \ref{corol a.2} and Lemmas \ref{lem a.4} and \ref{lem a.11}. We focus on the case $l=1,\dots,N$, which is completely analogue with respect to $l=0$, and we compute 
\begin{align*}
|(II.2)| &\le \left( \int_{\Omega_n} U_{i,n}^\frac{8N}{(N-2)(N+2)}    |P_n \psi_{i,n}^l-\psi_{i,n}^l|^\frac{2N}{N+2} \right)^\frac{N+2}{2N} \\
& \le C\left(\int_{\Omega_n} U_{i,n}^\frac{8N}{(N-2)(N+2)}  \left(\d_{i,n}^\frac{N-2}2 + \frac{\d_{i,n}^{\frac32N-4}   }{|x-a_i|^{N-2}} \right)^\frac{2N}{N+2}  \right)^{\frac{N+2}{2N}} \\
& \le C \d_{i,n}^\frac{N-2}2 \left( \int_{\Omega_n} U_{i,n}^\frac{8N}{(N-2)(N+2)} \right)^{\frac{N+2}{2N}} + C \d_{i,n}^{\frac32N-4}  \left( \int_{\Omega_n} \frac{U_{i,n}^\frac{8N}{(N-2)(N+2)}}{|x-a_i|^{\frac{2N(N-2)}{N+2}} }\right)^{\frac{N+2}{2N}} \le C \d_{i,n}^{N-3}. 
\end{align*}
Coming back to \eqref{265} (and recalling that $N=3,4$), we proved the validity of the second estimate in \eqref{estimates step 1}, as desired.\\
\textbf{Estimate of $(III)$:} we show that $(III) =o(\d_{i,n}^2)$ as $n \to \infty$, focusing here on the case $l=0$; the one $l=1,\dots,N$ is, once again, completely analogue. Recalling that $0 \le P_n U_{i,n} \le U_{i,n}$, using the H\"older and the Sobolev inequalities we compute
\begin{equation}\label{262}
\begin{split}
\d_{i,n}^2\Bigg| \int_{\Omega_n} &(P_n U_{i,n})^\frac{p-3}2  (P_n U_{j,n})^\frac{p+1}2    \phi_{i,n}  (P_n \psi_{i,n}^l) \Bigg|  \\
& = \d_{i,n}^2 \Bigg| \int_{\Omega_n} (P_n U_{i,n})^\frac{p-3}2  (P_n U_{j,n})^\frac{p+1}2      \phi_{i,n} \left[ \psi_{i,n}^l + \left(P_n \psi_{i,n}^l - \psi_{i,n}^l \right) \right]  \Bigg| \\
& \le \d_{i,n}^2 \int_{\Omega_n} U_{i,n}^\frac{p-3}2 U_{j,n}^\frac{p+1}2    |\phi_{i,n}| |\psi_{i,n}^l| + \d_{i,n}^2 \int_{\Omega_n} U_{i,n}^\frac{p-3}2 U_{j,n}^\frac{p+1}2    |\phi_{i,n}| |P_n \psi_{i,n}^l -\psi_{i,n}^l| \\
& \le C\d_{i,n}^2  \underbrace{\|\phi_{i,n}\| \left| U_{i,n}^\frac{p-3}2 U_{j,n}^\frac{p+1}2  |\psi_{i,n}^l| \right|_{\frac{2N}{N+2}}}_{=:(III.1)}  +  C \d_{i,n}^2 \underbrace{\|\phi_{i,n}\| \left| U_{i,n}^\frac{p-3}2 U_{j,n}^\frac{p+1}2  |P_n \psi_{i,n}^l-\psi_{i,n}^l| \right|_{\frac{2N}{N+2}} }_{=:(III.2)}
\end{split}
\end{equation}
In order to estimate $(III.2)$, we use Corollary \ref{corol a.2}, the ansatz \eqref{asympt expansion}, and the fact that $\|\phi_{i,n}\| = 1$:
\begin{equation}\label{791}
\begin{split}
(III.2)  & \le C \left| U_{i,n}^\frac{p-3}2 U_{j,n}^\frac{p+1}2  \left( \d_{i,n}^\frac{N-4}2 + \frac{\d_{i,n}^{\frac32N-4} }{|x-a_i|^{N-2}} \right) \right|_{\frac{2N}{N+2}} \\
  & \le C \d_{i,n}^\frac{N-4}2 \left( \int_{\Omega_n} U_{i,n}^\frac{2N(4-N)}{(N-2)(N+2)} U_{j,n}^\frac{2N^2}{(N-2)(N+2)} \right)^\frac{N+2}{2N}  \\
  & + C \d_{i,n}^{\frac32N-4}\left( \int_{\Omega_n} \frac{U_{i,n}^\frac{2N(4-N)}{(N-2)(N+2)} U_{j,n}^\frac{2N^2}{(N-2)(N+2)}}{|x-a_i|^{\frac{2N(N-2)}{(N+2)}  }} \right)^\frac{N+2}{2N}.
\end{split}
\end{equation}
By Lemmas \ref{lem a.4} and \ref{lem a.5},
\begin{align*}
\d_{i,n}^\frac{N-4}2  & \left( \int_{\Omega_n} U_{i,n}^\frac{2N(4-N)}{(N-2)(N+2)} U_{j,n}^\frac{2N^2}{(N-2)(N+2)} \right)^\frac{N+2}{2N} \\
& \le C \d_{i,n}^\frac{N-4}2 \left( \d_{i,n}^\frac{N(4-N)}{N+2} \d_{j,n}^\frac{N^2}{N+2} +  \d_{i,n}^\frac{N(4-N)}{N+2} \d_{j,n}^\frac{2N}{N+2} + C \d_{j,n}^\frac{N^2}{N+2} \int_{\Omega_n} U_{i,n}^\frac{2N(4-N)}{(N-2)(N+2)} \right)^{\frac{N+2}{2N}},
\end{align*}
and discussing separately the case $N=3$ and $N=4$, it is not difficult to check that in both cases
\[
\d_{i,n}^\frac{N-4}2   \left( \int_{\Omega_n} U_{i,n}^\frac{2N(4-N)}{(N-2)(N+2)} U_{j,n}^\frac{2N^2}{(N-2)(N+2)} \right)^\frac{N+2}{2N}  \le C \d_{j,n} =o(1)
\]
as $n \to \infty$. The second term on the right hand side in \eqref{791} can be controlled using Lemmas \ref{lem a.4} and \ref{lem a.9}, in a similar way:
\begin{align*}
\d_{i,n}^{\frac32N-4} & \left( \int_{\Omega_n} \frac{U_{i,n}^\frac{2N(4-N)}{(N-2)(N+2)} U_{j,n}^\frac{2N^2}{(N-2)(N+2)}}{|x-a_i|^{\frac{2N (N-2)}{(N+2)}  }} \right)^\frac{N+2}{2N} \\
& \le C \d_{i,n}^{\frac32N-4} \left(  \d_{i,n}^\frac{N(4-N)}{N+2} \d_{j,n}^\frac{N^2}{N+2} +  \d_{i,n}^\frac{N(4-N)}{N+2} \d_{j,n}^\frac{2N}{N+2} + C \d_{j,n}^\frac{N^2}{N+2} \d_{i,n}^{N- \frac{N(4-N)}{N+2}- \frac{2N(N-2)}{N+2}} \right)^\frac{N+2}{2N} = o(1)
\end{align*}
as $n \to \infty$. Therefore, $(III.2) \to 0$ as $n \to \infty$. Concerning the term $(III.1)$, by Lemmas \ref{lem a.7}, \ref{lem a.4} and \ref{lem a.5}, we obtain
\begin{align*}
|(III.1)| & \le C \d_{i,n}^{-1} \left| U_{i,n}^\frac{p-1}2 U_{j,n}^\frac{p+1}2 \right|_{\frac{2N}{N+2}} \\
& \le C \d_{i,n}^{-1} \left( \int_{\Omega_n} U_{i,n}^\frac{4N}{(N-2)(N+2)} U_{j,n}^\frac{2N^2}{(N-2)(N+2)}  \right)^\frac{N+2}{2N}  \\
& C \d_{i,n}^{-1} \left( \d_{i,n}^\frac{2N}{N+2} \d_{j,n}^\frac{N^2}{N+2} + \d_{i,n}^\frac{2N}{N+2} \d_{j,n}^\frac{2N}{N+2} + \d_{i,n}^\frac{2N}{N+2} \d_{j,n}^\frac{N^2}{N+2} \right)^\frac{N+2}{2N} \le C \d_{j,n}
= o(1)
\end{align*}
as $n \to \infty$.

To sum up we showed that $(III.1)$ and $(III.2)$ tend to $0$ as $n \to \infty$, and coming back to \eqref{262} we infer that $(III) = o(\d_{i,n}^2)$, as desired. \\
\textbf{Estimate of $(IV)$:} this is similar to the estimate of $(III)$, and this time we focus on $l=1,\dots,N$. We have
\begin{equation}\label{792}
\begin{split}
\d_{i,n}^2\Bigg| \int_{\Omega_n} &(P_n U_{i,n})^\frac{p-1}2  (P_n U_{j,n})^\frac{p-1}2    \phi_{i,n}  (P_n \psi_{i,n}^l) \Bigg|  \\
& \le C\d_{i,n}^2  \underbrace{\|\phi_{i,n}\| \left| U_{i,n}^\frac{p-1}2 U_{j,n}^\frac{p-1}2  |\psi_{i,n}^l| \right|_{\frac{2N}{N+2}}}_{=:(IV.1)}  +  C \d_{i,n}^2 \underbrace{\|\phi_{i,n}\| \left| U_{i,n}^\frac{p-1}2 U_{j,n}^\frac{p-1}2  |P_n \psi_{i,n}^l-\psi_{i,n}^l| \right|_{\frac{2N}{N+2}} }_{=:(IV.2)}
\end{split}
\end{equation}
The term $(IV.2)$ can be controlled using Corollary \ref{corol a.2}:
\begin{equation}\label{793}
\begin{split}
(IV.2)  & \le C \left| U_{i,n}^\frac{p-1}2 U_{j,n}^\frac{p-1}2  \left( \d_{i,n}^\frac{N-2}2 + \frac{\d_{i,n}^{\frac32N-4} }{|x-a_i|^{N-2}} \right) \right|_{\frac{2N}{N+2}} \\
  & \le C \d_{i,n}^\frac{N-2}2 \left( \int_{\Omega_n} U_{i,n}^\frac{4N}{(N-2)(N+2)} U_{j,n}^\frac{4N}{(N-2)(N+2)} \right)^\frac{N+2}{2N} \\
  & + C \d_{i,n}^{\frac32N-4}\left( \int_{\Omega_n} \frac{U_{i,n}^\frac{4N}{(N-2)(N+2)} U_{j,n}^\frac{4N}{(N-2)(N+2)}}{|x-a_i|^{\frac{2N(N-2)}{(N+2)}  }} \right)^\frac{N+2}{2N}.
\end{split}
\end{equation}
By Lemmas \ref{lem a.4} and \ref{lem a.5}, the first integral on the right hand side gives
\begin{align*}
\d_{i,n}^\frac{N-2}2 &  \left( \int_{\Omega_n} U_{i,n}^\frac{4N}{(N-2)(N+2)} U_{j,n}^\frac{4N}{(N-2)(N+2)} \right)^\frac{N+2}{2N}  \le C \d_{i,n}^\frac{N-2}2 \left( \d_{i,n}^\frac{2N}{N+2} \d_{j,n}^\frac{2N}{N+2} \right)^\frac{N+2}{2N} = o(1),
\end{align*}
and, by Lemmas \ref{lem a.4} and \ref{lem a.9}, the second integral gives
\begin{align*}
 \d_{i,n}^{\frac32N-4} & \left( \int_{\Omega_n} \frac{U_{i,n}^\frac{4N}{(N-2)(N+2)} U_{j,n}^\frac{4N}{(N-2)(N+2)}}{|x-a_i|^{\frac{2N (N-2)}{(N+2)}  }} \right)^\frac{N+2}{2N} \\
 & \le C \d_{i,n}^{\frac32N-4} \left(  \d_{i,n}^\frac{2N}{N+2} \d_{j,n}^\frac{2N}{N+2} +  \d_{i,n}^\frac{N(4-N)}{N+2} \d_{j,n}^\frac{2N}{N+2} \right)^\frac{N+2}{2N}= o(1) 
\end{align*}
 as $n \to \infty$. As far as $(IV.1)$ is concerned, using Lemmas \ref{lem a.7}, \ref{lem a.11} and \ref{lem a.9}, we infer that
\begin{align*}
|(IV.1)| & \le C \d_{i,n}^{-1} \left| U_{i,n}^{\frac{p-1}2 + \frac{N}{N-2}} U_{j,n}^\frac{p-1}2 |x_l - \xi_{i,l}^n|\right|_{\frac{2N}{N+2}} \\
& \le C \d_{i,n}^{-1} \left( \int_{\Omega_n} U_{i,n}^\frac{2N}{N-2} U_{j,n}^\frac{4N}{(N-2)(N+2)} |x_l - \xi_{i,l}^n|^\frac{2N}{N+2} \right)^\frac{N+2}{2N}  \\
& \le C \d_{i,n}^{-1} \left( \d_{i,n}^N \d_{j,n}^\frac{2N}{N+2} + \d_{j,n}^\frac{2N}{N+2} \d_{i,n}^{N-N+\frac{2N}{N+2}} \right)^\frac{N+2}{2N} 
= o(1)
\end{align*}
as $n \to \infty$. Altogether, we proved that $(IV.1)$ and $(IV.2)$ tend to $0$ as $n \to \infty$, so that $(IV) = o(\d_{i,n}^2)$ as $n \to \infty$.

This completes the proof of the validity of the estimates \eqref{estimates step 1}.

\section{}\label{app est step 2}

In this appendix we show that \eqref{eq phi tilde} gives \eqref{estimate step 2}. At first, we notice that
\[
\begin{split}
p \d_{\kappa,n}^2 \int_{\tilde \Omega_{\kappa,n}} (\widehat{P_n U_{\kappa,n}})^{p-1}  & \tilde \phi_{\kappa,n} \psi 
= p \d_{\kappa,n}^2 \int_{\tilde \Omega_{\kappa,n}} U_{\kappa,n}^{p-1}(\x_{\kappa,n}+\d_{\kappa,n} \, \cdot) \tilde \phi_{\kappa,n} \psi\\
& + p \d_{\kappa,n}^2 \int_{\tilde \Omega_{\kappa,n}} [(\widehat{P_n U_{\kappa,n}})^{p-1} - U_{\kappa,n}^{p-1}(\x_{\kappa,n}+\d_{\kappa,n} \, \cdot)] \tilde \phi_{\kappa,n} \psi.
\end{split}
\] 
Arguing as in Appendix \ref{app step 1 linear}, estimate of term $(II.1)$, it can be proved that the last integral on the right hand side tends to $0$ as $n \to \infty$. Moreover, the first integral can be explicitly computed: 
\[
U_{\kappa,n}^{p-1}(\x_{\kappa,n}+\d_{\kappa,n} y)  \psi(y) = \alpha_N^{p-1} \left( \frac{1}{1+|y|^2}\right)^2\psi(y) = U_{0,1}^{p-1}(y) \psi(y) \in L^{\frac{2N}{N+2}}(\R^N)
\]
(where we recall that $U_{1,0}$ denotes the standard bubble with $\delta=1$ and $\xi=0$), and hence by weak convergence $\tilde \phi_{\kappa,n} \wc \tilde \phi_\kappa$ in $L^{2^*}(\R^N)$
\begin{equation}\label{A_n}
\begin{split}
p \d_{\kappa,n}^2 \int_{\tilde \Omega_{\kappa,n}} (\widehat{P_n U_{\kappa,n}})^{p-1}   \tilde \phi_{\kappa,n} \psi 
&= p \d_{\kappa,n}^2 \int_{\tilde \Omega_{\kappa,n}} U_{\kappa,n}^{p-1}(\x_{\kappa,n}+\d_{\kappa,n} \, \cdot) \tilde \phi_{\kappa,n} \psi + o(1) \\
& = p \int_{\tilde \Omega_{\kappa,n}} U_{1,0}^{p-1} \tilde \phi_{\kappa,n} \psi  + o(1) = p \int_{\R^N} U_{1,0}^{p-1} \tilde \phi_{i} \psi  +o(1)
\end{split}
\end{equation}
as $n \to \infty$.

%Therefore
%\[
%\Bigg|\d_{\kappa,n}^2 \int_{\tilde \Omega_{\kappa,n}} U_{\kappa,n}^{p-1} (\x_{\kappa,n}+\d_{\kappa,n} \, \cdot) \tilde \phi_{\kappa,n} \psi \Bigg|  \le C \d_{\kappa,n}^2 \d_{\kappa,n}^2 \int_{\supp \psi}  |\psi| |\tilde \phi_{\kappa,n}| \to 0
%\]
%as $n \to \infty$.

Now we show that the remaining terms on the right hand side of \eqref{eq phi tilde} tend to $0$ as $n \to \infty$. We start observing that, if $j \neq \kappa$, for any compact set $K \subset \R^N$ there exists $C>0$ (depending on $K$) such that for sufficiently large $n$
\begin{equation}\label{5ott1}
\begin{split}
\inf_{y \in K}  |\x_{\kappa,n}+ \d_{\kappa,n} y-\xi_{j,n}| & \ge |\x_{\kappa,n}-\xi_{j,n}| - C   \d_{\kappa,n}\\
& \ge \frac12 |a_\kappa-a_j| - C \d_{\kappa,n} \ge \frac14 |a_\kappa-a_j|.
\end{split}
\end{equation}
Thus, for any $j \neq \kappa$
\[
\begin{split}
\d_{\kappa,n}^2 & \Bigg| \int_{\tilde \Omega_{\kappa,n}}  ( \widehat{P_n U_{\kappa,n}})^{\frac{p-3}{2}}  (\widehat{P_n U_{j,n}})^{\frac{p+1}{2}}  \tilde \phi_{\kappa,n} \psi \Bigg|  \\
&\le C \d_{\kappa,n}^2 \int_{\tilde \Omega_{\kappa,n}} U_{\kappa,n}^{\frac{p-3}{2}}(\xi_{\kappa,n} + \d_{\kappa,n} \, \cdot) U_{j,n}^{\frac{p+1}{2}} (\xi_{\kappa,n} + \d_{\kappa,n} \, \cdot) |\tilde \phi_{\kappa,n}||\psi| \\
& \le C \d_{\kappa,n}^2
% % |\tilde \phi_{i,n}|_{L^{2^*}(\tilde \Omega_{i,n})}  |\psi|_{L^{2^*}(\tilde \Omega_{i,n})} 
  \int_{\supp \psi} \left( \frac{\d_{j,n}}{\d_{j,n}^2 + |\x_{\kappa,n}+ \d_{\kappa,n} y-\xi_{j,n}|^2}\right)^\frac{N}{2} \left( \frac{\d_{\kappa,n}^{-1}}{1 + |y|^2}\right)^\frac{4-N}{2} |\tilde \phi_{\kappa,n}(y)| |\psi(y)|\,dy \\
&   \le C\d_{\kappa,n}^{\frac{N}{2}} \d_{j,n}^{\frac{N}{2}} \int_{\supp \psi} |\tilde \phi_{\kappa,n}(y)| |\psi(y)|\,dy   \le C \d_{\kappa,n}^{\frac{N}{2}} \d_{j,n}^{\frac{N}{2}} |\tilde \phi_{\kappa,n}|_{2^*} |\psi|_{\frac{2N}{N+2}} = o(1)
\end{split}
\]
as $n \to \infty$. 

Similarly, always for any $j \neq i$
\begin{align*}
\d_{\kappa,n}^2 & \Bigg|  \int_{\tilde \Omega_{\kappa,n}} (\widehat{P_n U_{\kappa,n}})^{\frac{p-1}{2}} (\widehat{P_n U_{j,n}})^{\frac{p-1}{2}}\tilde \phi_{j,n} \psi\Bigg| \\
& \le C \d_{\kappa,n}^2 \int_{\tilde \Omega_{\kappa,n}} \frac{\d_{j,n}}{\d_{j,n}^2 + |\x_{\kappa,n}+ \d_{\kappa,n} y-\xi_{j,n}|^2} \frac{\d_{\kappa,n}^{-1}}{1+ | y|^2}  |\tilde \phi_{j,n}(y)| |\psi(y)|\,dy \\
& \le C \d_{\kappa,n} \d_{j,n} \int_{\supp \psi} |\tilde \phi_{j,n}(y)| |\psi(y)|\,dy = o(1)
\end{align*}
as $n \to \infty$.

Finally, also the last term on the right hand side in \eqref{eq phi tilde} tends to $0$, since $\|\tilde h_{\kappa,n}\|_{H_0^1(\tilde \Omega_{\kappa,n})} = \|h_{\kappa,n}\|_{H_0^1(\Omega_n)}$, $\|\tilde w_{\kappa,n}\|_{H_0^1(\tilde \Omega_{\kappa,n})} = \|w_{\kappa,n}\|_{H_0^1(\Omega_n)}$, both $h_{\kappa,n}$ and $w_{\kappa,n}$ converge to $0$ strongly in $H_0^1(\Omega_n)$ ($h_{\kappa,n}$ by assumption, $w_{\kappa,n}$ by step 1), and $\tilde{\phi}_{\kappa,n} \wc \tilde{\phi}_i$ weakly.

\section{}\label{app estimate nonlinear}

In this appendix we prove the validity of the estimates \eqref{estimate nonlinear} and \eqref{674}. In order to ease the notation, we write $P$ and $U_i$ instead of $P_\eps$ and $U_{\d_i,\x_i}$.

\smallskip

We start with the proof of \eqref{estimate nonlinear}.

\noindent \textbf{Estimate of $|\tilde{R}^i_{\mf{d},\bs{\tau},\eps}|_{L^\frac{2N}{N+2}(\Omega_\eps)}$}. By Lemma \ref{lem a.3} and the fact that $0 \le P U_i \le U_i$,
\begin{align*}
|P U_i^p- U_i^p| &\le C \left( |U_i|^{p-1} |P U_i - U_i| + |P U_i - U_i|^p \right)  \le C |U_i|^{p-1} |P U_i - U_i|.
\end{align*}
Therefore, recalling that $\eps \simeq \d_i^2$ by \eqref{asympt expansion} (we fixed $\eta$ from the beginning), Corollary \ref{corol a.2} and Lemmas \ref{lem a.4} and \ref{lem a.11} give
\begin{equation}\label{171}
\begin{split}
\Bigg(\int_{\Omega_\eps} |P U_i^p & - U_i^p|^\frac{2N}{N+2}\Bigg)^\frac{N+2}{2N} \le 
C \Bigg( \int_{\Omega_\eps}  U_i^{\frac{8N}{(N-2)(N+2)}} |P U_i - U_i|^\frac{2N}{N+2} \Bigg)^\frac{N+2}{2N} \\
& \le  C \left(\int_{\Omega_\eps} U_i^{\frac{8N}{(N-2)(N+2)}} \left( \d_i^{\frac{N-2}{2}} + \frac{\d_i^{\frac{3}{2}(N-2)}}{|x-a_i|^{N-2}} \right)^{\frac{2N}{N+2}}\right)^{\frac{N+2}{2N}} \\
& \le C  \left( \d_{i}^{\frac{2N(N-2)}{N+2}} + \d_i^{\frac{3N(N-2)}{N+2} + \frac{N(2-N)}{N+2}   } \right)^\frac{N+2}{2N} \le C \d_{i}^{N-2} \le C \eps^\frac{N-2}2.
\end{split}
\end{equation}
We also have, by Lemmas \ref{lem a.4} and \ref{lem a.5}
\begin{align*}
\left| (P U_j)^{\frac{p+1}{2}} (P U_i)^{\frac{p-1}{2}} \right|_{\frac{2N}{N+2}} & \le    \left| U_j^{\frac{p+1}{2}}  U_i^{\frac{p-1}{2}} \right|_{\frac{2N}{N+2}} \\
& = \left(\int_{\Omega_\eps} U_j^\frac{2N^2}{(N+2)(N-2)} U_i^\frac{4N}{(N+2)(N-2)} \right)^\frac{N+2}{2N} \le C \d_i \d_j \le C \eps.
\end{align*}
This estimate and \eqref{171} imply that, for $\eps>0$ small enough, there exists $C>0$ such that
\begin{equation}\label{stima tilde R}
|\tilde{R}^i_{\mf{d},\bs{\tau},\eps}|_{L^\frac{2N}{N+2}(\Omega_\eps)} \le C \left(\eps^\frac{N-2}2 + \eps \right) \le C \eps^\frac{N-2}2,
\end{equation}
where we used the fact that $N=3,4$, so that $\eps \le \eps^\frac{N-2}2$ for small $\eps$.

\smallskip

\noindent \textbf{Estimate of $|\tilde{N}^i_{\mf{d},\bs{\tau},\eps}(\bs{\phi})|_{L^\frac{2N}{N+2}(\Omega_\eps)}$}. At first, with the aid of Lemma \ref{lem a.3}, we compute
\begin{equation}\label{471}
\begin{split}
|\tilde{P}^i_{\mf{d},\bs{\tau},\eps}(\bs{\phi}) |_{\frac{2N}{N+2}}  &\le C | (P U_i)^{p-2} \phi_i^2 + |\phi_i|^p| _{\frac{2N}{N+2}}   \le C | U_i^{\frac{6-N}{N-2}} \phi_i^2|_{\frac{2N}{N+2}} + |\phi_i|_{2^*}^\frac{N+2}{N-2} \\
& \le C |U_i|_{2^*}^\frac{6-N}{N-2} |\phi_i|_{2^*}^2 +  |\phi_i|_{2^*}^\frac{N+2}{N-2} \le C \|\phi_i\|^2 +  \|\phi_i\|^\frac{N+2}{N-2}  
\end{split}
\end{equation}
Regarding the remaining term $\tilde{Q}^i_{\mf{d},\bs{\tau},\eps}(\bs{\phi})$, we focus at first on $N=3$. Using Lemma \ref{lem a.6} \footnote{in the Lemma, we consider $x=\mu_i^{-\frac1{p-1}} P U_i$, $y = \mu_j^{-\frac1{p-1}} P U_j$, $h_1=\phi_i$, $k_1=\phi_j$, $h_2=0$, $k_2=0$.} and the usual arguments we deduce that 
\begin{subequations}\label{sub 1}
\begin{equation}
\begin{split}
|\tilde{Q}^i_{\mf{d},\bs{\tau},\eps} (\bs{\phi})|_{\frac65} & \le C 
\left| |\mu_j^{-\frac1{p-1}} P U_j+ |\phi_j| |^3|\phi_i|^2\right.  \\
& \hphantom{ \le C | ++} + |\mu_i^{-\frac1{p-1}} P U_i+ |\phi_i| | \, |\mu_j^{-\frac1{p-1}} P U_j+|\phi_j| |^2 |\phi_i|\, |\phi_j| 
\\ & \hphantom{ \le C | ++} \left. +   |\mu_i^{-\frac1{p-1}} P U_i + |\phi_i||^2|\mu_j^{-\frac1{p-1}} P U_j+ |\phi_j| ||\phi_j|^2 \right|_\frac65
\\
&  \le C\left[ | U_j^3 \phi_i^2 |_\frac65 + | \phi_j^3 \phi_i^2|_\frac65 + | U_i U_j^2 \phi_i \phi_j |_\frac65 + |U_i \phi_i \phi_j^3 |_\frac65+ |U_j^2 \phi_i^2 \phi_j|_\frac65\right. \\
& \hphantom{ \le C | ++} \left.+ | U_i^2  U_j \phi_j^2|_\frac65 + |U_i^2 \phi_j^3|_\frac65 + |U_j \phi_i^2 \phi_j^2 |_\frac65\right] \\
& \le C \Big[|U_j|_{2^*}^3 |\phi_i|_{2^*}^2+|\phi_i|_{2^*}^2 |\phi_j|_{2^*}^3 + |U_i|_{2^*} |U_j|_{2^*}^2 |\phi_i|_{2^*}|\phi_j|_{2^*}  \\
&  \hphantom{ \le C | ++}  + |U_i|_{2^*}  |\phi_i|_{2^*}|\phi_j|_{2^*}^3 +  |U_j|_{2^*}^2 |\phi_i|_{2^*}^2|\phi_j|_{2^*}  +  |U_i|_{2^*}^2 |U_j|_{2^*} |\phi_j|_{2^*}^2 \\
&  \hphantom{ \le C | ++}  + |U_i|_{2^*}^2 |\phi_j|_{2^*}^3 +  |U_j|_{2^*} |\phi_i|_{2^*}^2|\phi_j|_{2^*}^2\Big] \\
& \le C \left( \|\bs{\phi}\|^2 + \|\bs{\phi}\|^5\right) = C  \left( \|\bs{\phi}\|^2 + \|\bs{\phi}\|^\frac{N+2}{N-2}\right).
\end{split}
\end{equation}
Let us consider now the easier case $N=4$. We have
\begin{multline*}
|\mu_j^{-\frac{1}{2}} P U_j  + \phi_j|^{2} (\mu_i^{-\frac{1}{2}} P U_i + \phi_i)    -  (\mu_j^{-\frac{1}{2}} P U_j)^{2}(\mu_i^{-\frac{1}{2}} P U_i) -   (\mu_j^{-\frac{1}{2}} P U_{\d_j,\xi_j})^{2}    \phi_i \\
 - 2 (\mu_j^{-\frac{1}{2}} P U_{\d_j,\xi_j})(\mu_i^{-\frac{1}{2}} P U_i)  \phi_j = \mu_i^{-\frac12} P U_i\, \phi_j^2  + 2 \mu_j^{-\frac12} P U_j \,\phi_i \phi_j + \phi_i \phi_j^2
\end{multline*}
Therefore
\begin{equation}\label{472}
\begin{split}
|\tilde{Q}^i_{\mf{d},\bs{\tau},\eps}(\bs{\phi}) |_{{\frac43}} 
& \le C\sum_{i \neq j} \left(  |U_i\, \phi_j^2|_{\frac43}+ |U_j \phi_i \phi_j|_{\frac43}+|\phi_i\phi_j^2|_{\frac43}\right)\\
& \le C \sum_{i \neq j} \left(|U_i|_{2^*} |\phi_j|_{2^*}^2 
 + |U_j|_{2^*} |\phi_i|_{2^*} |\phi_j|_{2^*} + |\phi_i|_{2^*} |\phi_j|_{2^*}^2   \right) \\
& \le C \left(\|\bs{\phi}\|^2 + \|\bs{\phi}\|^3 \right) = C \left(\|\bs{\phi}\|^2 + \|\bs{\phi}\|^\frac{N+2}{N-2} \right).
\end{split}
\end{equation}
\end{subequations}
Collecting \eqref{471} and \eqref{sub 1}, we conclude that for every $i=1,\dots,m$, and provided that $\|\bs{\phi}\|<1$,
\begin{align*}
\|\tilde{N}^i_{\mf{d},\bs{\tau},\eps}(\bs{\phi})\|_\frac{2N}{N+2} & \le C \left(\|\bs{\phi}\|^2 + \|\bs{\phi}\|^\frac{N+2}{N-2} \right) \le C \|\bs{\phi}\|^2
\end{align*}
for some $C>0$ depending only on the data. 

\medskip

Now we pass to the proof of the validity of \eqref{674}

We apply again Lemma \ref{lem a.6} (notice that if $N=4$ then $p=3$, while if $N=3$ then $p=5$): thus
\begin{equation}\label{672}
\begin{split}
 |&  \tilde Q^i_{\mf{d},\bs{\tau},\eps}(\bs{\phi}^1) - \tilde Q^i_{\mf{d},\bs{\tau},\eps}(\bs{\phi}^2)|_{\frac{2N}{N+2}} \\
& \le C  \frac{p-3}2\sum_{j \neq i} \left| \big| U_j + |\phi^1_j| + |\phi^2_j| \big|^\frac{p+1}2 \big| |\phi^1_i| + |\phi^2_i| \big| \big| \phi^1_i - \phi_i^2\big| \right|_{\frac{2N}{N+2}} \\
& +  C\sum_{j \neq i} \left| \big| U_i + |\phi^1_i| + |\phi^2_i| \big|^\frac{p-3}2 \big| U_j + |\phi^1_j| + |\phi^2_j| \big|^\frac{p-1}2 \big| |\phi^1_j| + |\phi^2_j| \big| \big| \phi^1_i - \phi_i^2\big| \right|_{\frac{2N}{N+2}} \\
& +  C\sum_{j \neq i} \left| \big| U_i + |\phi^1_i| + |\phi^2_i| \big|^\frac{p-3}2 \big| U_j + |\phi^1_j| + |\phi^2_j| \big|^\frac{p-1}2 \big| |\phi^1_i| + |\phi^2_i| \big| \big| \phi^1_j - \phi_j^2\big| \right|_{\frac{2N}{N+2}} \\
& + C\sum_{j \neq i} \left| \big| U_i + |\phi^1_i| + |\phi^2_i| \big|^\frac{p-1}2 \big| U_j + |\phi^1_j| + |\phi^2_j| \big|^\frac{p-3}2 \big| |\phi^1_j| + |\phi^2_j| \big| \big| \phi^1_j - \phi_j^2\big| \right|_{\frac{2N}{N+2}},
\end{split}
\end{equation}
where we wrote $U_i$ for $U_i$. To estimate the right hand side, we apply the H\"older and the Sobolev inequalities. We notice at first that the first term survives only if $N=3$. In such case, for any $\bs{\phi}^1, \bs{\phi}^2 \in Y_\eps$ (recall that then $\|\bs{\phi}^1\|, \|\bs{\phi}^2\| \le C \eps^\frac{N-2}2$)
\begin{align*}
\Big| \big| U_j + |\phi^1_j| & + |\phi^2_j| \big|^3 \big| |\phi^1_i| + |\phi^2_i| \big| \big| \phi^1_i - \phi_i^2\big| \Big|_{\frac{6}{5}}  \le  \big| U_j + |\phi^1_j| + |\phi^2_j| \big|_{2^*}^3 \, \big| |\phi^1_i| + |\phi^2_i| \big|_{2^*} \, | \phi^1_i - \phi_i^2|_{2^*} \\
& \le C \left( \|\bs \phi^1\| + \|\bs{\phi}^2\| \right) \| \phi^1_i - \phi_i^2\| \le C \eps^\frac{N-2}2  \| \phi^1_i - \phi_i^2\|.
\end{align*}
Similarly
\begin{align*}
\Big| \big| U_i & + |\phi^1_i| + |\phi^2_i| \big|^\frac{p-3}2 \big| U_j + |\phi^1_j| + |\phi^2_j| \big|^\frac{p-1}2 \big| |\phi^1_j| + |\phi^2_j| \big| \big| \phi^1_i - \phi_i^2\big| \Big|_{\frac{2N}{N+2}} \\
& \le \left| U_i + |\phi^1_i| + |\phi^2_i| \right|_{2^*}^\frac{p-3}2 \left| U_j + |\phi^1_j| + |\phi^2_j| \right|_{2^*}^\frac{p-1}2 \big| |\phi^1_j| + |\phi^2_j| \big|_{2^*} \, | \phi^1_i - \phi_i^2|_{2^*} \\
& \le C \eps^\frac{N-2}2  \| \phi^1_i - \phi_i^2\|,
\end{align*}
and estimating the remaining terms on the right hand side in \eqref{672} in the very same way (the terms are obtained one from the other after a permutation of the exponents), we deduce that
\[
 |  \tilde Q^i_{\mf{d},\bs{\tau},\eps}(\bs{\phi}^1) - \tilde Q^i_{\mf{d},\bs{\tau},\eps}(\bs{\phi}^2)|_{\frac{2N}{N+2}}  \le C \eps^\frac{N-2}2 \|\bs{\phi}^1-\bs{\phi}^2\|
\]
for any $\bs{\phi}^1,\bs{\phi}^2 \in Y_\eps$, as desired.

\section{}\label{app derivative error}

In this appendix we prove the validity of estimate \eqref{stima derivative error}. 

In this proof we write $U_i$ for $U_{\d_i,\x_i}$, $\psi_i^h$ for $\psi_{\d_i,\x_i}^h$, and $P$ for $P_\eps$ in order to keep the notation as short as possible.

\noindent \textbf{Step 1)} We consider $\pa_{s_{k,h}} L^i_{\mf{d},\bs{\tau},\eps}(\bs{\phi})$, with $\bs{\phi} \in \mf{K}_{\mf{d},\bs{\tau},\eps}^\perp$. Let us focus at first on the case $k=i$. We have
\begin{align*}
\pa_{s_{i,h}} & L^i_{\mf{d},\bs{\tau},\eps}(\bs{\phi}) = \Pi_i^\perp \circ i^* \Bigg[ p (p-1) (P U_i)^{p-2} (P \psi_i^h) \phi_i + \\
& + \frac{(p-1)(p-3)}4  \sum_{j \neq i} \beta_{ij} (\mu_j^{-\frac1{p-1}} P U_j)^\frac{p+1}2 (\mu_i^{-\frac1{p-1}} P U_i)^\frac{p-5}2 \mu_i^{-\frac1{p-1}} (P \psi_i^h) \phi_i \\
& + \frac{(p+1)(p-1)}4 \sum_{j \neq i} \beta_{ij} (\mu_j^{-\frac1{p-1}} P U_j)^\frac{p-1}2 (\mu_i^{-\frac1{p-1}} P U_i)^\frac{p-3}2 \mu_i^{-\frac1{p-1}} (P \psi_i^h) \phi_j \Bigg].
\end{align*} 
Therefore, using the continuity of $\Pi^\perp$ and $i^*$, and the fact that $0 \le P U_i \le U_i$, we deduce that
\begin{align*}
\|\pa_{s_{i,h}}  L^i_{\mf{d},\bs{\tau},\eps}(\bs{\phi}) \| &\le C | U_i^{p-2} (P \psi_i^h) \phi_i|_{\frac{2N}{N+2}} + C \sum_{j \neq i} | U_i^{\frac{p-3}2} U_j^{\frac{p-1}2} (P \psi_i^h) \phi_j|_{\frac{2N}{N+2}}\\
& \qquad + C (p-3) \sum_{j \neq i} | U_i^{\frac{p-5}2} U_j^{\frac{p+1}2} (P \psi_i^h) \phi_i|_{\frac{2N}{N+2}} \\
& \le C \underbrace{| U_i^{p-2} \psi_i^h \phi_i|_{\frac{2N}{N+2}}}_{=:(I)} + C \sum_{j \neq i} \underbrace{| U_i^{\frac{p-3}2} U_j^{\frac{p-1}2}  \psi_i^h \phi_j|_{\frac{2N}{N+2}}}_{=:(II)}
\\
&\qquad + C (p-3) \sum_{j \neq i} \underbrace{| U_i^{\frac{p-5}2} U_j^{\frac{p+1}2} \psi_i^h \phi_i|_{\frac{2N}{N+2}}}_{=:(III)} + n.t.
\end{align*}
where $n.t$ stays for ``negligible terms", and denotes a quantity which is a small $o$ of the previous ones as $\eps \to 0$. The exact shape of these terms can be computed using Corollary \ref{corol a.2} (similarly as we have already done in Appendix \ref{app step 1 linear}), but it's not so important for what follows and hence it is omitted for the sake of brevity. Notice also that term ($III$) is present only in case $N=3$, in which case $p=5$.

Now, by the H\"older and the Sobolev inequalities, and using Lemmas \ref{lem a.4}-\ref{lem a.9}, we compute for $h=0$
\begin{equation}
\begin{split}
|(I)|  \le  C| U_i^{p-2} \psi_i^h |_{\frac{N}{2}} \|\phi_i\| \le C \d_i^{-1} \left( \int_{\Omega_\eps} U_i^{2^*} \right)^{\frac2N}  \|\phi_i\| \le C \d_i^{-1} \|\phi_i\|,
\end{split}
\end{equation}
\begin{equation}
\begin{split}
|(II)| &\le C| U_i^\frac{p-3}2 U_j^\frac{p-1}{2} \psi_i^h |_{\frac{N}{2}} \|\phi_j\| \le C \d_i^{-1} \left( \int_{\Omega_\eps} U_i^\frac{N}{N-2} U_j^\frac{N}{N-2} \right)^{\frac2N}  \|\phi_j\| \\
& \le C \d_j ( |\log \d_i| + |\log \d_j| )^{\frac2N}\|\phi_j\|,
\end{split}
\end{equation}
and in case $N=3$
\begin{equation}
\begin{split}
|(III)| & \le C \|\phi_i\| |U_j^3 \psi_i^h|_{\frac32} \le C \d_i^{-1} \|\phi_i\| \left( \int_{\Omega_\eps} U_i^{\frac32} U_j^\frac92\right)^{\frac23} \\
& \le C \d_i^{-1}  (\d_i \d_j)^{\frac12} \|\phi_i\|.
\end{split}
\end{equation}
Moreover, for $h=1,\dots,N$,
\begin{equation}
\begin{split}
|(I)| & \le C| U_i^{p-2} \psi_i^h |_{\frac{N}{2}} \|\phi_i\| \le C \d_i^{-1} \left( \int_{\Omega_\eps} U_i^{\frac{3N}{N-2}} |x_h-\xi_{i,h}|^\frac{N}2 \right)^{\frac2N}  \|\phi_i\| \le C \d_i^{-1} \|\phi_i\|,
\end{split}
\end{equation}
\begin{equation}
\begin{split}
|(II)| & \le C \|\phi_j\| | U_i^\frac{p-3}2 U_j^\frac{p-1}{2} \psi_i^h |_{\frac{N}{2}} \le C \d_i^{-1}\|\phi_j\|  \left( \int_{\Omega_\eps} U_i^\frac{2N}{N-2} U_j^\frac{N}{N-2} |x_h-\xi_{i,h}|^\frac{N}2  \right)^\frac2{N} \\
& \le C \d_i^{-1}\|\phi_j\| (\d_i \d_j),
\end{split}
\end{equation}
and for $N=3$
\begin{equation}
\begin{split}
|(III)| & \le C \|\phi_i\| |U_j^3 \psi_i^h|_{\frac32} \le C \d_i^{-1} \|\phi_i\| \left( \int_{\Omega_\eps} U_j^\frac92 U_i^\frac92 |x_h-\xi_{i,h}|^\frac32 \right)^\frac23 \le C \d_i^\frac12 \d_j^\frac12 \|\phi_i\|.
\end{split}
\end{equation}
Altogether, recalling the ansatz \eqref{asympt expansion} and that $\|\bs{\phi}\| \le C \eps^\frac{N-2}2$ by Lemma \ref{lem: exist nonlin}, we conclude that for every $i=1,\dots,m$ and $h=0,\dots,N$
\[
\|\pa_{s_{i,h}}  L^i_{\mf{d},\bs{\tau},\eps}(\bs{\phi}) \|_{H_0^1(\Omega_\eps)} \le C \eps^\frac{N-3}2,
\]
as desired. 

Now, if $k \neq i$, we have still to consider
\begin{align*}
\pa_{s_{k,h}} &  L^i_{\mf{d},\bs{\tau},\eps}(\bs{\phi}) \\
& = \Pi_i^\perp \circ i^* \Bigg[\frac{(p-1)(p+1)}4 \beta_{ki} (\mu_k^{-\frac1{p-1}} P U_k)^\frac{p-1}2 (\mu_i^{-\frac1{p-1}} P U_i)^\frac{p-3}2 \mu_k^{-\frac1{p-1}} P \psi_k^h \, \phi_i \\
&  + \frac{(p-1)(p+1)}4 \beta_{ki} (\mu_k^{-\frac1{p-1}} P U_k)^\frac{p-3}2 (\mu_i^{-\frac1{p-1}} P U_i)^\frac{p-1}2 \mu_k^{-\frac1{p-1}} P \psi_k^h \, \phi_k \Bigg].
 \end{align*}
This term can be treated exactly as the one with $k=i$ (actually the computations are easier).

\smallskip

\noindent \textbf{Step 2)} We consider $\pa_{s_{k,h}} R^i_{\mf{d},\bs{\tau},\eps}(\bs{\phi})$, and again we present the details only for the (a bit harder) case $k=i$. We have
\begin{align*}
\pa_{s_{i,h}} &R^i_{\mf{d},\bs{\tau},\eps} = \Pi_i^\perp \circ i^* \Bigg[ \mu_i^{-\frac1{p-1}} p \left((P U_i)^{p-1}  P \psi_i^h - U_i^{p-1} \psi_i^h\right) \\
& + \frac{p-1}2\sum_{j \neq i} \beta_{ij} (\mu_j^{-\frac1{p-1}} P U_j)^\frac{p+1}2  (\mu_i^{-\frac1{p-1}} P U_i)^\frac{p-3}2 (\mu_i^{-\frac1{p-1}} P \psi_i^h)  \Bigg], 
\end{align*}
and hence
\begin{align*}
\|\pa_{s_{i,h}} &R^i_{\mf{d},\bs{\tau},\eps}\| \le C  |(P U_i)^{p-1} P \psi_i^h - U_i^{p-1} \psi_i^h|_{\frac{2N}{N+2}}  + C\sum_{j \neq i} |U_j^{\frac{p+1}2} U_i^\frac{p-3}2 P \psi_i^h|_{\frac{2N}{N+2}} \\
& \le  C| P U_i^{p-1} | P \psi_i^h - \psi_i^h | |_{\frac{2N}{N+2}}+  C | |P U_i^{p-1} -U_i^{p-1}| \psi_i^h|_\frac{2N}{N+2}  + C\sum_{j \neq i} |U_j^{\frac{p+1}2} U_i^\frac{p-3}2  \psi_i^h|_{\frac{2N}{N+2}} + n.t. \\
& \le C \underbrace{| U_i^{p-1} | P \psi_i^h - \psi_i^h | |_{\frac{2N}{N+2}}}_{=:(I)} + C  \underbrace{| \psi_i^h U_i^{p-2} |P U_i -U_i||_\frac{2N}{N+2}}_{=:(II)}  + C\sum_{j \neq i}  \underbrace{|U_j^{\frac{p+1}2} U_i^\frac{p-3}2  \psi_i^h|_{\frac{2N}{N+2}}}_{=:(III)} + n.t.,
\end{align*}
where we used Lemma \ref{lem a.3} and the fact that $0 \le P U_i \le U_i$. 

As in the first step, using the lemmas collected in Appendix \ref{app tools}, we estimate separately the three terms on the right hand side, starting from the case $h=0$:
\begin{equation}
\begin{split}
|(I)| & \le C\left|  U_i^{\frac{4}{N-2}} \left( \d_i^{\frac{N-4}2} + \frac{\d_i^{\frac32N-4}}{|x-a_i|^{N-2}} \right) \right|_{\frac{2N}{N+2}} \\
& \le C \d_i^\frac{N-4}2 \left( \int_{\Omega_\eps} U_i^\frac{8N}{(N+2)(N-2)} \right)^{\frac{N+2}{2N}} + C \d_i^{\frac32N-4} \left( \int_{\Omega_\eps} \frac{U_i^\frac{8N}{(N+2)(N-2)}}{|x-a_i|^{\frac{2N(N-2)}{N+2}}} \right)^{\frac{N+2}{2N}} \\
& \le C \d_i^{\frac{N-4}2 + \frac{N-2}2} + C \d_i^{\frac32N-4 -\frac{N-2}2} \le C \d_i^{N-3},
\end{split}
\end{equation}
\begin{equation}
\begin{split}
|(II)| & \le C \d_i^{-1}\left| U_i^{p-1} \left( \d_i^\frac{N-2}2 + \frac{\d_i^{\frac{3N}2-3}  }{|x-a_i|^{N-2}} \right)\right|_{\frac{2N}{N+2}} \\
& \le C \d_i^\frac{N-4}2 \left( \int_{\Omega_\eps} U_i^\frac{8N}{(N+2)(N-2)} \right)^{\frac{N+2}{2N}} + C \d_i^{\frac32N-4} \left( \int_{\Omega_\eps} \frac{U_i^\frac{8N}{(N+2)(N-2)}}{|x-a_i|^{\frac{2N(N-2)}{N+2}}} \right)^{\frac{N+2}{2N}}  \le C \d_i^{N-3},
\end{split}
\end{equation}
\begin{equation}
\begin{split}
|(III)| & \le C \d_i^{-1} |U_j^{\frac{p+1}2} U_i^{\frac{p-1}2}|_{\frac{2N}{N+2}} \le C \d_i^{-1} \left( \d_i^{\frac{2N}{N+2}} \d_j^{\frac{2N}{N+2}} \right)^{\frac{N+2}{2N}} \le C \d_j.
\end{split}
\end{equation}
Analogously, for $h =1,\dots,N$, we have
\begin{equation}
\begin{split}
|(I)| & \le C\left|  U_i^{\frac{4}{N-2}} \left( \d_i^{\frac{N-2}2} + \frac{\d_i^{\frac32N-4}}{|x-a_i|^{N-2}} \right) \right|_{\frac{2N}{N+2}} \le C \d_i^{\frac{N-2}2 + \frac{N-2}2} + C \d_i^{\frac32N-4 -\frac{N-2}2} \le C \d_i^{N-3},
\end{split}
\end{equation}
\begin{equation}
\begin{split}
|(II)| & \le C \d_i^{-1}\left| U_i^\frac6{N-2} |x_h-\xi_{i,h}|\left( \d_i^\frac{N-2}2 + \frac{\d_i^{\frac32N-3}  }{|x-a_i|^{N-2}}\right)\right|_{\frac{2N}{N+2}} \\
& \le C \d_i^{\frac{N-4}2} \left( \int_{\Omega_\eps} U_i^\frac{12N}{(N+2)(N-2)} |x_h-\xi_{i,h}|^{\frac{2N}{N+2}} \right)^{\frac{N+2}{2N}} \\
&+ C \d_i^{\frac32N-4} \left( \int_{\Omega_\eps} U_i^\frac{12N}{(N+2)(N-2)}  \frac{|x_h-\xi_{i,h}|^{\frac{2N}{N+2}}}{|x-a_i|^{\frac{2N(N-2)}{N+2}}}\right)^{\frac{N+2}{2N}} \\
& \le C  \d_i^{\frac{N-4}2 + \frac{N-2}2} + C \d_i^{\frac32N-4 + \frac{2-N}2} \le C \d_i^{N-3},
\end{split}
\end{equation}
\begin{equation}
\begin{split}
|(III)| & \le C \d_i^{-1} |U_j^{\frac{N}{N-2}} U_i^{\frac{4}{N-2}} |x_h-\xi_{i,h}| |_{\frac{2N}{N+2}} \\
& \le C \d_i^{-1} \left( \int_{\Omega_\eps} U_i^{\frac{8N}{(N+2)(N-2)}} U_j^{\frac{2N^2}{(N+2)(N-2)}} |x_h-\xi_{i,h}|^\frac{2N}{N+2}  \right)^\frac{N+2}{2N}\\
& \le C \d_i^{-1} \Bigg( \d_i^{\frac{4N}{N+2}} \d_j^\frac{N^2}{N+2} + \d_i^{\frac{4N}{N+2}} \int_{\Omega_\eps} U_j^{\frac{2N^2}{(N+2)(N-2)}} 
%\\
%& \qquad \qquad 
+
 \d_j^\frac{N^2}{N+2} \int_{\Omega_\eps} U_i^{\frac{8N}{(N+2)(N-2)}} |x_h-\xi_{i,h}|^\frac{2N}{N+2} \Bigg)^\frac{N+2}{2N} \\
& \le C \d_i^{-1} \left( \d_i^2 \d_j + \d_i^2 \d_j^\frac{N}2 |\log \d_j|^{\frac{N+2}{2N} } \right) = o(\eps^\frac{N-3}2).
\end{split}
\end{equation}
Notice that the logarithm appears in the computations only for $N=4$, see Lemma \ref{lem a.11}. 

Collecting together the previous estimates, we deduce that $\| \pa_{s_{i,h}} R^i_{\mf{d},\bs{\tau},\eps}\| \le C \eps^\frac{N-3}2$, as desired.

\smallskip

\noindent \textbf{Step 3)} We consider $\pa_{s_{k,h}} N^i_{\mf{d},\bs{\tau},\eps}(\bs{\phi})$, with $\bs{\phi} \in \mf{K}_{\mf{d},\bs{\tau},\eps}^\perp$. Once again, we focus on the case $k=i$, which presents the heaviest calculations. We have
\begin{align*}
\pa_{s_{k,h}} & N^i_{\mf{d},\bs{\tau},\eps}(\bs{\phi}) \\
& = \Pi_i^\perp \circ i^* \Bigg[ \left( f'( \mu_i^{-\frac{1}{p-1}} P U_i + \phi_i) -  f'( \mu_i^{-\frac{1}{p-1}} P U_i)   - f''( \mu_i^{-\frac{1}{p-1}} P U_i ) \phi_i \right) \mu_i^\frac{p-2}{p-1} P \psi_i^h \\
& + \frac{p-1}2\sum_{j \neq i} \beta_{ij} |\mu_j^{-\frac{1}{p-1}} P U_j + \phi_j|^{\frac{p+1}{2}} |\mu_i^{-\frac{1}{p-1}} P U_i + \phi_i|^{\frac{p-3}{2}} \mu_i^{-\frac{1}{p-1}} P \psi_i^h  \\
& - \frac{p-1}2 \sum_{j \neq i} \beta_{ij} (\mu_j^{-\frac{1}{p-1}} P U_j)^{\frac{p+1}{2}}(\mu_i^{-\frac{1}{p-1}} P  U_i)^{\frac{p-3}{2}} \mu_i^{-\frac{1}{p-1}} P \psi_i^h \\
& - \frac{(p-1)(p-3)}{4} \sum_{j \neq i} \beta_{ij} (\mu_j^{-\frac{1}{p-1}} P U_j)^{\frac{p+1}{2}}   
(\mu_i^{-\frac{1}{p-1}} P U_i)^{\frac{p-5}{2}} \mu_i^{-\frac{1}{p-1}} P \psi_i^h \,  \phi_i \\
& - \frac{(p+1)(p-1)}{4} \sum_{j \neq i} \beta_{ij} (\mu_j^{-\frac{1}{p-1}} P U_j)^{\frac{p-1}{2}}   
(\mu_i^{-\frac{1}{p-1}} P U_i)^{\frac{p-3}{2}} \mu_i^{-\frac{1}{p-1}} P\psi_i^h \, \phi_j \Bigg].
\end{align*}
Therefore, if $N=3$ we can use the Lagrange theorem to deduce that \begin{align*}
\|\pa_{s_{k,h}} & N^i_{\mf{d},\bs{\tau},\eps}(\bs{\phi}) \|  \le C \big| (U_i^2  + \phi_i^2) \phi_i^2 P \psi_i^h\big|_{\frac{6}{5}}  
 \\
& + C \sum_{j \neq i} \bigg( \big| |U_i + |\phi_i|| |U_j + |\phi_j||^2 \phi_j \,P \psi_i^h \big|_\frac65 +  \big| |U_j+ |\phi_j|| \phi_i \,P \psi_i^h \big|_\frac65 \\
& \hphantom{+ C \sum_{j \neq i} \bigg(  } \ + | U_j^3 \phi_i \, P \psi_i^h|_{\frac{6}{5}}  +  |U_i U_j^2 \phi_j \, P \psi_i^h|_{\frac{6}{5}} \bigg); 
%
%
%
%\sum_{j \neq i} |U_i U_j^\frac{p+1}2 P \psi_i^h|_{\frac{2N}{N+2}} \\
%&+  C \sum_{j \neq i} |U_i^\frac{p-3}{2} \phi_j^\frac{p+1}2 P \psi_i^h|_{\frac{2N}{N+2}}   + C \sum_{j \neq i} |\phi_i^\frac{p-3}{2} U_j^\frac{p+1}2 P \psi_i^h|_{\frac{2N}{N+2}}  + C \sum_{j \neq i} |\phi_i^\frac{p-3}{2} \phi_j^\frac{p+1}2 P \psi_i^h|_{\frac{2N}{N+2}} \\
%& + C \frac{p-3}2 \sum_{j \neq i} |U_i^\frac{p-5}2 U_j^\frac{p+1}2 \phi_i P \psi_i^h|_\frac{2N}{N+2}  + C \sum_{j \neq i} |U_i^\frac{p-1}2 U_j^\frac{p+1}2 \phi_j P \psi_i^h|_\frac{2N}{N+2} \\
%& = C |f'''( \mu_i^{-\frac{1}{p-1}} P U_i + t \phi_i) \phi_i^2 \psi_i^h|_{\frac{2N}{N+2}}  + C \sum_{j \neq i} |U_i^\frac{p-3}{2} U_j^\frac{p+1}2  \psi_i^h|_{\frac{2N}{N+2}} \\
%&+  C \sum_{j \neq i} |U_i^\frac{p-3}{2} \phi_j^\frac{p+1}2  \psi_i^h|_{\frac{2N}{N+2}}   + C \sum_{j \neq i} |\phi_i^\frac{p-3}{2} U_j^\frac{p+1}2  \psi_i^h|_{\frac{2N}{N+2}}  + C \sum_{j \neq i} |\phi_i^\frac{p-3}{2} \phi_j^\frac{p+1}2  \psi_i^h|_{\frac{2N}{N+2}} \\
%& + C \frac{p-3}2 \sum_{j \neq i} |U_i^\frac{p-5}2 U_j^\frac{p+1}2 \phi_i  \psi_i^h|_\frac{2N}{N+2}  + C  \sum_{j \neq i} |U_i^\frac{p-1}2 U_j^\frac{p+1}2 \phi_j  \psi_i^h|_\frac{2N}{N+2} + n.t.
\end{align*}
if on the other hand $N=4$, then by direct computations
\begin{align*}
\|\pa_{s_{k,h}} & N^i_{\mf{d},\bs{\tau},\eps}(\bs{\phi}) \|  \le C |\phi_i^2 P \psi_i^h|_\frac43 \\
& + C \sum_{j \neq i} \bigg(  | U_j \, \phi_j\, P \psi_i^h|_\frac43+ | \phi_j^2 P \psi_i^h|_\frac43  + | U_j \, \phi_j \, P \psi_i^h |_\frac43 \bigg).
\end{align*}
In both cases, using the lemmas in Appendix \ref{app tools} as in the previous steps, it is not difficult to check that 
\[
\|\pa_{s_{k,h}}  N^i_{\mf{d},\bs{\tau},\eps}(\bs{\phi}) \| \le C \eps^{\frac{N-3}2}.
\]
This completes the proof of \eqref{stima derivative error}.

%\bibliography{bibliography}
%\bibliographystyle{abbrv}

\end{document}